\documentclass[12pt]{amsart}

\usepackage[leqno]{amsmath}
\usepackage{amssymb}
\usepackage{amsxtra}
\usepackage{amscd}
\usepackage{hyperref}
\usepackage{graphicx, color}
\usepackage{dcpic}

\newtheorem{thm}{Theorem}
\newtheorem{prop}[thm]{Proposition}

\newtheorem{defn}[thm]{Definition}

\newtheorem{lemma}[thm]{Lemma}
\newtheorem{cor}[thm]{Corollary}

\newcommand{\R}{\mathbb{R}}

\newcommand{\F}{\mathbb{F}}
\newcommand{\I}{\ensuremath{\mathbb{I}}}
\newcommand{\N}{\mathbb{N}}
\newcommand{\half}{\mathbb{H}}

\newcommand{\Z}{\mathbb{Z}}
\newcommand{\Zmod}[1]{\Z/{#1}\Z}

\newcommand{\disp}[1]{\displaystyle{#1}}

\newcommand{\graph}[1]{\Gamma_{L}}
\newcommand{\circles}[1]{\ensuremath{\mathrm{\textsc{cir}}(#1)}}

\newcommand{\lra}{\longrightarrow}
\newcommand{\inlinediag}[2][0.33]{\includegraphics[scale=#1]{./images/#2}}

\newcommand{\tangle}[1]{\mathcal{#1}}

\newcommand{\cross}[1]{\ensuremath{\mathrm{\textsc{cr}}(#1)}}
\newcommand{\merge}[1]{\ensuremath{\mathrm{\textsc{merge}}(#1)}}
\newcommand{\leftMerge}[1]{\ensuremath{\overleftarrow{\mathrm{\textsc{merge}}}(#1)}}
\newcommand{\rightMerge}[1]{\ensuremath{\overrightarrow{\mathrm{\textsc{merge}}}(#1)}}
\newcommand{\fission}[1]{\ensuremath{\mathrm{\textsc{divide}}(#1)}}
\newcommand{\leftFission}[1]{\ensuremath{\overleftarrow{\mathrm{\textsc{divide}}}(#1)}}
\newcommand{\rightFission}[1]{\ensuremath{\overrightarrow{\mathrm{\textsc{divide}}}(#1)}}

\newcommand{\dec}[1]{\ensuremath{\mathrm{\textsc{dec}}(#1)}}
\newcommand{\interior}[1]{\ensuremath{\mathrm{\textsc{interior}}(#1)}}
\newcommand{\free}[1]{\ensuremath{\mathrm{\textsc{free}}(#1)}}

\newcommand{\state}[1]{\ensuremath{\mathrm{\textsc{State}}(#1)}}
\newcommand{\actor}[1]{\ensuremath{\mathrm{\textsc{active}}(#1)}}

\newcommand{\resolution}[1]{\ensuremath{\mathrm{\textsc{res}}(#1)}}

\newcommand{\lefty}[1]{\overleftarrow{#1}}
\newcommand{\righty}[1]{\overrightarrow{#1}}

\newcommand{\leftHalf}{\overleftarrow{\half}}
\newcommand{\rightHalf}{\overrightarrow{\half}}

\newcommand{\bridgeGraph}[1]{\Gamma_{#1}}

\newcommand{\cleaved}[1]{\mathcal{C\!L}_{#1}}
\newcommand{\cleave}[1]{\widehat{\mathcal{C\!L}}_{#1}}
\newcommand{\leftBridges}[1]{\ensuremath{\overleftarrow{\mathrm{\textsc{Br}}}(#1)}}
\newcommand{\rightBridges}[1]{\ensuremath{\overrightarrow{\mathrm{\textsc{Br}}}(#1)}}
\newcommand{\bridges}[1]{\ensuremath{\mathrm{\textsc{Bridge}}(#1)}}

\newcommand{\startCircle}[1]{C_{a}(#1)}
\newcommand{\terminalCircle}[1]{C_{b}(#1)}

\newcommand{\leftnorm}[1]{\lefty{l}(#1)}
\newcommand{\leftComplex}[1]{\langle\!\!\!\langle\,#1\,]\!]} 
\newcommand{\rightComplex}[1]{[\![\, #1 \,\rangle\!\!\!\rangle}
\newcommand{\complex}[1]{\ensuremath{\langle\!\!\!\langle\,#1\,\rangle\!\!\!\rangle}}

\title{A type $A$ structure in Khovanov Homology}
\author{Lawrence P. Roberts}
\thanks{This research was supported by a research grant from Research Grants Committee of the University of Alabama, Tuscaloosa}

\begin{document}
\begin{abstract}    
Inspired by bordered Floer homology, we describe a type $A$ structure on a Khovanov homology for a tangle which complements the type $D$ structure previously defined by the author. The type $A$ structure is a differential module over a certain algebra. This can be paired with the type $D$ structure to recover the Khovanov chain complex. The homotopy type of the type $A$ structure is a tangle invariant, and homotopy equivalences of the type $A$ structure result in chain homotopy equivalences on the Khovanov chain complex. We can use this to simplify computations and introduce a modular approach to the computation of Khovanov homologies. This approach adds to the literature even in the case of a connect sum, where the techniques here will allow an exact computation of Khovanov homology from the stuctures for two tangles coming from the summands. Several examples are included, showing in particular how we can compute the correct torsion summands for the Khovanov homology of the connect sum. A lengthy appendix is devoted to establishing the theory of these structures over a characterstic zero ring. 
\end{abstract}
\maketitle

\section{Introduction}
\noindent In a previous paper, \cite{typeD}, we described an algebra $\mathcal{B}\Gamma_{n}$ for a set of $2n$ points $P_{2n}$ ordered along a line (summarized in the next section) and a type $D$ structure $\rightComplex{\righty{T}}$ for an outside tangle $\righty{T}$ whose endpoints are these $2n$ points. An outside tangle is one with a diagram in an oriented half-plane whose boundary contains $P_{2n}$ but provides $P_{2n}$ with the opposite linear ordering when inherited from the boundary orientation. In this paper, we consider inside tangles: tangles where the orientation on the boundary equips $P_{2n}$ with the same ordering. We will picture these as lying on the left side of the $y$-axis in $R^{2n}$. For example, the following is an inside tangle $\lefty{T}$ over $P_{4}$ when the plane has its usual orientation
$$
\inlinediag[0.5]{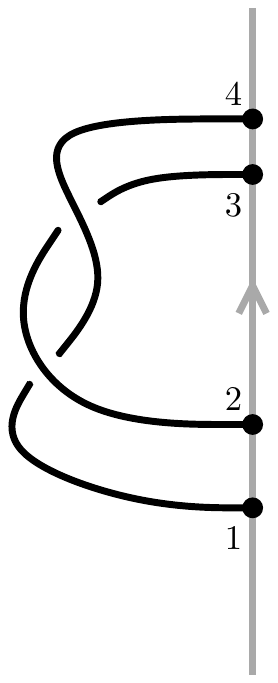}
$$
These tangles will be taken with an orientation, although we suppress that data for the introduction. To such a tangle we will associate a bigraded module $\leftComplex{\lefty{T}}$ and a differential $d_{APS}$, which is a modified version of the differential defined by M. Asaeda, J. Przytycki, and A. Sikora in \cite{APS} for their tangle homology. It is modified to have more generators, in a manner similar to Khovanov's invariant for tangles in \cite{Khta}. \\
\ \\
\noindent From there we define a bigrading preserving right action $\leftComplex{\righty{T}} \otimes \mathcal{B}\Gamma_{n} \rightarrow \leftComplex{\righty{T}}$ which is compatible with $d_{APS}$ by a certain Leibniz identity.  This will make $\leftComplex{\righty{T}}$ into a differential right module over $\mathcal{B}\Gamma_{n}$. If we consider this within a suitable category of right $A_{\infty}$-modules we have a notion of homotopy equivalence of right modules. We will then show that Reidemeister moves on the diagram $\lefty{T}$ will produce homotopy equivalent $A_{\infty}$-modules. We do this over $\Z$ with a somewhat different sign convention than usual, and a good bit of this paper is taken up by ensuring that the sign choices will work (the reader should consider that there are different sign conventions that can be followed in the Khovanov construction, and that these will produce distinct ``even'' and ``odd'' versions -- we only consider the original, ``even,'' version here). Following the conventions of bordered Floer homology, \cite{Bor1}, we will call this a type A structure. \\
\ \\
\noindent We arrange these constructions so that the following argument will work. Consider the following knot $K$ cleaved transversely in half by the $y$-axis:
$$
\inlinediag[0.5]{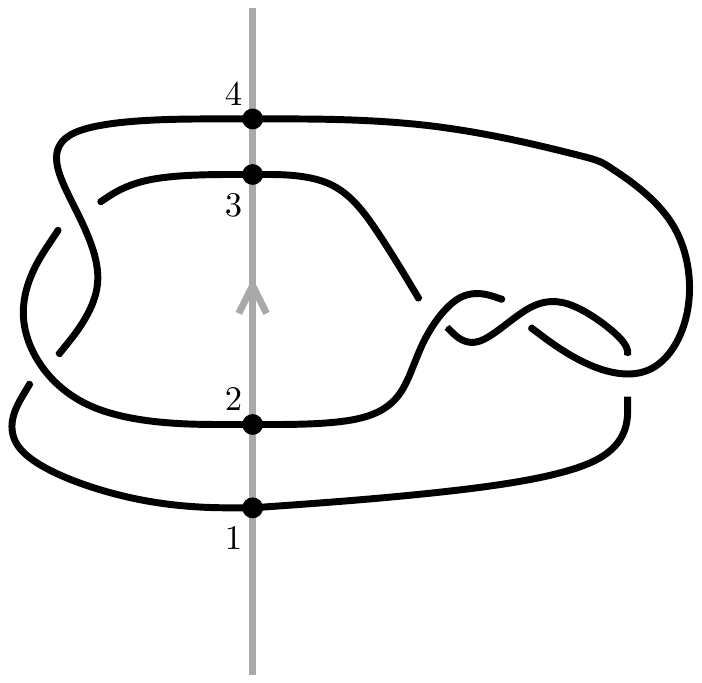}
$$
On the left side of the $y$-axis we recognize the inside tangle $\lefty{T}$. On the right side, there is an outside tangle $\righty{T}$:
$$
\inlinediag[0.5]{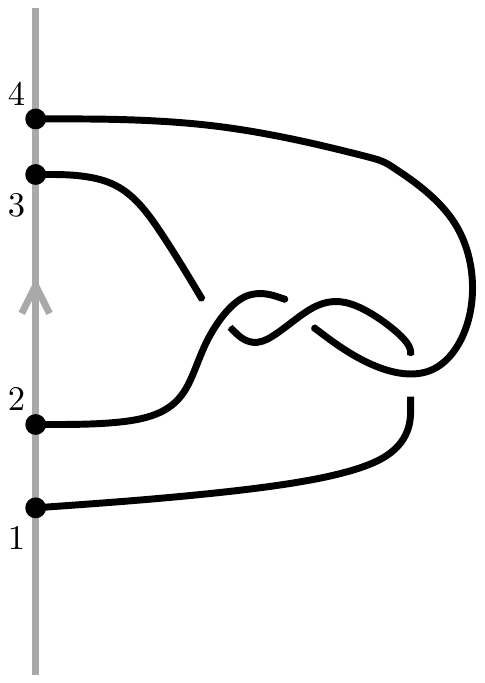}
$$
\noindent The Khovanov complex $\complex{K}$ is generated by states consisting of a smoothing at each crossing of a diagram for $K$  and a decoration of $\{+,-\}$ attached to each planar circle resulting from the smoothing. Such a state $\xi$ might look like
$$
\inlinediag[0.5]{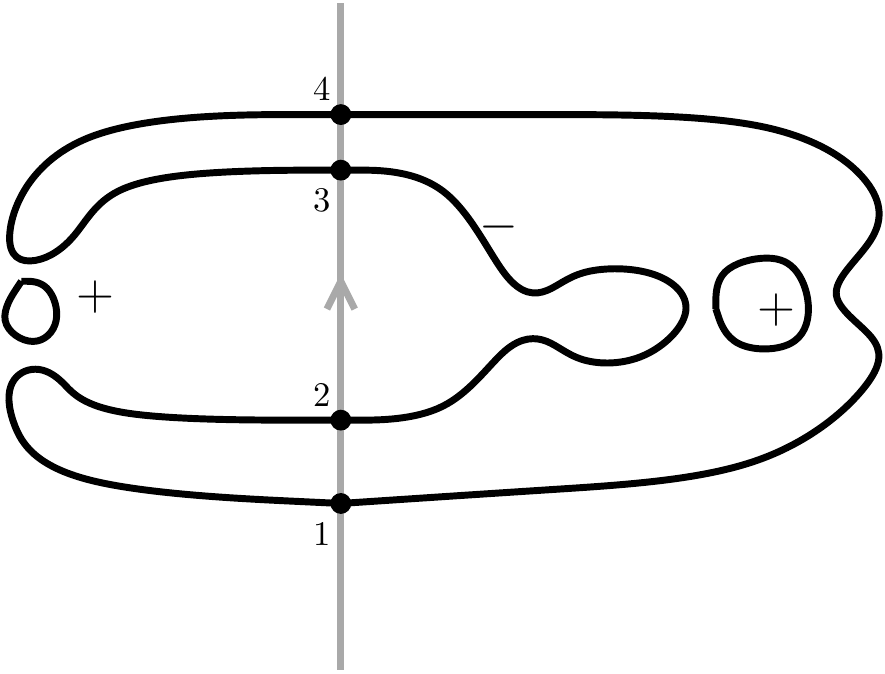}
$$
We can similarly divide this resolution along the $y$-axis; however, we do this in a less obvious way. The left side will be the diagram obtained by forgetting the circles on the right which do not intersect the $y$-axis. We similarly describe the right side:
$$
\inlinediag[0.5]{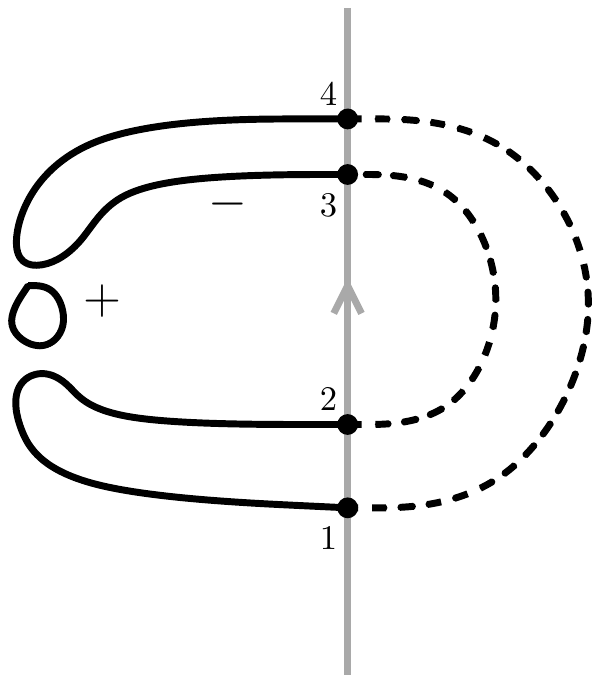} \hspace{1in} \inlinediag[0.5]{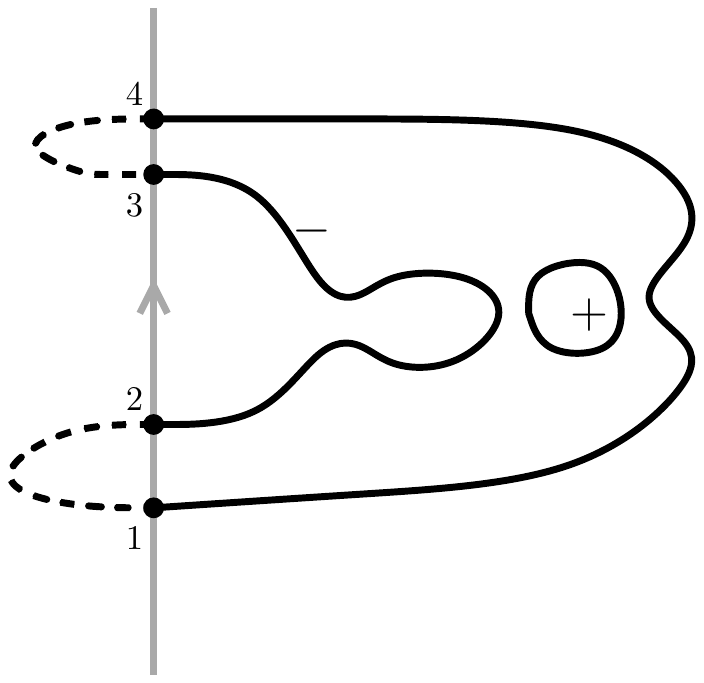}
$$
These are states, $\lefty{\xi}$ and $\righty{\xi}$, generating summands in $\leftComplex{\lefty{T}}$ and $\rightComplex{\righty{T}}$, respectively. To obtain the resolution of $K$ we will consider these to be glued along their common cleaved link:
$$
\inlinediag[0.5]{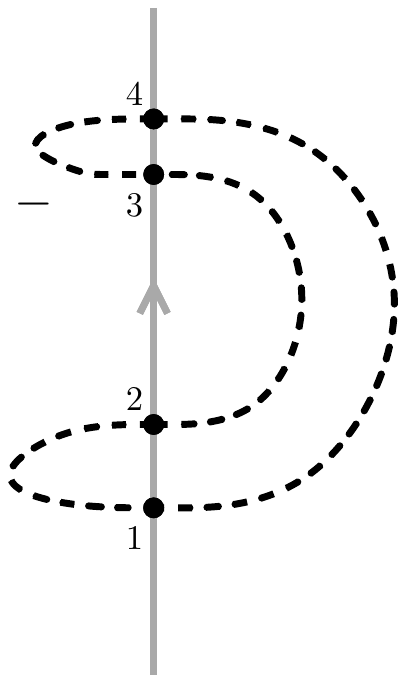}
$$
The latter diagram corresponds to an idempotent in $\mathcal{B}\Gamma_{n}$ which acts on the two states as the identity. These idempotents will be orthogonal in $\mathcal{B}\Gamma_{n}$, and if we let $\mathcal{I}$ be the idempotent subalgebra, then $\leftComplex{\lefty{T}} \otimes_{\mathcal{I}} \rightComplex{\righty{T}}$ will be isomorphic to $\complex{K}$, and $\lefty{\xi} \otimes \righty{\xi}$ will represent $\xi$ in this decomposition. The Khovanov differential can then be decomposed into the contribution of the crossings on the right and left. However, these contributions can change the cleaved link, and the corresponding idempotent. We record the changed in the cleaved link with the algebra $\mathcal{B}\Gamma_{n}$. For the crossings on the right we obtain a map $\righty{\delta}: \rightComplex{\righty{T}} \longrightarrow \mathcal{B}\Gamma_{n} \otimes_{\mathcal{I}} \rightComplex{\righty{T}}$ which satisfies the requirements of a type D structure, \cite{Bor1}. The crossings on the left give rise to the type $A$ structure. \\
\ \\
\noindent Following the constructions in \cite{Bor1} we can combine the type $A$ structure on $\leftComplex{\lefty{T}}$ and the type $D$ structure on $\rightComplex{\righty{T}}$ into a chain complex $\leftComplex{\lefty{T}} \boxtimes \rightComplex{\righty{T}}$ with underlying module $\leftComplex{\lefty{T}} \otimes_{\mathcal{I}} \rightComplex{\righty{T}}$ and differential
$$
\partial^{\boxtimes}(x \otimes y) =  d_{APS}(x) \otimes |y| + (m_{2} \otimes \I)(x \otimes \righty{\delta}(y))
$$ 
where $m_{2}$ is the action on $\leftComplex{\lefty{T}}$. We then show that $\complex{T}  \cong (\leftComplex{\lefty{T}} \boxtimes \rightComplex{\righty{T}})$.\\
\ \\
\noindent Furthermore, changing either $\leftComplex{\lefty{T}}$ by a homotopy equivalence (of type A structures) or $\rightComplex{\righty{T}}$ (of type $D$ structures) changes $\leftComplex{\lefty{T}} \boxtimes \rightComplex{\righty{T}}$ by a chain homotopy equivalence. Thus we can construct and simplify the two factors independently of each other, and then combine them using the $\boxtimes$-construction.  \\
\ \\
\noindent This provides a fully modular approach to constructing Khovanov homology at the level of bigraded modules. In particular we can compute the structures for tangles and then combine them. For example, in section \ref{sec:examplesA} we will compute the type $A$ structures for tangles underlying the three Reidemeister moves, and simplify them, to see that they are homotopy equivalent to the type $A$ structure after applying the move. The pairing through $\boxtimes$ immediately implies the Reidemeister invariance of Khovanov homology. In short, we obtain a more theoretical and convenient means for understanding local modifications of link diagrams and their effects on the global Khovanov homology. \\
\ \\
\noindent Once these constructions are understood, it is straightforward to define a type $DA$-bimodule for an $(m,n)$-tangle. At least for $\Z/2\Z$-coefficients, this bimodule will satisfy all the properties required for the algebra developed in the bordered Floer theory of such bimodules. In particular, we will be able to understand their Hochschild homologies directly. For now, we content ourselves with an example. \\
\ \\
\noindent{\bf Example:} Consider the following connect sum:
$$
\inlinediag[0.5]{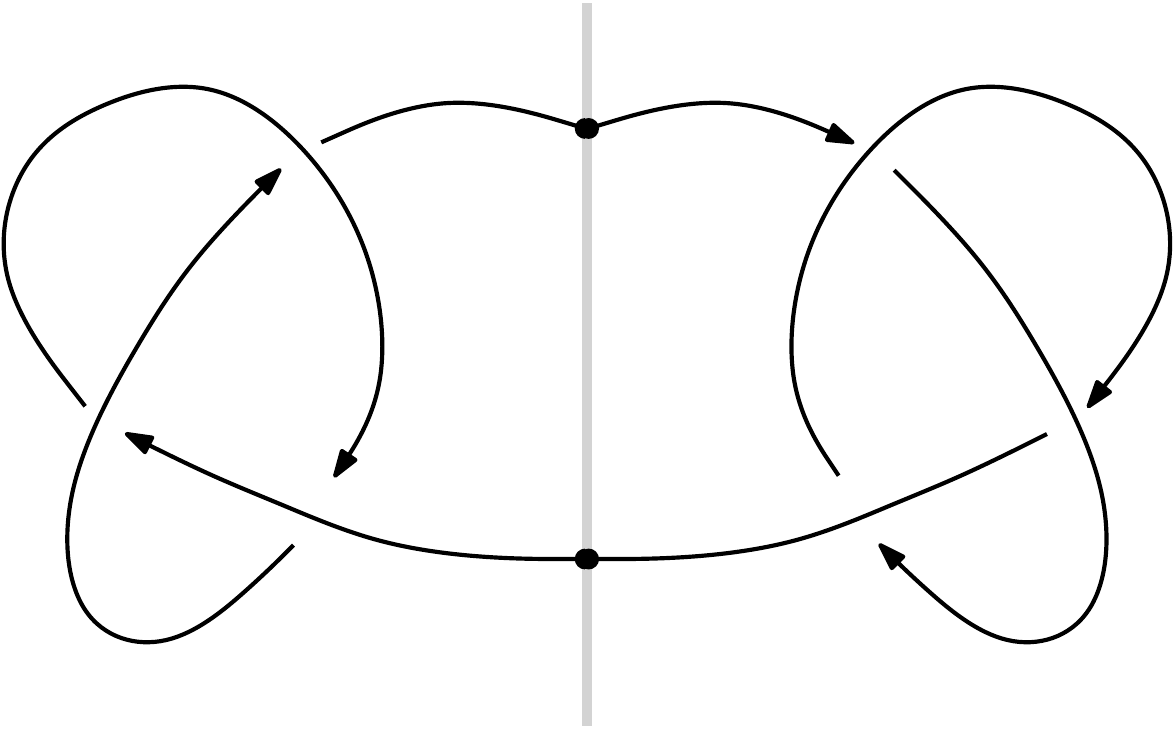}
$$
It was shown in \cite{typeD} that $\mathcal{B}\Gamma_{1}$ is the quiver algebra for
$$
\inlinediag{bridgeGraph1}
$$
Not all the $\mathcal{B}\Gamma_{n}$ are quiver algebras, however, and describing them takes some work. In \cite{typeD} we showed that, for this example, $\rightComplex{\righty{T}}$ is homotopy equivalent to 
\begin{equation}
\begin{split}
\righty{\delta}(s_{(-3,-15/2)}^{+}) &= 2\,\righty{e_{C}} \otimes s_{(-2,-13/2)}^{-} + \lefty{e_{C}} \otimes s_{(-3,-17/2)}^{-} \\
\righty{\delta}(s_{(-2,-11/2)}^{+}) &= - \lefty{e_{C}} \otimes s_{(-2,-13/2)}^{-}\\
\righty{\delta}(s_{(0,-3/2)}^{+})\ &= - \lefty{e_{C}} \otimes s_{(0,-5/2)}^{-}\\
\end{split}
\end{equation}
where the $+$ and $-$ superscripts identify which idempotent acts as the identity on the generator.  Furthermore, the type $A$ structure $\leftComplex{\lefty{T}}$ is homotopy equivalent to one which also has six generators: $t^{+}_{(0,5/2)}$, $t^{+}_{(2,13/2)}$, $t^{+}_{(3,17/2)}$, $t^{-}_{(0,3/2)}$, $t^{-}_{(2,11/2)}$, $t^{-}_{(3,15/2)}$. For these generators $d_{APS} \equiv 0$, the action of $\righty{e_{C}}$ is given by
$t^{+}_{(0,5/2)} \rightarrow t^{-}_{(0,3/2)}$, $t^{+}_{(2,13/2)} \rightarrow t^{-}_{(2,11/2)}$, $t^{+}_{(3,17/2)} \rightarrow t^{-}_{(3,15/2)}$, and the action of $\lefty{e_{C}}$ is given by $t^{+}_{(2,13/2)} \rightarrow 2\cdot t^{-}_{(3,15/2)}$. The complex $\leftComplex{\lefty{T}} \boxtimes \rightComplex{\righty{T}}$ 
can be computed exactly (see section \ref{sec:expairing} for the details). It has homology with free part
$$
\Z_{(-3,7)} \oplus \Z_{(-2,-3)} \oplus \Z_{(-1,-3)} \oplus \Z^{2}_{(0,-1)} \oplus \Z^{2}_{(0,1)} \oplus \Z_{(1,3)} \oplus \Z_{(2,3)} \oplus \Z_{(3,7)}
$$
and torsion part
$$
\big(\Z/2\Z)_{(-2,-5)} \oplus \big(\Z/2\Z)_{(0,-1)} \oplus \big(\Z/2\Z)_{(1,1)} \oplus \big(\Z/2\Z)_{(3,5)}
$$
which is the Khovanov homology of this knot. \\
\ \\
\noindent {\em Degree shift convention:} If $M$ is a $\Z$-graded module, $M[n]$ is the graded module with $(M[n])_{i} = M_{i-n}$, i.e. the module found by shifting the homogeneous elements of $M$ up $n$ levels. If $m \in M$, the corresponding element in $M[n]$ will be denoted $m[n]$. Thus $\mathrm{gr}(m[n]) = \mathrm{gr}(m) + n$. \\
\ \\
\noindent {\bf Note:} After posting this paper to the arXiv, the author was informed by Cotton Seed that he had independently discovered a similar constructions of a type A structure in Khovanov homology. The author would like to thank Andy Manion for finding an error in the original version of this paper. 
\section{The algebra from cleaved links}\label{sec:algebra}

\noindent We summarize the construction of the algebra $\mathcal{B}\Gamma_{n}$ from \cite{typeD}. 

\subsection{Cleaved planar links}

\noindent Let $P_{n}$ be the set of points $p_{1}=(0,1),\ldots, p_{2n}=(0,2n)$ on the $y$-axis of $\R^{2}$, ordered by the second coordinate. We denote the closed half-plane $(-\infty,0] \times \R \subset \R^{2}$ by $\leftHalf$ while  $\rightHalf = [0,\infty) \times \R$. 

\begin{defn}
A {\em $n$-cleaved link} $L$ is an embedding of circles in $\R^{2}$ such that 
\begin{enumerate}
\item the circles of $L$ are disjoint and transverse to the $y$-axis,
\item each point in $P_{n}$ is on a circle in $L$,
\item each circle in $L$ contains at least two points in $P_{n}$
\end{enumerate}
The set of circle components of $L$ will be denoted $\circles{L}$.
\end{defn}

\noindent We take two $n$-cleaved links to be equivalent if they are related by an isotopy of $\R^{2}$ which pointwise fixes the $y$-axis, and by reversing the orientation on each circle. We will denote the equivalence classes by $\cleave{n}$

\begin{defn}
The {\em constituents} of an $n$-cleaved link $L$ are the planar matchings:
\begin{equation}
\lefty{L} = \leftHalf \cap L  \hspace{1in}
\righty{L} = \rightHalf \cap L 
\end{equation}
\end{defn}

\begin{defn}
A {\em bridge} for a cleaved link $L$ is an embedding $\gamma : [0,1] \rightarrow \R^{2}\backslash \big(\{0\}\times \R\big)$ such that 
\begin{enumerate}
\item $\gamma(0)$ and  $\gamma(1)$ are on {\em distinct} arcs of $\lefty{L}$ or $\righty{L}$
\item the image under $\gamma$ of $(0,1)$ is disjoint from $L$ 
\end{enumerate} 
\end{defn}

\noindent By definition, a bridge $\gamma$ has image in either $\rightHalf$ or $\leftHalf$. We call this half plane the {\em location} of the bridge. Bridges are also considered up to isotopy fixing the $y$-axis. 

\begin{defn}
The equivalence classes of bridges for a cleaved link $L$ will be denoted $\bridges{L}$. $\bridges{L} = \leftBridges{L} \cup \rightBridges{L}$ where $\leftBridges{L}$ consists of those equivalence classes in $\leftHalf$ and $\rightBridges{L}$ consists of those classes in $\rightHalf$.
\end{defn}

\noindent For each class of bridges $\gamma \in \bridges{L}$ we can construct a new cleaved planar link. 

\begin{defn}
Let $L$ be an equivalence class of cleaved links and let $\gamma \in \bridges{L}$. $L_{\gamma}$ is the equivalence class of cleaved links found by surgery along $\gamma$. 
\end{defn}

\noindent $L_{\gamma}$ has a special bridge $\gamma^{\dagger}$ introduced by the surgery. More specifically, there is a neighborhood of $\gamma$ homeomorphic to $[-1,1] \times [-1,1]$ which intersects $L$ along $\{\pm 1\}\times[-1,1]$, and for which $\gamma$ is $[-1,1] \times \{0\}$ (i.e. the core). $L_{\gamma}$ results from removing these two arcs from $L$ and replacing them with $[-1,1] \times \{\pm 1\}$. $\gamma^{\dagger}$ is then the bridge for $L_{\gamma}$ defined by the arc $\{0\} \times [-1,1]$.  

\begin{defn}
\begin{enumerate}
\item The support of $\gamma \in \bridges{L}$ is the set of {\em three} circles in $L$ and $L_{\gamma}$ which contain the feet of $\gamma$ and $\gamma^{\dagger}$
\item $\merge{L}$ is the subset of $\gamma \in \bridges{L}$ where surgery on $\gamma$ merges two circles $\{\startCircle{\gamma}, \terminalCircle{\gamma}\}$. In this case, $C_{\gamma}$ is the circle in $\circles{L_{\gamma}}$ which contains both feet of $\gamma^{\dagger}$.
\item $\fission{L}$ is the subset of $\gamma \in \bridges{L}$ where surgery on $\gamma$ divides a circle $C$ of $L$. In this case, $C^{a}_{\gamma}$ and $C^{b}_{\gamma}$ are the circles in $\circles{L_{\gamma}}$ which contain the feet of $\gamma^{\dagger}$.   
\end{enumerate}
\end{defn}

\begin{prop}
Given any bridge $\gamma$ for $L$, $\bridges{L}\backslash\{\gamma\}$ can be decomposed as a disjoint union $B_{\pitchfork}(L,\gamma) \cup B_{|\,|}(L,\gamma)$, where
\begin{enumerate}
\item $B_{\pitchfork}(T,\gamma)$ consists of the classes of bridges all of whose representatives intersect $\gamma$
\item $B_{|\,|}(L,\gamma)$ consists of the classes of bridges containing a representative which does not intersect $\gamma$,
\end{enumerate}  
Furthermore, we can divide $B_{|\,|}(L,\gamma)$ into the disjoint union $B_{d}(L,\gamma) \cup B_{s}(L,\gamma) \cup B_{o}(L,\gamma)$ where
\begin{enumerate}
\item $B_{d}(L,\gamma)$ consists of those bridges neither of whose ends is on an arc with $\gamma$,
\item $B_{s}(L,\gamma)$ consists of those bridges with a single end on the same arc as $\gamma$ and lying on the {\em same} side of the arc as $\gamma$
\item $B_{o}(L,\gamma)$ consists of those bridges with a single end on the same arc as $\gamma$ and lying on the {\em opposite} side of the arc as $\gamma$
\end{enumerate}
\end{prop}

\noindent If $\eta \in B_{\circ}(L,\gamma)$ where $\circ$ represents a specific choice of one of the above then $\gamma \in B_{\circ}(L,\eta)$. Furthermore, if $\delta$ has a different location than $\gamma$ then $\delta \in B_{d}(L,\gamma)$. We consider how these sets change under surgery on $\gamma$.

\begin{prop} 
Surgery on $\gamma$ induces an identification $B_{d}(L,\gamma)$ with $B_{d}(L_{\gamma}, \gamma^{\dagger})$ and a $2:1$ map $B_{s}(L,\gamma) \lra B_{o}(L_{\gamma},\gamma^{\dagger})$. 
\end{prop}

\noindent Dually there is a $2:1$ map $B_{s}(L_{\gamma}, \gamma^{\dagger}) \lra B_{o}(L,\gamma)$.\\
\ \\
\noindent{\bf Proof:} Let $\eta \in B_{|\,|}(L,\gamma)$. Pick a representative arc for $\eta$ which does not intersect the representative arc for $\gamma$. Then $\eta$ also represents a bridge in $B_{|\,|}(L_{\gamma}, \gamma^{\dagger})$, and vice-versa. If $\eta \in B_{d}(L,\gamma)$, any isotopy of the representaive arc occurs in a region disjoint from $\gamma$ and its endpoints, since the isotopy will occur along arcs disjoint from those intersecting $\gamma$. This isotopy also survives into $(L_{\gamma}, \gamma^{\dagger})$. Reversing this construction for $B_{d}(L_{\gamma}, \gamma^{\dagger})$ proves the identification. Note that an isotopy of $\eta \in B_{s}(L,\gamma)$ missing $\gamma$ can likewise be pushed forward. However, for each $\eta$ we can slide $\eta$ over $\gamma$ to get another bridge $\eta'\in B_{s}(L,\gamma)$. In $L_{\gamma}$ $\eta \simeq \eta'$ and both are on the opposite side of $\gamma^{\dagger}$. By looking at a local model, this is the only type of collision, so the map is $2:1$ on $B_{s}(L,\gamma)$. We can apply the same argument to $(L_{\gamma},\gamma^{\dagger})$ to obtain the $2:1$ map in the other direction. $\Diamond$. \\
\ \\
\noindent A {\em decoration} for an $n$-cleaved link $L$ is a map $\sigma\!: \circles{L} \lra \{+,-\}$. 

\begin{defn}
$\cleaved{n}$ is the set of decorated, $n$-cleaved links:
\begin{equation}
\cleaved{n}= \big\{\,(L,\sigma)\,\big|\,L\in \cleave{n}, \sigma\mathrm{\ is\ a\ decoration\ for\ }L\big\}
\end{equation} 
\end{defn}

\noindent We will often restrict a decoration $\sigma$ of $L$ to give decorations on the arcs of its constituents $\lefty{L}$ and  $\righty{L}$. In addition, we will need to following statistic for a decorated, cleaved link:
\begin{equation}
\iota(L,\sigma) = \#\big\{C \in \circles{L} \big| \sigma(C) = +\big\} - \#\big\{C \in \circles{L} \big| \sigma(C) = -\big\}
\end{equation}

\subsection{The algebra $\mathcal{B}\Gamma_{n}$}

\noindent We will describe $\mathcal{B}\Gamma_{n}$ by generators and relations. First, there is an idempotent for each decorated, cleaved link in $\cleaved{n}$. We will denote the idempotent corresponding to $(L,\sigma)$ by $I_{(L,\sigma)}$. The idempotents will be orthogonal to each other. 

\begin{defn}
$\mathcal{I}_{n}$ is the sub-algebra generated by the idempotents $I_{(L,\sigma)}$.
\end{defn}

\noindent For each $(L,\sigma) \in \cleaved{n}$ we specify certain elements in $I_{(L,\sigma)}\mathcal{B}\Gamma_{n}$. $\mathcal{B}\Gamma_{n}$ will be freely generated by the idempotents and these elements, subject to the relations described below.
\begin{enumerate}
\item For each circle $C \in \circles{L}$ with $\sigma(C) = +$ there are two elements $\righty{e_{C}}$ and $\lefty{e_{C}}$ in 
$I_{(L,\sigma)}\mathcal{B}\Gamma_{n}$. Furthermore, $\righty{e_{C}}I_{(L,s_{C})} = \righty{e_{C}}$ for the decoration with $s_{C}(C) = -$ and $s_{C}(D) = s(D)$ for each $D \in \circles{L}\backslash\{C\}$, while $\righty{e_{C}}I_{(L',s')} = 0$ for every other idempotent. The same relations hold for $\lefty{e_{C}}$. These types of elements are called {\em decoration elements}, while the $C$ above is called the {\em support} of the element.
\item Let $\gamma \in \bridges{L}$, then there is a {\em bridge element} $e_{(\gamma; \sigma , \sigma_{\gamma})}$ with $I_{(L,\sigma)}e_{(\gamma; \sigma , \sigma_{\gamma})} = e_{(\gamma; \sigma , \sigma_{\gamma})}I_{(L_{\gamma},\sigma_{\gamma})} = e_{(\gamma; \sigma , \sigma_{\gamma})}$ in each of the following cases, based on the decorations,
\begin{enumerate}
\item\label{bridge1} when  $\gamma \in \merge{L}$ and $\sigma$ and $\sigma_{\gamma}$ restrict to the support of $\gamma$ as one of 
$$
\begin{array}{llcl}
\sigma(\startCircle{\gamma})=+ & \sigma(\terminalCircle{\gamma})= + &\hspace{1in} &\sigma_{\gamma}(C_{\gamma})= + \\
\sigma(\startCircle{\gamma})=- & \sigma(\terminalCircle{\gamma})= + &\  &\sigma_{\gamma}(C_{\gamma})= - \\
\sigma(\startCircle{\gamma})=+ & \sigma(\terminalCircle{\gamma})= - &\ &\sigma_{\gamma}(C_{\gamma})= - \\
\end{array}
$$ 
and $s(D)=s_{\gamma}(D)$ on every circle not in the support of $\gamma$;
\item\label{bridge2} when $\gamma \in \fission{L}$, $C \in \circles{L}$ is the circle containing both feet of $\gamma$, and
\begin{enumerate}
\item  $\sigma(C) = +$,  if $\sigma$ and $\sigma_{\gamma}$ restrict to the support of $\gamma$ as either of
$$
\begin{array}{llcl}
\sigma(C)=+ &\hspace{1in} & \sigma_{\gamma}(C^{a}_{\gamma})= + & \sigma_{\gamma}(C^{b}_{\gamma})= - \\
\sigma(C)=+ &\ & \sigma_{\gamma}(C^{a}_{\gamma})= - & \sigma_{\gamma}(C^{b}_{\gamma})= + \\
\end{array}
$$
\item $\sigma(C) = -$, if $\sigma$ and $\sigma_{\gamma}$ restrict to the support of $\gamma$ as 
$$
\begin{array}{llcl}
\sigma(C)=- & \hspace{1in} & \sigma^{-}_{\gamma}(C^{a}_{\gamma})= - & \sigma^{-}_{\gamma}(C^{b}_{\gamma})= - \\
\end{array}
$$
\end{enumerate}
and $s(D)=s_{\gamma}(D)$ on every circle not in the support of $\gamma$.
\end{enumerate}
\end{enumerate}

\noindent In \cite{typeD} we note that with these generators and idempotents
\begin{prop}
$\mathcal{B}\bridgeGraph{n}$ is finite dimensional
\end{prop}

\noindent Furthermore, $\mathcal{B}\bridgeGraph{n}$ can be given a bigrading, \cite{typeD}. On the generating elements the bigrading is specified by setting
$$
\begin{array}{lcl}
\mathrm{I}_{(L,\sigma)} & \longrightarrow & (0,0)\\
\righty{e_{C}} & \longrightarrow & (0,-1)\\
\lefty{e_{C}} & \longrightarrow & (1,1)\\
\righty{e_{\gamma}} &\longrightarrow & (0,-1/2)\\
\lefty{e_{\gamma}} &\longrightarrow & (1,1/2)\\
\end{array}
$$
\noindent On every other element it is computed by extending the above homomorphically. The first entry of this bigrading will be denoted $\leftnorm{\alpha}$, while the second element will be denoted $q(\alpha)$.\\
\ \\
\noindent We now turn to describing the relations between these generators. Each of the relations is homogeneous with respect to the bigrading. First, there are a number of ``graded commutativity'' relations, based on the first entry in the bigrading:
\begin{equation}
e_{\alpha}e_{\beta'} = (-1)^{\leftnorm{e_{\alpha}}\leftnorm{e_{\beta}}}e_{\beta}e_{\alpha'} 
\end{equation}
This graded commutativity occurs in the following cases, assuming that $I_{(L,\sigma)}e_{\alpha} \neq 0$ and $I_{(L,\sigma)}e_{\beta} \neq 0$, 
\begin{enumerate}
\item If $e_{\alpha}$ and $e_{\beta}$ are decoration elements for distinct circles $C$ and $D$ in $(L,\sigma)$ with $\sigma(C)=\sigma(D)=+$, and $e_{\alpha'}$ is the decoration element for $D$ in $(L,\sigma_{C})$, while $e_{\beta'}$ is the
decoration element for $C$ in $(L,\sigma_{D})$.
\item If $e_{\alpha}=e_{(\gamma,\sigma,\sigma')}$ for a bridge $\gamma$ in $(L,\sigma)$ and $e_{\beta}$ is a decoration element for $C \in \circles{L}$, with $C$  not in the support of $\gamma$, while $e_{\alpha'} = e_{(\gamma,\sigma_{C},\sigma_{C}')}$ and $e_{\beta'}$ is the decoration element for $C$ in $(L_{\gamma},s')$. Due to the disjoint support, there will always be a pair of such elements.
\item If $e_{\alpha}= e_{(\gamma,\sigma,\sigma')}$ and $e_{\beta}=e_{(\eta,\sigma,\sigma'')}$ are bridge elements for distinct bridges $\gamma$ and $\eta$ in $(L,\sigma)$, with $\eta \in B_{d}(L,\gamma)$ and $e_{\beta'}= e_{(\eta,\sigma',\sigma''')}$ and $e_{\alpha'}=e_{(\gamma,\sigma'',\sigma''')}$ for some decoration $\sigma'''$ on $L_{\gamma,\eta}$.
\item If $e_{\alpha}= e_{(\righty{\gamma},\sigma,\sigma')}$ and $e_{\beta}=e_{(\righty{\eta},\sigma,\sigma'')}$ are bridge elements for distinct {\em right} bridges $\gamma$ and $\eta$ in $(L,\sigma)$, and $e_{\beta'}= e_{(\righty{\delta},\sigma',\sigma''')}$ and $e_{\alpha'}=e_{(\righty{\omega},\sigma'',\sigma''')}$ %for $\righty{\delta} \in B_{|\,|}(L_{\righty{\gamma}}, \righty{\gamma}^{\dagger})$, $\righty{\omega} \in B_{|\,|}(L_{\righty{\eta}}, \righty{\gamma}^{\dagger})$
, such that $L_{\gamma, \delta} = L_{\eta,\omega}$, and some comptatible decoration $\sigma'''$.
\item If $e_{\alpha}= e_{(\lefty{\gamma},\sigma,\sigma')}$ and $e_{\beta}=e_{(\lefty{\eta},\sigma,\sigma'')}$ are bridge elements for distinct {\em left} bridges in $(L,\sigma)$, with $\lefty{\eta} \in B_{o}(L,\lefty{\gamma})$, and $e_{\beta'}= e_{(\lefty{\delta},\sigma',\sigma''')}$, $e_{\alpha'}=e_{(\lefty{\omega},\sigma'',\sigma''')}$ with $L_{\gamma, \delta} = L_{\eta,\omega}$, and some comptatible decoration $\sigma'''$.
\end{enumerate}

\noindent We also note that the type (decoration vs. bridge) and the location are the same for $e_{\alpha}$ and $e_{\alpha'}$ as well as for the pair $e_{\beta}$ and $e_{\beta'}$. In fact, in all these relations elements of the algebra from $\rightHalf$ will act like even elements for the $\Z/2\Z$-grading from $\leftnorm{\alpha}$, while elements from $\leftHalf$ act like {\em odd} elements.\\
\ \\
\noindent {\bf Other bridge relations:} Suppose $\gamma \in \leftBridges{L}$ and $\eta \in B_{\pitchfork}(L_{\gamma},\gamma^{\dagger})$, then 
$$
e_{(\gamma,\sigma,\sigma')}e_{(\eta,\sigma',\sigma'')} = 0
$$
whenever $\sigma'$ and $\sigma''$ are compatible decorations.\\
\ \\
\noindent Furthermore, suppose that there is a circle $C \in \circles{L}$ with $\sigma(C)=+$, and there are elements $\righty{e}_{(\gamma,\sigma,\sigma')}$ and $\righty{e}_{(\gamma^{\dagger},\sigma',\sigma_{C})}$ for a bridge $\gamma \in \rightBridges{L}$ then
\begin{equation}\label{rel:special}
\righty{e}_{(\gamma,\sigma,\sigma')}\righty{e}_{(\gamma^{\dagger},\sigma',\sigma_{C})} = \righty{e_{C}}
\end{equation}
Such a circle $C$ is unique for the choice of $\gamma$ and $\sigma'$ and is called the {\em active circle} for $\gamma$.  \\
\ \\
\noindent Finally, suppose that $\lefty{\alpha} \in B_{s}(L,\lefty{\gamma})$. Let $\lefty{\beta}$ be the bridge obtained by sliding $\lefty{\alpha}$ over $\lefty{\gamma}$. Let $\lefty{\delta}$ be the image of $\lefty{\alpha}$ and $\lefty{\beta}$ in $L_{\lefty{\gamma}}$, $\lefty{\zeta}$ be the image of $\lefty{\alpha}$ and $\lefty{\gamma}$ in $L_{\lefty{\beta}}$ and $\lefty{\eta}$ be the image of $\lefty{\beta}$ and $\lefty{\gamma}$ in $L_{\lefty{\alpha}}$. Then
$$
\lefty{e_{\alpha}}\lefty{e_{\eta}} + \lefty{e_{\beta}}\lefty{e_{\zeta}} + \lefty{e_{\gamma}}\lefty{e_{\delta}} = 0
$$
whenever there are comptaible decorations on $L_{\alpha}$, $L_{\beta}$, $L_{\gamma}$, and $L_{\alpha,\eta}=L_{\beta,\zeta} = L_{\gamma,\delta}$, for the paired edges to exist.\\
\ \\
\noindent{\bf Relations for decoration edges:} When the support of $e_{C}$ is not disjoint from that of $\righty{e}_{(\gamma,\sigma,\sigma_{\gamma})}$ the relations are different depending upon the location of $e_{C}$.

\begin{enumerate}
\item {\bf The relations for $\righty{e_{C}}$:} Suppose that $\gamma \in \merge{L}$ merges $C_{1}$ and $C_{2}$ to get $C \in \circles{L_{\gamma}}$, and $\sigma(C_{1}) = \sigma(C_{2}) = +$, then
\begin{equation}
\righty{e_{C_{1}}}m_{(\gamma,\sigma_{C_{1}},\sigma_{C})} = \righty{e_{C_{2}}}m_{(\gamma,\sigma_{C_{2}},\sigma_{C})} = m_{(\gamma,\sigma,\sigma_{\gamma})} \righty{e_{C}}
\end{equation}
Note that if $\sigma(C_{i}) = -$ for either $i=1$ or $2$, then there is no relation imposed. \\
\ \\
Dually, if surgery on $\gamma \in \fission{L}$ divides circle $C \in \circles{L}$ into $C_{1}$ and $C_{2}$ in $\circles{L_{\gamma}}$, and $\sigma$ assigns $+$ to $C$, then 
\begin{equation}
\righty{e_{C}}f_{(\gamma,\sigma_{C},\sigma_{C,\gamma})} =  f_{(\gamma,\sigma,\sigma^{1}_{\gamma})}\righty{e_{C_{1}}} = f_{(\gamma,\sigma,\sigma^{2}_{\gamma})} \righty{e_{C_{2}}}
\end{equation}
where $\sigma^{i}_{\gamma}$ assigns $+$ to $C_{i}$ and $-$ to $C_{3-i}$. 

\item {\bf The relations for $\lefty{e_{C}}$:} Suppose that $\gamma \in \rightMerge{L}$ merges $C_{1}$ and $C_{2}$ to get $C \in \circles{L_{\gamma}}$, and $\sigma(C_{1}) = \sigma(C_{2}) = +$, then:
\begin{equation}\label{rel:lefty1}
\lefty{e_{C_{1}}}m_{(\gamma,\sigma_{C_{1}},\sigma_{C})} + \lefty{e_{C_{2}}}m_{(\gamma,\sigma_{C_{2}},\sigma_{C})} - m_{(\gamma,\sigma,\sigma_{\gamma})} \lefty{e_{C}} = 0
\end{equation}
and when $\sigma(C) = +$ and $\gamma \in \rightFission{L}$ divides $C$ into $C_{1}$ and $C_{2}$
\begin{equation}\label{rel:lefty2}
\lefty{e_{C}}f_{(\gamma,\sigma_{C},\sigma_{C,\gamma})} +  f_{(\gamma,\sigma,\sigma^{1}_{\gamma})}\lefty{e_{C_{1}}} - f_{(\gamma,\sigma,\sigma^{2}_{\gamma})} \lefty{e_{C_{2}}} = 0
\end{equation}
whereas if $\gamma \in \leftMerge{L}$ merges $C_{1}$ and $C_{2}$ to get $C \in \circles{L_{\gamma}}$, and $\sigma(C_{1}) = \sigma(C_{2}) = +$, then:
\begin{equation}\label{rel:lefty3}
\lefty{e_{C_{1}}}m_{(\gamma,\sigma_{C_{1}},\sigma_{C})} + \lefty{e_{C_{2}}}m_{(\gamma,\sigma_{C_{2}},\sigma_{C})} + m_{(\gamma,\sigma,\sigma_{\gamma})} \lefty{e_{C}} = 0
\end{equation}
and when $\sigma(C) = +$ and $\gamma \in \leftFission{L}$ divides $C$ into $C_{1}$ and $C_{2}$
\begin{equation}\label{rel:lefty4}
\lefty{e_{C}}f_{(\gamma,\sigma_{C},\sigma_{C,\gamma})} +  f_{(\gamma,\sigma,\sigma^{1}_{\gamma})}\lefty{e_{C_{1}}} + f_{(\gamma,\sigma,\sigma^{2}_{\gamma})} \lefty{e_{C_{2}}} = 0
\end{equation}
\end{enumerate}

\subsection{A differential on $\mathcal{B}\bridgeGraph{n}$}

\noindent Surgery along a bridge $\gamma \in \leftBridges{L}$ followed surgery on $\gamma^{\dagger}$ {\em does not} correspond to a relation (compare relation \ref{rel:special}). Instead these products occur in a differential on $\mathcal{B}\bridgeGraph{n}$.

\begin{prop}\cite{typeD}
Let $(L,\sigma) \in \cleaved{n}$ such that there is a circle $C \in \circles{L}$ with $\sigma(C) = +$. Let $\lefty{e_{C}}$ be
the decoration element corresponding to $C$. Let
\begin{equation}\label{eq:dga}
d_{\Gamma_{n}}(\lefty{e_{C}}) = -\sum e_{(\gamma,\sigma,\sigma_{\gamma})}e_{(\gamma^{\dagger},\sigma_{\gamma},\sigma_{C})}
\end{equation}
where the sum is over all $\gamma \in \leftBridges{L}$ with $C$ as active circle, and all decorations $\sigma_{\gamma}$ which define
compatible elements. Let $d_{\Gamma_{n}}(e) = 0$ for every other generator $e$ (including idempotents). Then $d_{\Gamma_{n}}$ can be extended to a $(1,0)$ differential on bigraded algebra $\mathcal{B}\bridgeGraph{n}$ which satisfies the following Leibniz identity:
\begin{equation}\label{eqn:Leibniz}
d_{\Gamma_{n}}(\alpha\beta) = (-1)^{\leftnorm{\beta}}(d_{\Gamma_{n}}(\alpha)\big)\beta + \alpha\big(d_{\Gamma_{n}}(\beta)\big)
\end{equation}
\end{prop}

\noindent $(\mathcal{B}\bridgeGraph{n}, d_{\Gamma_{n}})$ denotes this differential, bigraded $\Z$-algebra.  \\
\ \\
\subsection{Example:} $(\mathcal{B}\bridgeGraph{1}, d_{\Gamma_{1}})$: $P_{1}$ consists of two points, and there is only one planar matching in $\leftHalf$ and $\rightHalf$. Consequently, the only $1$-cleaved link is a circle intersecting the $y$-axis in two points. Thus, there are two vertices in $\bridgeGraph{1}$: when this circle is decorated with a $+$ and when it is decorated with a $-$. We will call these $C^{\pm}$. There are no bridges in either $\leftHalf$ or $\rightHalf$, so the only edges are $\lefty{e}_{C}: C^{+} \longrightarrow C^{-}$ and  $\righty{e}_{C}: C^{+} \longrightarrow C^{-}$. Thus $\bridgeGraph{1}$ looks like  

\begin{center}
\includegraphics[scale=0.5]{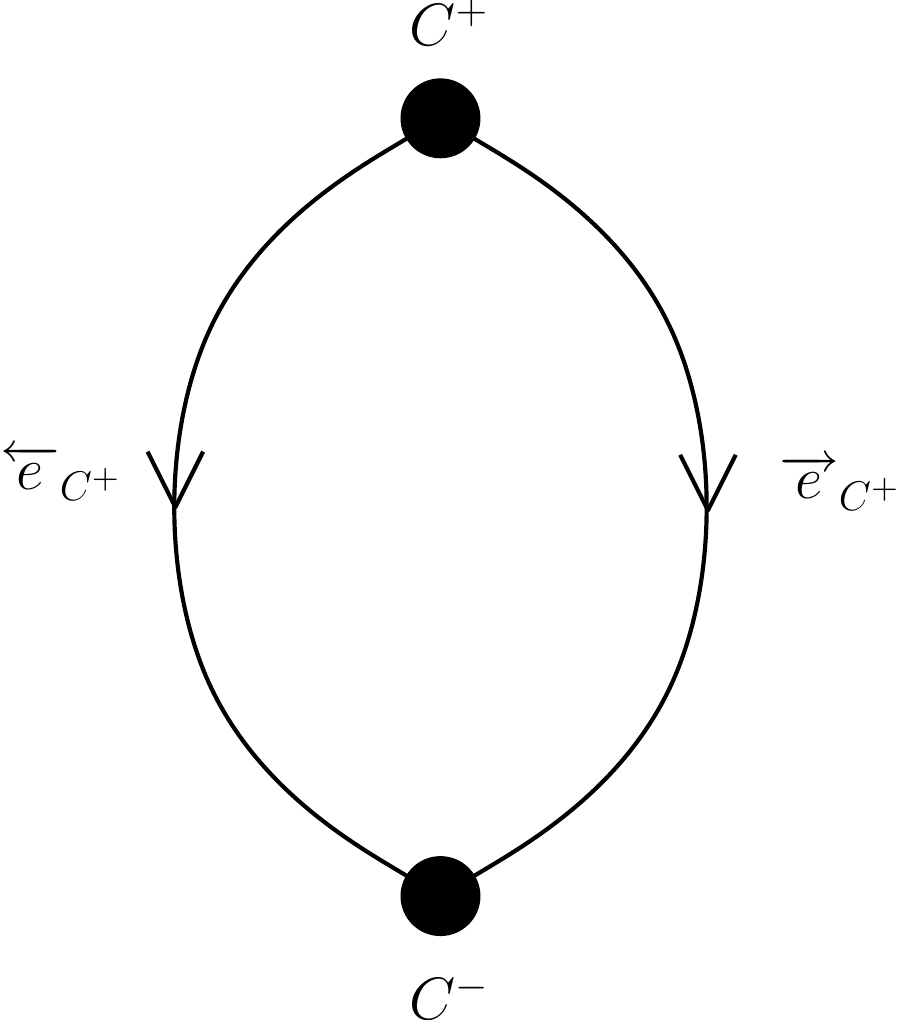}
\end{center}

\noindent Thus, $\mathcal{B}\bridgeGraph{1}$ consists of four elements $I_{C^{+}}, I_{C^{-}}$ in grading $(0,0)$, $\lefty{e}_{C}$ in grading $(1,1)$, and $\righty{e}_{C}$ in grading $(0,-\frac{1}{2})$. The product of any two of these is trivial except for the actions of the idempotents: $I_{C^{+}}\lefty{e}_{C} = \lefty{e}_{C} = \lefty{e}_{C}I_{C^{-}}$, and similarly for $\righty{e}_{C}$. The differential $d_{\Gamma_{1}} \equiv 0$ since its image is in the set generated by paths of bridge edges.\\
\ \\
\noindent For more detail about $(\mathcal{B}\bridgeGraph{2}, d_{\Gamma_{2}})$ see the examples in section 2 of \cite{typeD}.

\section{Tangles and Resolutions}\label{sec:APS}

\noindent In this section we recall the notions of tangles and resolutions used in \cite{typeD}, adapting them to the case at hand. For more detail please consult \cite{typeD}.\\
\ \\
\noindent  Let $\lefty{\mathbb{R}^{3}} = \leftHalf \times \R$ be the half space corresponding to $\leftHalf \subset \R^2$ under the standard projection $\pi$ to the $xy$-plane. 

\begin{defn}
An (inside) tangle $\lefty{\tangle{T}}$ is a smooth, proper embedding of 
\begin{enumerate}
\item[] i) $n$ copies of the interval $[0,1]$, and 
\item[] ii) $k$ copies of $S^{1}$
\end{enumerate}
in $\lefty{\mathbb{R}^{3}}$, whose boundary is the set of $2n$ points $P_{n}$ in $\partial \leftHalf$. $\lefty{\tangle{T}_{1}}$ and $\lefty{\tangle{T}_{2}}$ are equivalent if there is an isotopy of $\lefty{\mathbb{R}^{3}}$ taking $\lefty{\tangle{T}_{1}}$ to $\lefty{\tangle{T}_{2}}$ and pointwise fixing the boundary $\partial\lefty{\mathbb{R}^{3}}$.  
\end{defn}

\noindent As usual, we will study $\lefty{\tangle{T}}$ through its tangle diagrams in $\leftHalf$. Different diagrams for $\righty{\tangle{T}}$ are related by sequences of Reidemeister moves, and planar isotopies, in the interior of $\leftHalf$. We will denote a tangle diagram for a tangle by the corresponding roman letter:  $\lefty{T}$ will be a diagram for $\lefty{\tangle{T}}$. \\
\ \\
\noindent The crossings of $\lefty{T}$ form a set $\cross{\lefty{T}}$. We will orient $\lefty{\tangle{T}}$ and use the usual convention for positive and negative crossings: \\
$$
\inlinediag[0.5]{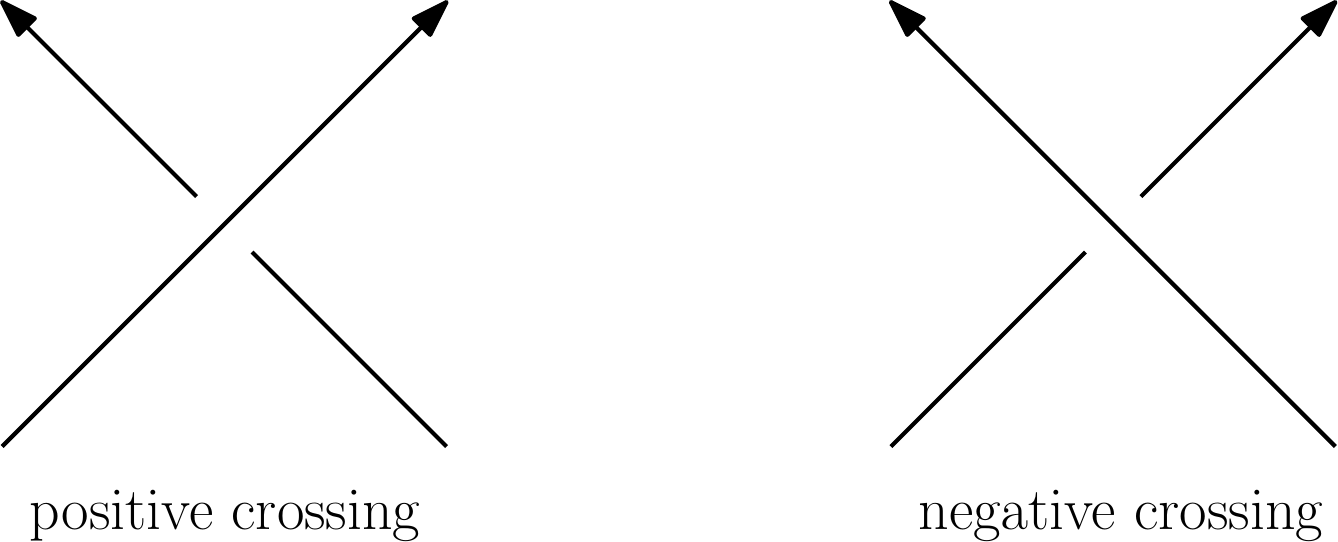} 
$$
\noindent The number of positive/negative crossings will be denoted $n_{\pm}(\lefty{T})$.

\subsection{Resolutions}

\begin{defn}
A resolution $r$ of $\lefty{T}$ is a pair $(\rho, \lefty{m})$ where $\rho: \cross{\lefty{T}} \lra \{0,1\}$, and $\righty{m}$ is a planar matching of $P_{2n}$ embedded in $\rightHalf$. The resolution diagram, $r(\lefty{T})$ is the crossingless, planar link in $\leftHalf$ obtained by 1) gluing $\lefty{T} \subset \leftHalf$ to $\righty{m} \subset \rightHalf$, and 2) locally replacing (disjoint) neighborhoods of each crossing $c \in \cross{\lefty{T}}$ using the following rule:\\
$$
\inlinediag[0.5]{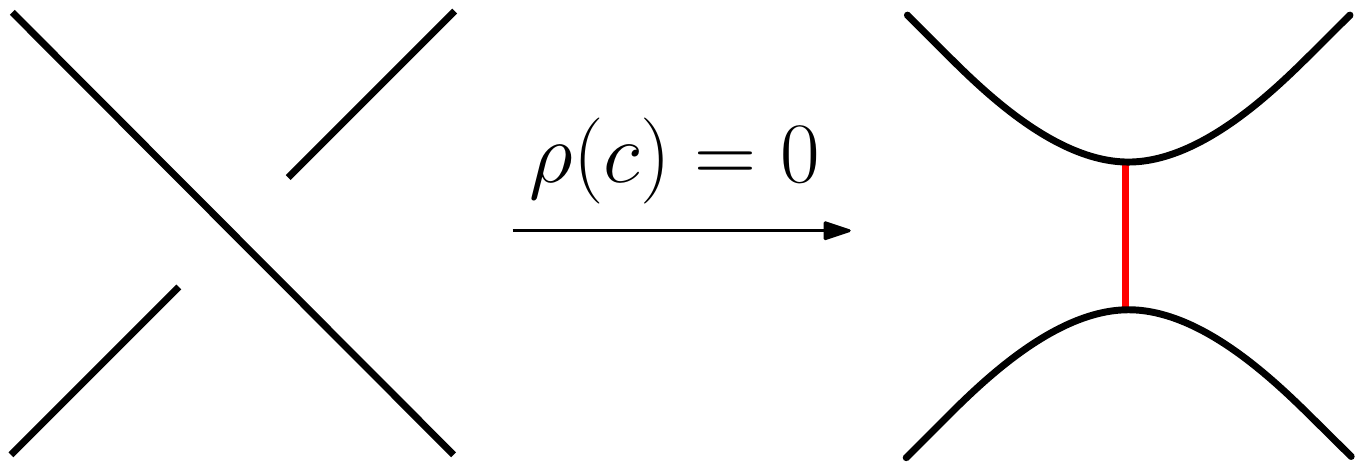} \hspace{.5in} \inlinediag[0.5]{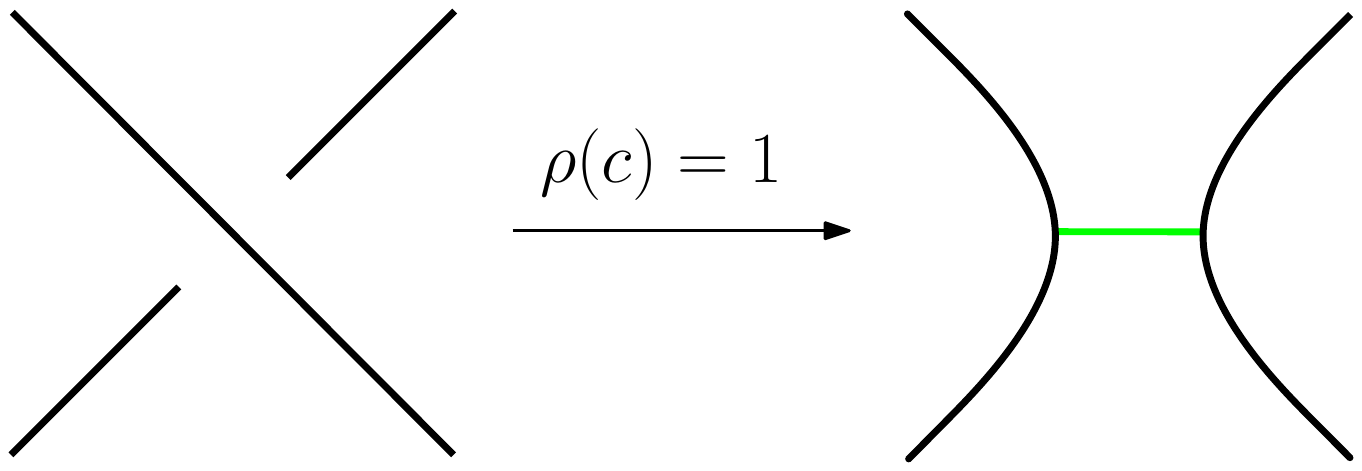} 
$$ 
\noindent The set of resolutions will be denoted $\resolution{\lefty{T}}$. 
\end{defn}

\noindent The local arcs introduced by $\lefty{T} \leftarrow \rho(\lefty{T})$ are called {\em resolution bridges}. The resolution bridge for a crossing $c \in \cross{\lefty{T}}$ will be denoted $\gamma_{r, c}$ (or just $\gamma_{c}$ when the resolution is understood). If $\rho(c) = 0$ we will call $\gamma_{c}$ an {\em active} bridges for $r$, while if $\rho(c)=1$ it will be called {\em inactive}. The active bridges for $r$ are the elements of the set $\actor{r}$. We will denote by $\gamma_{c}(r)$ the resolution obtained by surgering the diagram for $r$ along $\gamma_{c}$. The resolution bridges at $c$ for $\gamma_{c}(r)$ will be denoted by $\gamma_{c}^{\dagger}$ when considered from $r$. \\
\ \\

\noindent A resolution diagram $r(\lefty{T})$ consists of a planar diagram of circles, which we divide into two groups: 1) the {\em free circles} which are contained in $\mathrm{int\,}\leftHalf$ and are the elements of a set $\free{r}$, and 2) the {\em cleaved circles} which cross the $y$-axis, and which determine an element $\mathrm{cl}(r) \in \cleave{n}$. 

\begin{defn}
A {\em state} for $\lefty{T}$ is a pair $(r, s)$ where
\begin{enumerate}
\item $r$ is a resolution of $\lefty{T}$,
\item $s$ is an assignment of an element of $\{+,-\}$ to each circle of $r(\lefty{T})$. This assignment will be called a decoration on $r(\lefty{T})$. 
\end{enumerate}
The states for $\lefty{T}$ will be denoted $\state{\lefty{T}}$.
\end{defn}

\begin{defn}
The boundary of a state $(r, s)$ for $\lefty{T}$ is the element $\partial(r,s) =(\mathrm{cl}(r),\sigma) \in \cleaved{n}$ where $\sigma= s|_{\mathrm{cl}(r)}$. 
\end{defn}

\subsection{A bigraded module spanned by the states}

\begin{defn}
For a state $(r,s) \in \state{\lefty{T}}$ with $r = (\rho, \righty{m})$, let
\begin{enumerate}
\item $h(r) = \sum_{c \in \cross{\lefty{T}}} \rho(c)$
\item $q(r,s) = \sum_{C \in \free{r}} s(C)$
\item $\iota(r,s) = \iota(L,\sigma)$ where $(L,\sigma) = \partial(r,s)$
\end{enumerate}
\end{defn}

\noindent Let $R$ be a ring, and pick $(L,\sigma) \in \cleaved{n}$. Let
$$
\lefty{CK}(\lefty{T}, L, \sigma) \cong \bigoplus_{\partial(r,s) = (L,\sigma)} R\cdot(r,s)
$$
where $(r,s)$ occurs in bigrading $(h(r) - n_{-}, h(r) + q(r,s) + 1/2\iota(r,s)+ n_{+} - 2 n_{-})$. The first entry will be called the {\em homological grading} of the state, while the second is its {\em quantum grading}. 

\begin{defn}
The type $A$ module for an inside tangle $\lefty{T}$ is
$$
\leftComplex{\lefty{T}} = \bigoplus_{(L,\sigma) \in \cleaved{n}} \lefty{CK}(\lefty{T}, L, \sigma)
$$
\end{defn}

\noindent There is a right action of the idempotent algebra $\mathcal{I}_{n} \subset \mathcal{B}\Gamma_{n}$ on $\leftComplex{\lefty{T}}$:
$$
(r,s) \cdot I_{(L,\sigma)}  = \left\{\begin{array}{ll} (r,s) & \partial(r,s) = (L,\sigma) \\ 0 & \mathrm{else} \end{array} \right.
$$
Thus $I_{(L,\sigma)}$ acts non-trivially only on the summand $\lefty{CK}(\lefty{T}, L, \sigma)$.\\
\ \\
\noindent In addition the construction of M. Asaeda, J. Przytycki, and A. Sikora in \cite{APS} endows $\leftComplex{\lefty{T}}$ with a $(1,0)$ differential, $d_{APS}$, described presently. With this differential, $\leftComplex{\lefty{T}}$ becomes a chain complex, and the main result of \cite{APS} implies that the (bigraded) homology of $\leftComplex{\lefty{T}}$ is an isotopy invariant of $\lefty{\mathcal{T}}$, up to (bigraded) isomorphism. 

\subsection{The differential from \cite{APS}}

To define the differential we first {\em order the crossings} of $\lefty{T}$. Then for $(r,s) \in \state{\lefty{T}}$, we define
$$
d_{APS}(r,s) = \sum_{\gamma \in \actor{r}} (-1)^{I(\rho,\gamma)}\,D_{\gamma,\rho}(r,s)
$$
where 1) $r = (\rho, \righty{m})$, 2) $I(\rho,\gamma) = \sum_{c_{\gamma} < c'} \rho(c')$ is the number of $\rho$-inactive crossings which occur after the crossing $c$ corresponding to $\gamma$, and 3) $D_{\gamma,\rho}$ is a map defined at each active arc. This map is prescribed by the following recipe:\\
\ \\
\begin{enumerate}
\item (Khovanov Case i:) Suppose surgery on $\gamma$ merges the {\em free} circles $C_{1}$ and $C_{2}$ in $\rho$ to get a {\em free circle} $C$ in $\gamma(\rho)$. 
\begin{enumerate}
\item if $s(C_{1})=s(C_{2}) =+$, then $D_{\gamma,\rho}(r,s) = (\gamma(r), s')$ where $s'(C) = +$ and $s'(D) = s(D)$ for every circle $D \neq C_{1},C_{2},C$, free or cleaved;
\item if either $s(C_{1})=-$, $s(C_{2}) =+$ or $s(C_{1})=-$, $s(C_{2})=+$, then $D_{\gamma,\rho}(r,s) = (\gamma(r), s'')$ where $s''(C) = -$ and $s''(D) = s(D)$ for every circle $D \neq C_{1},C_{2},C$, free or cleaved.
\item if $s(C_{1})=s(C_{2})=-$ then $D_{\gamma,\rho}(r,s) = 0$
\end{enumerate}
\item (Khovanov Case ii:) Suppose  $\gamma$ has both feet on the same {\em free circle} $C$ in $r$, then surgering $C$ along $\gamma$ produces two new {\em free circles} $C_{1}$ and $C_{2}$ in $\gamma(r)$, and
\begin{enumerate}
\item if $s(C)=+$ then $D_{\gamma,\rho}(r,s) = (\gamma(r),s_{+-}) + (\gamma(r),s_{-+})$ where $s_{+-}(C_{1}) = +$, $s_{-}(C_{2})=-$ and $s_{+-}(D)=s(D)$ for every other circle in $\circles{r}$. $s_{-+}$ is defined similarly, with the roles of $C_{1}$ and $C_{2}$ reversed. 
\item if $s(C)=-$ then $D_{\gamma,\rho}(r,s) = (\gamma(r),s_{--}) $ where $s_{--}(C_{1}) = s_{--}(C_{2})=-$ and $s_{--}(D)=s(D)$ for every other circle in $\circles{r}$.
\end{enumerate}
\item Suppose $\gamma$ has both feet on the same arc $A$ in $\righty{H} \cap r$, then $\gamma(r)$ will have a new {\em free} circle component $C$. Then $D_{\gamma,\rho}(r,s) = (\gamma(r), s')$  where $s'(C) = -$ and $s'(D) = s(D)$ for every other circle in $\circles{\gamma(r)}$.
\item If $\gamma$ has one foot on a cleaved circle $C$ and the other foot on a {\em free} circle $D$ then surgery on $\gamma$ will merge $D$ into $C$, leaving the other circles unchanged. If $s(D)=+$ then $D_{\gamma,\rho}(r,s) = (\gamma(r), s')$ with $s'(C') = s(C')$ for every other circle in $\circles{r}$, including $C$, while if $s(D)=-$ then $D_{\gamma,\rho}(r,s) = 0$.
\item In every case not covered on this list, $D_{\gamma,\rho}(r,s) = 0$. 
\end{enumerate} 

\subsection{Some classes of active bridges}

\noindent For a state $(r,s)$, we can use $s$ to group the bridges in $\actor{r}$ into (overlapping) classes: 
\begin{enumerate}
\item $\interior{r,s}$ is the subset $\actor{r}$ consisting of those $\gamma$ for which $D_{\gamma, \rho}(r,s) \neq 0$. That is
\begin{enumerate}
\item if both feet of $\gamma$ are on elements of $\free{r}$, or
\item one foot of $\gamma$ is on $C \in \mathrm{cl}(r)$ and the other foot is on $C' \in \free{r}$ with $s(C') = +$, or
\item both feet are on the same arc of $C \cap \leftHalf$ for some $C \in \mathrm{cl}(r)$
\end{enumerate}
\item $\dec{r,s}$ is the subset $\actor{r}$ consisting of those $\gamma$ where
\begin{enumerate}
\item both feet are on the same arc of $C \cap \leftHalf$ for some $C \in \mathrm{cl}(r)$ with $s(C) = +$, or
\item one foot of $\gamma$ is on $C \in \mathrm{cl}(r)$ with $s(C) = +$ and the other foot is on $C' \in \free{r}$ with $s(C') = -$
\end{enumerate}
\item $\leftBridges{r}$ is the subset $\actor{r}$ consisting of those $\gamma$ such that either
\begin{enumerate}
\item $\gamma$ has one foot on $C_{1} \in \mathrm{cl}(r)$ and the other on a distinct circle $C_{2} \in \mathrm{cl}(r)$, or
\item $\gamma$ has both feet on some $C \in \mathrm{cl}(r)$, but they are on different arcs of $C \cap \leftHalf$.
\end{enumerate}
\end{enumerate}

\noindent If $r = (\rho, \lefty{m})$ we will let $\bridges{r} = \leftBridges{r}\cup \rightBridges{\lefty{m}}$ and $\rightBridges{r} = \rightBridges{\righty{m}}$. There is a natural map $\bridges{r} \longrightarrow \bridges{\mathrm{cl}(r)}$.

\section{The type $A$-structure for an inside tangle}

\noindent Given a diagram $\lefty{T}$ of an inside tangle $\lefty{\tangle{T}}$, we describe a type $A$-structure on $\leftComplex{\lefty{T}}$ over $\mathcal{B}\Gamma_{n}$. This structure is specified by two {\em bigrading preserving} maps
\begin{align}
m_{1}&: \leftComplex{\lefty{T}} \longrightarrow \leftComplex{\lefty{T}}[(-1,0)]\\
m_{2}&: \leftComplex{\lefty{T}} \otimes_{\mathcal{I}} \mathcal{B}\bridgeGraph{n} \longrightarrow \leftComplex{\lefty{T}}
\end{align}
Let $\xi=(r,s)$ be a generator of $\leftComplex{\lefty{T}}$ with $\partial \xi = (L,\sigma)$, and let $e \in \mathcal{B}\Gamma_{n}$ be a generator. \\
\ \\
\noindent For $m_{1}$ we let $m_{1}(\xi) = d_{APS}(\xi)$, the differential on $\leftComplex{\lefty{T}}$. $d_{APS}$ maps $(r,s)$ in bigrading $(h,q)$ to an element in $(h+1,q)$. This is bigrading preserving into $\leftComplex{\lefty{T}}[(-1,0)]$. \\
\ \\
\noindent To define the action $m_{2}$ we start by describing the action of the generators of $\mathcal{B}\Gamma_{n}$. $m_{2}(\xi \otimes_{\mathcal{I}} e)$ is computed by 
\begin{enumerate}
\item For the idempotents, $m_{2}(\xi \otimes I_{(L,\sigma)}) = \xi \cdot I_{(L,\sigma)}$, the idempotent action defined above. 
\item when $e = \righty{e_{C}}$ for some $C \in \partial \xi$ with $\sigma(C) = +$, then $m_{2}(\xi \otimes_{\mathcal{I}} \righty{e_{C}}) = (r,s_{C})$ where $s_{C}(C) = -$ but equals $s$ on all other circles in $\circles{r}$. 
\item when $e = \lefty{e_{C}}$ for some $C \in \partial \xi$ with $\sigma(C) = +$ then
$$
m_{2}(\xi \otimes_{\mathcal{I}} \lefty{e_{C}}) = \disp{\sum_{\gamma \in \dec{(r,s),C}} (-1)^{I(r,\cross{\gamma})} (r_{\gamma}, s_{\gamma})}
$$
where $\dec{(r,s),C}$ are those active arcs which can change the decoration on $C$, $r_{\gamma}$ is the result of surgery on $\gamma$, and $s_{\gamma}$ is the new decoration with $s'(C) = -$ (and $s'(D) = +$ if a new circle is created). 
\item when $e = e_{\eta,\sigma,\sigma'}$ for some $\eta \in \leftBridges{L}$ then
$$
m_{2}(\xi \otimes_{\mathcal{I}} e_{\eta,\sigma,\sigma'}) = \disp{\sum_{\gamma \in \actor{r},\mathrm{cl}(\gamma) = \eta} (-1)^{I(r,\cross{\gamma})} (r_{\gamma}, s_{\gamma})}
$$
where $r_{\gamma}$ is surgery along $\gamma$ and $s_{\gamma}$ is the decoration on $r_{\gamma}$ which equals $s$ on $\free{r}$ and $\sigma'$ on $\mathrm{cl}(r)$.  
\item when $e = e_{\eta,\sigma,\sigma'}$ for some $\eta \in \rightBridges{L}$ and $r = (\rho, \righty{m})$, let $r' = (\rho, \righty{m_{\eta}})$ and $s'$ equals $\sigma'$ on the cleaved circles and $s$ on $\free{r'}$. Then $m_{2}(\xi \otimes_{\mathcal{I}} e_{\eta,\sigma,\sigma'}) = (r',s')$.    
\item in all other cases $m_{2}(\xi \otimes_{\mathcal{I}} e) = 0$. Note that $(r,s) \otimes_{\mathcal{I}} e_{1} = 0$ unless $\partial(r,s)$ is the source of $e_{1}$ since otherwise $I_{\partial(r,s)} \cdot e_{1} = 0$.
\end{enumerate}

\begin{prop}
$m_{2}$ is bi-grading preserving
\end{prop}

\noindent{\bf Proof:} Following the same order as above:
\begin{enumerate}
\item when $e = \righty{e_{C}}$, if $\xi$ is in bigrading $(h,q)$, then $\xi \otimes_{\mathcal{I}} \righty{e_{C}}$ is in in bigrading $(h,q) + (0,-1) = (h,q-1)$, whereas $m_{2}(\xi \otimes_{\mathcal{I}} \righty{e_{C}})$ is in bigrading $(h, q - 1)$ since we changed a $(0,+1/2)$ cleaved circle to a $(0,-1/2)$ cleaved circle. 
\item when $e = \lefty{e_{C}}$ the bigrading of  $\xi \otimes_{\mathcal{I}} \lefty{e_{C}}$ is $(h,q) + (1,1)$. For $m_{2}(\xi \otimes_{\mathcal{I}} \lefty{e_{C}})$ we consider the bigrading in two cases. If we merge a $-$ free circle, then $q = \overline{q} + 1/2 - 1$ and the bigrading of $(r_{\gamma}, s_{\gamma})$ is $(h,\overline{q}-1/2) + (1,1) = (h+1,q+1)$. If we divide a $+$ cleaved circle then the bigrading change is from $(h,\overline{q} + 1/2) + (1,1) = (h+1, q+ 3/2)$ to $(h,\overline{q} + 1 - 1/2) + (1,1) =(h+1,q+3/2)$. In either case, there is a $(0,0)$ change in bigrading.
\item when $e = e_{\eta,\sigma,\sigma'}$ for some $\eta \in \leftBridges{L}$: if $\eta$ merges two plus circles then $\xi \otimes_{\mathcal{I}} e_{\eta,\sigma,\sigma'}$ has bidgrading $(h,\overline{q}+1/2+1/2) + (1,1/2) = (h+1,\overline{q}+3/2)$, while
$m_{2}(\xi \otimes_{\mathcal{I}} e_{\eta,\sigma,\sigma'})$ has $(h,\overline{q}+1/2) + (1,1)$, since we change the resolution at a crossing. If $\eta$ merges a $+$ and a $-$ we have $(h,\overline{q}+1/2-1/2) + (1,1/2) = (h+1,\overline{q}+1/2)$ before and
$(h,\overline{q}-1/2) + (1,1)$ after. If $\eta$ divides a $+$ circle then we start with $(h,\overline{q}+1/2)+(1,1/2)$ and end with $(h,\overline{q}+1/2-1/2) + (1,1)$, while if $\eta$ divides a $-$ circle we start with $(h,\overline{q}-1/2)+(1,1/2)$ and end with $(h,\overline{q}-1/2-1/2) + (1,1)$. Each of these is a $(0,0)$ change.   
\item when $e = e_{\eta,\sigma,\sigma'}$ for some $\eta \in \rightBridges{L}$ and $r = (\rho, \righty{m})$:  if surgery on $\eta$ merges two $+$ cleaved circles, then the bigrading of $\xi \otimes_{\mathcal{I}} e_{\eta,\sigma,\sigma'}$ is $(h,\overline{q} + 1/2 + 1/2) + (0,-1/2)$, while that of $m_{2}(\xi \otimes_{\mathcal{I}} e_{\eta,\sigma,\sigma'})$ is $(h,\overline{q}+1/2)$ (as there is no crossing change). Likewise for a $+$ and $-$ circle: $(h,\overline{q} + 1/2 - 1/2) + (0,-1/2) \rightarrow (h,\overline{q}-1/2)$, while for a divide of a $+$ circle: $(h,\overline{q} + 1/2) + (0,-1/2) \rightarrow (h,\overline{q}+1/2-1/2)$, and a divide of a $-$ circle: $(h,\overline{q} - 1/2) + (0,-1/2) \rightarrow (h,\overline{q}-1/2-1/2)$. In all cases there is a $(0,0)$ change in bigrading.  
\end{enumerate}

\noindent This specifies $m_{2}$ on the generators of $\mathcal{B}\Gamma_{n}$. To define $m_{2}$ on all elements we impose the following relation. If $p_{1}, p_{2} \in \mathcal{B}\Gamma_{n}$ we define 
$$
\widetilde{m}_{2}(\xi \otimes p_{1}p_{2}) = \widetilde{m}_{2}(\widetilde{m}_{2}(\xi \otimes p_{1}) \otimes p_{2}) 
$$
with $\widetilde{m}_{2}$ equal to $m_{2}$, as defined above,  on the idempotents and generators. 

\begin{prop}
If two products of the generators $p_{1}$ and $p_{2}$ define  equal elements in $\widetilde{B}\Gamma_{n}$, then
$\widetilde{m}_{2}(\xi \otimes p_{1}) = \widetilde{m}_{2}(\xi \otimes p_{2})$ for every $\xi \in \leftComplex{\lefty{T}}$.   
\end{prop}

\noindent Thus the rules above fully specify $m_{2}: \leftComplex{\lefty{T}} \otimes_{\mathcal{I}} \mathcal{B}\bridgeGraph{n} \longrightarrow \leftComplex{\lefty{T}}$.\\
\ \\
\begin{proof}
It suffices to show that $m_{2}(\xi \otimes \rho) = 0$ whenever $\rho$ is a relation defining $\mathcal{B}\bridgeGraph{n}$. We start with the relations from disjoint supports. Suppose $C$ and $D$ are distinct cleaved circles with $\sigma(C) = \sigma(D) = +$. Then 
\begin{equation}
\begin{split}
\widetilde{m}_{2}(\xi \otimes \big(\righty{e_{C}}\righty{e_{D}} - \righty{e_{D}}\righty{e_{C}}\big)) &= \widetilde{m}_{2}(\widetilde{m}_{2}(\xi \otimes \righty{e_{C}}) \otimes \righty{e_{D}}) -
\widetilde{m}_{2}(\widetilde{m}_{2}(\xi \otimes \righty{e_{D}}) \otimes \righty{e_{C}})\\
 &= (r,s_{C,D}) - (r,s_{C,D}) = 0 \\
\end{split}
\end{equation}
\end{proof}
On the other hand: 
\begin{equation}
\begin{split}
\widetilde{m}_{2}(\xi \otimes \big(\lefty{e_{C}}\righty{e_{D}} - \righty{e_{D}}\lefty{e_{C}}\big)) &= 
\widetilde{m}_{2}(\widetilde{m}_{2}(\xi \otimes \lefty{e_{C}}) \otimes \righty{e_{D}}) -
\widetilde{m}_{2}(\widetilde{m}_{2}(\xi \otimes \righty{e_{D}}) \otimes \lefty{e_{C}})\\
&= \disp{\sum_{\gamma \in \dec{(r,s),C}} (-1)^{I(r,\cross{\gamma})} (r_{\gamma}, s_{\gamma,D})} -
\disp{\sum_{\gamma \in \dec{(r,s),C}} (-1)^{I(r,\cross{\gamma})} (r_{\gamma}, s_{D,\gamma})} \\ & = 0
\end{split}
\end{equation}
To compute $\widetilde{m}_{2}(\xi \otimes \big(\lefty{e_{C}}\lefty{e_{D}} + \lefty{e_{D}}\lefty{e_{C}}\big))$ note
that each $\widetilde{m}_{2}(\widetilde{m}_{2}(\xi \otimes \lefty{e_{C}}) \otimes \lefty{e_{D}})$  and
$\widetilde{m}_{2}(\widetilde{m}_{2}(\xi \otimes \lefty{e_{D}}) \otimes \lefty{e_{C}})$ are sums over pairs of edges
$\gamma \in \dec{r,s,C}$ and $\gamma' \in \dec{r,s,D}$. In one case we sum over $(\gamma,\gamma')$ pairs and in the other $(\gamma', \gamma)$ pairs. In each case we obtain $(r_{\gamma,\gamma'}, s_{\gamma,\gamma'})$ with $s_{\gamma,\gamma'}$ uniquely determined by the requirement that $-$'s decorate $C$ and $D$. Thus we need only look at the signs: for $(\gamma,\gamma')$ we have $(-1)^{I(r,\cross{\gamma})+ I(r_{\gamma}, \cross{\gamma'})}$ which is $-(-1)^{I(r,\cross{\gamma'})+ I(r_{\gamma'}, \cross{\gamma})}$. Consequently, they cancel in the sum. \\
\ \\
\noindent Now suppose that $C_{1}$ and $C_{2}$ are cleaved circles in $r$ with $s(C_{1}) = s(C_{2})= +$. Let $\beta$
be an active arc which merges $C_{1}$ and $C_{2}$ to get $C$ and maps to $\gamma \in \leftBridges{L}$. We can partition the active arcs $\alpha$ which contribute to $\dec{r_{\beta},s_{\beta},C}$ into the three sets: $\dec{r,s,C_{1}}$, $\dec{r,s,C_{2}}$, and $\alpha$ which also map to $\gamma$. To obtain
$\widetilde{m}_{2}(\xi \otimes m_{\gamma}\lefty{e_{C}})$ we sum over all such $\beta$ and $\alpha$ arcs : $\sum_{(\beta,\alpha)}(-1)^{I(\beta) + I(r_{\beta},\alpha)}(r_{\beta,\alpha}, s_{\beta\alpha})$. For $\alpha$ isotopic as a bridge to $\beta$ the term for $(\alpha,\beta)$ occurs in this sum, with sign  $(-1)^{I(\alpha) + I(r_{\alpha},\beta)}$. This cancels the term from $(\beta,\alpha)$. Thus 
\begin{equation}
\begin{split}
\widetilde{m}_{2}(\xi \otimes m_{\gamma}\lefty{e_{C}}) &= \sum_{\beta,\alpha \in \dec{r,s,C_{1}}}(-1)^{I(\beta) + I(r_{\beta},\alpha)}(r_{\beta,\alpha}, s_{\beta\alpha}) + \sum_{\beta,\alpha \dec{r,s,C_{2}}}(-1)^{I(\beta) + I(r_{\beta},\alpha)}(r_{\beta,\alpha}, s_{\beta\alpha})\\
&= -\sum_{\alpha \in \dec{r,s,C_{1}}, \beta}(-1)^{I(\alpha) + I(r_{\alpha},\beta)}(r_{\alpha,\beta}, s_{\alpha\beta}) - \sum_{\alpha \in \dec{r,s,C_{2}}, \beta}(-1)^{I(\alpha) + I(r_{\alpha},\beta)}(r_{\alpha,\beta}, s_{\alpha\beta})\\
&= - \widetilde{m}_{2}(\xi \otimes \lefty{e_{C_{1}}}m_{\gamma}) - \widetilde{m}_{2}(\xi \otimes \lefty{e_{C_{2}}}m_{\gamma})
\end{split}
\end{equation} 
which verifies that $\widetilde{m}_{2}$ is compatible with relation \ref{rel:lefty3}. Exactly the same argument applies to
for $\gamma \in \rightBridges{L}$, although we no longer need to sum over the representatives of $\gamma$ since there is only one such bridge. More significantly, all the terms occur with sign $(-1)^{I(\alpha)}$ since surgery on $\gamma$ does not affect the signs. The conclusion becomes
$$
\widetilde{m}_{2}(\xi \otimes m_{\gamma}\lefty{e_{C}}) =  \widetilde{m}_{2}(\xi \otimes \lefty{e_{C_{1}}}m_{\gamma}) + \widetilde{m}_{2}(\xi \otimes \lefty{e_{C_{2}}}m_{\gamma})
$$
which is compatible with relation \ref{rel:lefty1}. The case where surgery on $\gamma$ is divisive follows from the same line of reasoning.\\
\ \\
\noindent For $\righty{e_{C}}$ and $\gamma$ merging $C_{1}$ and $C_{2}$ the situation is easier. First, suppose $\gamma \in \rightBridges{L}$. Then
$$
\widetilde{m}_{2}(\xi \otimes m_{\gamma}\righty{e_{C}}) = (r, s_{\gamma,C})
$$
while
$$
\widetilde{m}_{2}(\xi \otimes \righty{e_{C_{1}}}m_{\gamma}) = \widetilde{m}_{2}((r,s_{C_{1}}) \otimes m_{\gamma}) = (r,s_{C,\gamma})
$$
As these are equal, and as $C_{2}$ plays a symmetric role, $\widetilde{m}_{2}$ is compatible with this type of relation. Again, the case for dividing is similar. For $\righty{e_{C}}$ and $\gamma \in \leftBridges{L}$, we have
$$
\widetilde{m}_{2}(\xi \otimes m_{\gamma}\righty{e_{C}}) = \sum_{\alpha}(-1)^{I(\alpha)}(r_{\alpha}, s_{\alpha,C})
$$
where the sum is over active arcs which map to $\gamma$. On the other hand,
while
$$
\widetilde{m}_{2}(\xi \otimes \righty{e_{C_{1}}}m_{\gamma}) = \widetilde{m}_{2}((r,s_{C_{1}}) \otimes m_{\gamma}) = \sum_{\alpha}(-1)^{I(\alpha)}(r_{\alpha},s_{C,\alpha})
$$
Thus $\widetilde{m}_{2}$ is compatible with $m_{\gamma}\righty{e_{C}} = \righty{e_{C_{1}}}m_{\gamma}$ for all $\gamma \in \bridges{L}$. \\
\ \\
\noindent Suppose that $\gamma$ and $\gamma'$ are in $\bridges{L}$ and that there is a commuting square
$$
\begin{CD}
(L,\sigma) @>e_{(\gamma,\sigma,\sigma_{01})}>> (L_{\gamma},\sigma_{01}) \\
@Ve_{(\gamma',\sigma,\sigma_{10})}VV 				@VVe_{(\gamma',\sigma_{01},\sigma'')}V\\
(L_{\gamma'}, \sigma_{10}) @>e_{(\gamma,\sigma_{01},\sigma'')}>> (L_{\gamma,\gamma'}, \sigma'')\\
\end{CD}$$
Then if
\begin{enumerate}
\item Both $\gamma$ and $\gamma'$ are in $\rightBridges{L}$ we need to see $\widetilde{m}_{2}(\xi \otimes \big(e_{\gamma}e_{\gamma'} - e_{\gamma'}e_{\gamma}\big)) = 0$. However, both terms resulting from expanding $\widetilde{m}$ will
equal $(r_{\gamma\gamma'}, s'')$ where $r_{\gamma,\gamma'}$ is identical to $r$ in $\leftHalf$ but equals $L_{\gamma\gamma'}$ in $\rightHalf$, and $s''$ is $s$ on $\free{r}$ but $\sigma''$ on $\mathrm{cl}(r)$. Since both
terms are identical, the difference will be zero and $\widetilde{m}_{2}$ is compatible with this case.
\item If $\gamma$ in $\leftBridges{L}$ but $\gamma' \in \rightBridges{L}$, then we need $\widetilde{m}_{2}(\xi \otimes \big(e_{\gamma}e_{\gamma'} - e_{\gamma'}e_{\gamma}\big)) = 0$. The action of $e_{\gamma}$ followed by $e_{\gamma'}$ (or vice-versa) will give $\sum_{\alpha} (-1)^{I(r,\alpha)}(r_{\alpha,\gamma'}, s'')$ where the sum is over all active arcs for $r$ which have image $\gamma$ in $\mathrm{cl}(r)$. Since surgery on $\gamma'$ does not affect the sign, we see that the two terms will cancel, and $\widetilde{m}_{2}$ is compatible with this case.
\item Suppose both $\gamma$ and $\gamma'$ are in $\leftBridges{L}$, then $\widetilde{m}_{2}(\xi \otimes e_{\gamma}e_{\gamma'})$ is the sum over pairs of active arcs $(\alpha,\alpha')$ for $r$ which map to $\gamma$ and $\gamma'$ when considered as bridges. Each pair also contributes to $\widetilde{m}_{2}(\xi \otimes e_{\gamma'}e_{\gamma})$ but in the reversed order $(\alpha',\alpha)$. The decorations of the result are determined by $\sigma''$, so we need only check the signs of each term. The sign for $(\alpha,\alpha')$ is $(-1)^{I(r,\cross{\alpha})+ I(r_{\alpha}, \cross{\alpha'})}$ while that for $(\alpha',\alpha)$ is $(-1)^{I(r,\cross{\alpha'})+ I(r_{\alpha'}, \cross{\alpha})}$. Due to the ordering of the crossings, one of these will be $+1$ and the other $-1$. Consequently, $\widetilde{m}_{2}(\xi \otimes e_{\gamma}e_{\gamma'}) = -\widetilde{m}_{2}(\xi \otimes e_{\gamma'}e_{\gamma})$ which is compatible with the relation for $\leftBridges{L}$.
\end{enumerate}

\noindent Note that a similar argument to that in (1) works for all pairs of right bridge edges that form commutative squares, so it will be omitted. For left bridge elements, there are two further cases to consider:\\
\ \\
\noindent A) Suppose $\delta \in \lefty{B}_{o}(L,\gamma)$ and $\delta_{1}, \delta_{2} \in B_{s}(L_{\gamma}, \gamma^{\dagger})$ map to $\delta$ under surgery along $\gamma^{\dagger}$. Likewise, suppose $\gamma_{1}$ and $\gamma_{2}$ map to $\gamma$ under surgery along $\delta^{\dagger}$. If we orient the mutual arc for $\delta$ and $\gamma$, then $\delta_{1}$ corresponds to the version of $\delta$ before $\gamma$ along the mutual arc, and $\delta_{2}$ corresponds to that after $\gamma$. Likewise for $\gamma_{i}, i=1,2$ and $\delta$. There can then be (anti-)commutative squares
$$
\begin{CD}
(L,\sigma) @>e_{(\gamma,\sigma,\sigma_{01})}>> (L_{\gamma},\sigma_{01}) \\
@Ve_{(\delta,\sigma,\sigma_{10})}VV 				@VVe_{(\delta_{1},\sigma_{01},\sigma'')}V\\
(L_{\delta}, \sigma_{10}) @>e_{(\gamma_{2},\sigma_{01},\sigma'')}>> (L_{\gamma,\delta_{1}}, \sigma'')\\
\end{CD}
 \hspace{1in}
\begin{CD}
(L,\sigma) @>e_{(\gamma,\sigma,\sigma_{01})}>> (L_{\gamma},\sigma_{01}) \\
@Ve_{(\delta,\sigma,\sigma_{10})}VV 				@VVe_{(\delta_{2},\sigma_{01},\sigma'')}V\\
(L_{\delta}, \sigma_{10}) @>e_{(\gamma_{1},\sigma_{01},\sigma'')}>> (L_{\gamma,\delta_{2}}, \sigma'')\\
\end{CD}
$$ 
where $L_{\gamma,\delta_{2}}$ and $L_{\gamma,\delta_{1}}$ have different left matchings. Two resolution arcs $a_{1}$ and $a_{2}$ for $\xi$, one corresponding to $\gamma$ and one corresponding to $\delta$, will be counted in the action of either $\gamma$ followed by $\delta_{1}$, or $\gamma$ followed by $\delta_{2}$ (but not both). Reversing the order means contributing to $\delta$ followed by $\gamma_{2}$, or $\delta$ followed by $\gamma_{1}$. As these contibute with the usual sign conventions the contributions of the pair will cancel in either the action of $e_{\gamma}e_{\delta_{1}} + e_{\delta}e_{\gamma_{2}}$ or $e_{\gamma}e_{\delta_{2}} + e_{\delta}e_{\gamma_{1}}$, which verifies that the action respects the anti-commutativity in this case. \\
\ \\
\noindent B) Now suppose that $\lefty{\alpha} \in B_{s}(L,\lefty{\gamma})$ and that $\lefty{\beta}$ is the bridge obtained by sliding $\lefty{\alpha}$ over $\lefty{\gamma}$. Let $\lefty{\delta}$ be the image of $\lefty{\alpha}$ and $\lefty{\beta}$ in $L_{\lefty{\gamma}}$, $\lefty{\zeta}$ be the image of $\lefty{\alpha}$ and $\lefty{\gamma}$ in $L_{\lefty{\beta}}$ and $\lefty{\eta}$ be the image of $\lefty{\beta}$ and $\lefty{\gamma}$ in $L_{\lefty{\alpha}}$. The action of $\lefty{e_{\alpha}}\lefty{e_{\eta}}$ is the sum over active resolution arcs $a_{1}$ and $a_{2}$ for $\xi$ with $a_{1}$ representing $\lefty{\alpha}$ and $a_{2}$ representing $\lefty{\eta}$ in $L_{\alpha}$. $a_{2}$ thus represents one of $\beta$ or $\gamma$ in $L$. Reversing the order thus gives a contribution to either $\lefty{e_{\beta}}\lefty{e_{\zeta}}$ or $\lefty{e_{\gamma}}\lefty{e_{\delta}}$. However, pairs representing $\beta$ and $\gamma$ also contribute to $\lefty{e_{\beta}}\lefty{e_{\zeta}}$, while their reverse contributes to $\lefty{e_{\gamma}}\lefty{e_{\delta}}$. Thus, if we consider all ordered pairs $(a_{1},a_{2})$ of resolution arcs which represent pairs of $\alpha$, $\beta$, or $\gamma$, in either order, we will count both $(a_{1},a_{2})$ and $(a_{2},a_{1})$ for each such pair, and they will contribute with opposite signs (due to the Khovanov sign conventions) in the action of $\lefty{e_{\alpha}}\lefty{e_{\eta}} + \lefty{e_{\beta}}\lefty{e_{\zeta}} + \lefty{e_{\gamma}}\lefty{e_{\delta}}$ on $\xi$. Consequently, all the terms in this action will cancel, verifying that it acts as $0$.\\
\ \\
\noindent For the additional bridge relations, suppose $\gamma \in \rightBridges{L}$ then there is a relation $e_{\gamma,\sigma,\sigma'}e_{\gamma^{\dagger}\sigma',\sigma_{C}}= \righty{e_{C}}$. In this case
\begin{equation}
\begin{split}
\widetilde{m}_{2}(\xi \otimes \big(e_{\gamma,\sigma,\sigma'}e_{\gamma^{\dagger}\sigma',\sigma_{C}} - \righty{e_{C}}\big)) &=
((\rho, \righty{m}_{\gamma\gamma^{\dagger}}), s_{C}) - (r, s_{C})\\
&= ((\rho, \righty{m}), s_{C}) - (r, s_{C}) = (r,s_{C}) - (r,s_{C}) = 0\\
\end{split}
\end{equation}
where $\righty{m}_{\gamma\gamma^{\dagger}} = \righty{m}$ follows from the result that surgery on a bridge $\gamma$ for $L$, followed by surgery on $\gamma'$, recovers $L$. \\
\ \\
\noindent Now suppose that $\gamma \in \leftBridges{L}$ and $\eta \in \leftBridges{L_{\gamma}}$ intersects $\gamma^{\dagger}$ non-trivially. We need to see that $\widetilde{m}_{2}(\xi \otimes e_{\gamma}e_{\eta}) = 0$ since $e_{\gamma}e_{\eta}=0$. However,
in $\widetilde{m}_{2}(\widetilde{m}_{2}(\xi \otimes e_{\gamma}) \otimes e_{\eta})$ the action of $e_{\eta}$ will result in a sum over active arcs in $r_{\gamma}$ which map to $\eta$ in $\mathrm{cl}(r)$. Each active arc comes from a crossing, and thus must already be present in $r$ for it to be present in $r_{\gamma}$. This excludes there being any active arc for $\eta$ in $r_{\gamma}$. Consequently the sum is $0$ and we have verified that $\widetilde{m}_{2}$ is compatible with this relation.

\begin{prop}
For $\xi=(r,s)$ a generator of $\leftComplex{\lefty{T}}$ and $\rho_{1},\rho_{2} \in \mathcal{B}\bridgeGraph{n}$. The maps $m_{1}$ and $m_{2}$ above satisfy:
\begin{align}
0 &= m_{1}(m_{1}(\xi)) \\
0 &= (-1)^{\lefty{l}(\rho_{1})}m_{2}(m_{1}(\xi) \otimes \rho_{1}) +  m_{2}(\xi \otimes \mu_{\gamma}(\rho_{1})) - m_{1}(m_{2}(\xi \otimes \rho_{1}))\\
0 &= m_{2}(m_{2}(\xi \otimes \rho_{1}) \otimes \rho_{2}) - m_{2}(\xi \otimes \rho_{1}\rho_{2})
\end{align}
\end{prop}

\noindent{\bf Note:} These are the relations for $\leftComplex{\lefty{T}}$ to be an $A_{\infty}$-module over the differential graded algebra $\mathcal{B}\bridgeGraph{n}$, as in \cite{Bor1}, with $m_{n} = 0$ for $n \geq 3$.\\
\ \\
\begin{proof}
That $m_{1}(m_{1}(\xi)) = 0$ is a byproduct of $m_{1} = d$ being a differential (see also the proof that $\righty{\delta}$  is a $D$-structure for an outside tangle). Furthermore, that $m_{2}$ defines a right action follows from defining $\widetilde{m}_{2}$ to be a right action, which descends to $m_{2}$ after we see $m_{2}$ is well-defined. Thus, we need only verify that $d$ and $m_{2}$ are compatible with $\mu_{\Gamma}$ through the equation 
$$
d(m_{2}(\xi \otimes \rho_{1})) = (-1)^{\lefty{l}(\rho_{1})}m_{2}(d(\xi) \otimes \rho_{1}) +  m_{2}(\xi \otimes \mu_{\gamma}(\rho_{1}))
$$
It suffices to prove this for $\rho_{1}$ of length $0$ or $1$ since we can bootstrap the relation for longer words using 
\begin{equation}
\begin{split}
d(\xi \cdot(\alpha\beta)) &= d((\xi \cdot \alpha) \cdot \beta) \\
&= (-1)^{\lefty{l}(\beta)}d(\xi \cdot \alpha)\cdot \beta + (\xi\cdot \alpha)\cdot \mu_{\gamma}(\beta)\\
&= (-1)^{\lefty{l}(\beta)+\lefty{l}(\alpha)}\big(d(\xi) \cdot \alpha\big)\cdot \beta + (-1)^{\lefty{l}(\beta)}\big(\xi \cdot \mu_{\gamma}(\alpha)\big)\cdot \beta + \xi \cdot\big(\alpha \mu_{\gamma}(\beta) \big)\\
&= (-1)^{\lefty{l}(\beta)+\lefty{l}(\alpha)}\big[d(\xi) \cdot (\alpha\beta)\big] + \xi \cdot \big[(-1)^{\lefty{l}(\beta)}\mu_{\gamma}(\alpha)\beta +  \alpha \mu_{\gamma}(\beta)\big)\big]\\
&= (-1)^{\lefty{l}(\beta)+\lefty{l}(\alpha)}\big[d(\xi) \cdot (\alpha\beta) + \xi \cdot \mu_{\gamma}(\alpha\beta)\big]
\end{split}
\end{equation}
For length $0$ we have $\rho_{1} = I_{(L,\sigma)}$ for some idempotent. If $\partial \xi \neq (L,\sigma)$ then both sides are zero since 1) $m_{2}(\xi \otimes I_{(L,\sigma)}) = m_{2}(0) = 0$, 2) $d(\xi)$ has the same boundary as $\xi$ so $d(\xi) \otimes I_{(L,\sigma)} = 0$, and 3) $\mu_{\gamma}(I_{(L,\sigma)}) = 0$ for every idempotent. On the other hand, if $\partial \xi = (L,\sigma)$ the last term vanishes, and
$$
d(m_{2}(\xi \otimes I_{(L,\sigma)})) = d(\xi) = m_{2}(d(\xi) \otimes I_{(L,\sigma)})
$$
For length one words, we need to check when $\rho_{1} = \righty{e_{C}}, \lefty{e_{C}},$ or  $e_{\gamma}$ for $\gamma \in \bridges{L}$ where $\partial \xi = (L,\sigma)$. \\
\ \\
\noindent We know $\mu_{\gamma}(\righty{e_{C}})=0$ and $\lefty{l}(\righty{e_{C}})=0$, so for $\righty{e_{C}}$ we need only verify that $d(m_{2}(\xi \otimes \righty{e_{C}})) = m_{2}(d(\xi) \otimes \righty{e_{C}})$. If $\xi$ has $\sigma(C) = -$ then both are $0$, whereas if $\sigma(C) = +$ then
both equal $\Sigma_{\alpha} (-1)^{I(\alpha)}(r_{\alpha}, s_{\alpha,C})$ where the sum is over all active, non-bridging arcs $\alpha$ and $s_{\alpha,C}$ is any decoration compatible with $d$, $r_{\alpha}$, and assigning $C$ a $-$. For $e_{\gamma}$ with $\gamma \in \rightBridges{L}$, the only difference is that the sum is over terms $(r_{\alpha, \gamma}, s_{\alpha,\gamma})$.\\
\ \\
\noindent For $e_{\gamma}$ with $\gamma \in \leftBridges{L}$, $\mu_{\Gamma}$ still vanishes but $\lefty{l}=1$. We then have
\begin{equation}
d(m_{2}(\xi \otimes e_{\gamma})) = \sum_{\alpha,\beta}(-1)^{I(r,\alpha) + I(r_{\alpha},\beta)} (r_{\alpha,\beta}, s_{\alpha,\beta})
\end{equation} 
where the sum is over all active arcs $\alpha$ which map to $\gamma$ and all active arcs $\beta$ which contribute to $d$ (as well as all compatible decorations on $r_{\alpha,\beta}$. On the other hand,
$$
m_{2}(d(\xi)\otimes e_{\gamma}) = \sum_{\beta,\alpha}(-1)^{I(r,\beta) + I(r_{\beta},\alpha)} (r_{\beta,\alpha}, s_{\beta,\alpha})
$$
due to the ordering of the crossings the signs will be different for each $(\alpha,\beta)$ term, so
$$
d(m_{2}(\xi \otimes e_{\gamma}))= - m_{2}(d(\xi)\otimes e_{\gamma}) = (-1)^{\lefty{l}(e_{\gamma})}m_{2}(d(\xi)\otimes e_{\gamma})
$$
\ \\
\noindent We are left with verifying the formula for $\lefty{e_{C}}$. We start with
$$
d(m_{2}(\xi \otimes \lefty{e_{C}})) = \sum_{\alpha,\beta} (-1)^{I(r,\alpha) + I(r_{\alpha},\beta)}(r_{\alpha,\beta}, s_{\alpha,\beta})
$$
where the sum is over all active arcs $\alpha \in \dec{r,s,C}$ and $\beta$ contributing to $d$ on $r_{\alpha}$. Furthermore,
$$
m_{2}(d(\xi)\otimes \lefty{e_{C}}) = \sum_{\beta',\alpha'}(-1)^{I(r,\beta') + I(r_{\beta'},\alpha')} (r_{\beta',\alpha'}, s_{\beta',\alpha'})
$$ 
where the sum is over all $\beta'$ contributing to $d$ on $r$ and all $\alpha$ contributing to $\dec{r_{\beta'},s_{\beta'},C}$.
For pairs $(\alpha,\beta)$ and $(\beta,\alpha)$ occurring in both sums, the coefficient of one is minus the coefficient of the other. However, there are also terms that do not cancel. These correspond to $(\alpha,\beta)$ which become bridges when reversed, and correspond to a $\gamma,\gamma^{\dagger}$ pair with $C$ as its active circle. Due to the reversal, these will be counted with opposite signs from the count above. Let $R$ be the sum over reversible pairs, $\Psi_{1}$ be the part of $d(m_{2})$ which comes from pairs that reverse to bridges, and $\Psi_{2}$ be the part that comes from pairs in $m_{2}(d)$ which reverse to 
bridges. If we sum of the $\gamma\gamma^{\dagger}$ pairs we will get $-\Psi_{1}$ and $-\Psi_{2}$. So,
$$
d(m_{2}(\xi \otimes \lefty{e_{C}})) = R + \Psi_{1} = -(-R + \Psi_{2}) + \Psi_{2} + \Psi_{1} = - m_{2}(d(\xi)\otimes \lefty{e_{C}}) - (-\Psi_{1} - \Psi_{2})
$$
where $-\Psi_{1} - \Psi_{2} = m_{2}(\xi \otimes \sum_{\gamma} e_{\gamma}e_{\gamma^{\dagger}})$ where the sum is over all
$\gamma$ with active circle $C$. As this sum is just the action of $-\mu_{\Gamma}(\lefty{e_{C}})$ we obtain the relation
\begin{equation}
\begin{split}
d(m_{2}(\xi \otimes \lefty{e_{C}})) &=  - m_{2}(d(\xi)\otimes \lefty{e_{C}}) - m_{2}(\xi \otimes -\mu_{\Gamma}(\lefty{e_{C}})) \\
&=(-1)^{\lefty{l}(\lefty{e_{C}})}\big(m_{2}(d(\xi)\otimes \lefty{e_{C}})\big) + m_{2}(\xi \otimes \mu_{\Gamma}(\lefty{e_{C}}))
\end{split}
\end{equation}
as required. \\
\ \\
\noindent We have now verified that the action of the length one words is compatible with the (right) Leibniz relation, and thus using the bootstrap, that $d$ is a (right) differential on the module $\leftComplex{\lefty{T}}$.
\end{proof}

\section{Simplifying type $A$-structures and Reidemeister Invariance}\label{sec:invariance}

\subsection{Algebra}

\noindent In this section we redefine $A_{\infty}$-algebras and module to be  consistent with our sign conventions. We begin with some notation for handling signs and gradings

\begin{defn}
Let $W = W_{0} \oplus W_{1}$ be a $\Zmod{2}$-graded module.  $|\I_{W}| : W \rightarrow W$ is the signed identity defined by
linearly extending
$$
|\I_{W}|(w) = (-1)^{\mathrm{gr}(w)}w
$$
for homogeneous $w \in A$.  
\end{defn}

\noindent{\bf Note:} By $|\I|^{j}$ we mean the composition of $|\I|$ with itself $j$ times. Furthermore, by $|\I|^{j \otimes n}$ we will mean the function $|\I|^{j}\otimes \cdots \otimes |\I|^{j}$, where there are $n$ factors. For an element $\alpha$, $|\I|^{j}(\alpha)$ will be shortened to $|\alpha|^{j}$. Thus, on a homogeneous element $\alpha$, $|\alpha|^{j} = (-1)^{j\mathrm{gr}(\alpha)}\alpha$ , and $\big||\alpha|^{j}\big|^{k} = |\alpha|^{j+k}$.   \\
\ \\

\begin{defn}
An $A_{\infty}$-algebra $A$ over a ring $R$ is a graded module $A$ equipped with maps $\mu_{n}:A^{\otimes n} \rightarrow A[n-2]$ for each $n \in \N$ which
satisfy the relation
$$
0 = \sum_{\footnotesize\begin{array}{c} i\!+\!j\!=\!n\!+\!1 \\ k\!\in\!\{1,\!\ldots\!,\!n\!-\!j\!+\!1\!\} \end{array}\normalsize} (-1)^{j(i+1)+(k+1)(j+1)}\mu_{i}(\I^{\otimes (k-1)} \otimes \mu_{j} \otimes |\I|^{j \otimes (n-k-j+1)} )  
$$
\end{defn} 

\begin{defn}
A right module over a $\Zmod{2}$-graded differential $R$-algebra $(A, \mu_{1}, \mu_{2})$ is an $R$-module $M$ together
with maps $m_{1}: M \rightarrow M[-1]$ and $m_{2}: M \otimes_{R} A \rightarrow M$ such that
\begin{align}
0 &= m_{1} \circ m_{1}  \\
0 &= m_{2}(m_{1} \otimes |\I|) +  m_{2}(\I \otimes \mu_{1}) - m_{1}(m_{2})\\
0 &= m_{2}(m_{2} \otimes \I) - m_{2}(\I \otimes \mu_{2})
\end{align}
\end{defn}

\noindent A right module as above is a special case of the $A_{\infty}$-modules found in \cite{Bor1} (defined using the sign conventions in this paper):

\begin{defn}[\cite{Bor1}]
A right $A_{\infty}$-module $M$ over an $A_{\infty}$-algebra $A$ is a set of maps $\{m_{i}\}_{i\in \N}$ with
$m_{i}: M \otimes A^{\otimes (i-1)} \rightarrow M[i-2]$, and satisfying the following relations for each $n \geq 1$:
\begin{equation}
\begin{split}
0 = \sum_{i+j=n+1} (-1)^{j(i+1)}m_{i}(m_{j}& \otimes |\I|^{j \otimes (i-1))})\\
& + \sum_{i+j=n+1,k>0}(-1)^{k(j+1)+j(i+1)} m_{i}(\I^{\otimes\,k} \otimes \mu_{j} \otimes |\I|^{j \otimes (i-k-1)})
\end{split}
\end{equation}
$M$ is said to be strictly unital if for any $\xi \in M$,  $m_{2}(\xi \otimes 1_{A}) = \xi$, but for $n > 1$, $m_{n}(\xi \otimes a_{1} \otimes a_{2} \otimes \cdots \otimes a_{n-1}) = 0$ if any $a_{i} = 1_{A}$.
\end{defn} 

\noindent Our right modules correspond to $m_{i} = 0$ for $i \geq 2$. Nevertheless, we will think of these as objects in the category of right $A_{\infty}$-modules. The morphisms in this category, again ignoring signs, are

\begin{defn}[\cite{Bor1}]
An $A_{\infty}$-morphism $\Psi$ of right $A$-modules $M$ and $M'$ is a set of maps
$\psi_{i}: M \otimes A^{\otimes\,(i-1)} \longrightarrow M'[i-1]$ for $i \in \N$,  satisfying
\begin{equation}
\begin{split}
\sum_{i+j=n+1}(-1)^{(i+1)(j+1)}m'_{i}(\psi_{j}&\otimes |\I|^{(j+1)\otimes(i-1)}) =\\
&= \sum_{i+j=n+1}(-1)^{j(i+1)}\psi_{i}(m_{j} \otimes |\I|^{j \otimes(i-1)})\\
& \hspace{.5in} + \sum_{i+j=n+1,k>0} (-1)^{j(i+1)+k(j+1)}\psi_{i}(\I^{\otimes\,k} \otimes \mu_{j} \otimes |\I|^{j \otimes (i-k-1)})
\end{split}
\end{equation} 
$\Psi$ is {\em strictly unital} if $\psi_{i}(\xi \otimes a_{1} \otimes \cdots \otimes a_{i-1}) = 0$ when $a_{j} = 1_{A}$ for some $j$ and $i > 1$. The {\em identity} morphism $I_{M}$ is the collection of maps $i_{1}(\xi) = \xi$, $i_{j} = 0$ for $j > 1$
\end{defn}

\begin{defn}[\cite{Bor1}]
Let $\Psi$ be an $A_{\infty}$-morphism from $M$ to $M'$, and let $\Phi$ be an $A_{\infty}$-morphism from $M'$ to $M''$. The composition $\Phi \ast \Psi$ is the morphism whose component maps for $n \geq 1$ are
$$
(\Phi \ast \Psi)_{n} = \sum_{i+j=n+1}(-1)^{(i+1)(j+1)}\phi_{i}(\psi_{j} \otimes |\I|^{(j+1) \otimes (i-1)})
$$
\end{defn}

\begin{defn}[\cite{Bor1}]
Let $\Psi, \Phi$ be $A_{\infty}$-morphisms from $M$ to $M'$. $\Psi$ and $\Phi$ are homotopic if there is a set of maps
$\{h_{i}\}$ with $h_{i} : M \otimes A^{\otimes\,(i-1)} \longrightarrow M'[i]$ such that
\begin{equation}
\begin{split}
\psi_{i} - \phi_{i} = \sum_{i+j=n+1}(-1)^{(i+1)j}m'_{i}&(h_{j}\otimes |\I|^{j \otimes (i-1)})\\
& + \sum_{i+j=n+1} (-1)^{(i+1)j}h_{i}(m_{j} \otimes |\I|^{j \otimes (i-1)}) \\
& \hspace{0.5in} + \sum_{i+j=n+1,k>0}(-1)^{k(j+1) + j(i+1)}h_{i}(\I^{\otimes\,k} \otimes \mu_{j} \otimes |\I|^{j \otimes (i-k-1)})
\end{split}
\end{equation}
and for $i > 1$, $h_{i}(\xi \otimes a_{1} \otimes \cdots \otimes a_{i-1}) = 0$ when $a_{j} = 1_{A}$ for some $j$.
\end{defn} 

\noindent The sign convention used in the previous definitions is the one in Keller \cite{BKel} with the Koszul sign rule
$$
(f \otimes g)(x \otimes y) = (-1)^{|f||y|}(f(x) \otimes g(y))
$$
Thus, as can be checked directly, the composition of morphisms is a morphism for this sign convention, and homotopy of morphisms is an equivalence relation (or see the appendix). With these definitions, we are equipped to consider right $A_{\infty}$-modules up to homotopy equivalence. The following is our version of a standard result in the study of $A_{\infty}$-modules: 

\begin{prop}\label{prop:typeAcancel}
Let $(M, \{m_{i}\})$ be a strictly unital, right $A_{\infty}$-module over $(A, \{\mu_{i}\})$, and let $(\overline{M},\overline{m}_{1})$ be a chain complex. Suppose there exist chain maps $\iota : (\overline{M}, \overline{m}_{1}) \longrightarrow (M,m_{1})$ and $\pi : (M,m_{1}) \longrightarrow (\overline{M}, \overline{m}_{1})$, and a map $H : M \longrightarrow M[1]$ satisfying
\begin{align}
\pi \circ \iota &= \I_{\overline{M}}\\
\iota \circ \pi - \I_{M} &= m_{1}\circ H + H \circ m_{1}\\
H \circ \iota &= 0\\
\pi \circ H &= 0 \\
H^{2} &= 0
\end{align}
Then there are maps $\overline{m}_{i}: \overline{M} \otimes A^{\otimes\,(i-1)} \rightarrow \overline{M}$ for $i \geq 2$ such that
$\{\overline{m}_{i}\}_{i=1}^{\infty}$ defines a strictly unital right $A_{\infty}$-module structure on $\overline{M}$. This structure is homotopy equivalent to $(M, \{m_{i}\})$ through strictly unital morphisms which extend $\pi$ and $\iota$. 
\end{prop}

\noindent The proof supplies an explicit formula for computing $\overline{m}_{i}$ and the morphisms in the homotopy equivalence. First, we introduce some notation to simplify the formulas.

\begin{defn}
For positive integers $i_{1}, \ldots, i_{k}$ let $$N(i_{1}, \ldots, i_{k}) = \sum_{j} (i_{j} - 1)$$ and 
$$\alpha(i_{1}, \ldots, i_{k}) = \sum_{1 \leq r < s \leq k} (i_{r}-1)(i_{s}-1)$$
\end{defn}

\begin{defn}
Let $i_{j} \geq 2$ be integers for $j=1,\ldots, k$. By $[i_{1}, \ldots, i_{k}]$  we will mean the composition 
$$
(m_{i_{1}})(H \otimes |\I|^{\otimes\,(i_{1}-1)})(m_{i_{2}} \otimes |\I|^{i_{2} \otimes (i_{1} - 1)}) \cdots (H \otimes |\I|^{\otimes(I-i_{k})})(m_{i_{k}} \otimes |\I|^{i_{k} \otimes (I-i_{k})})
$$
where we alternate between applying $m_{i_{j}}$ to the first $i_{j}$ entries in the tensor product, and applying $H$ to the first factor in the tensor product. 
\end{defn}

\noindent Using this notation, we can define the action, morphisms, and homotopy. First, for $n \geq 2$ define a map $M \otimes A^{\otimes (n-1)} \longrightarrow M[n-2]$ by
$$
\Sigma_{n} = \sum_{\footnotesize\begin{array}{c} N(i_{1},\!i_{2},\!\ldots,\!i_{k})=\!n-1 \\ i_{j} \geq 2 \end{array}\normalsize} (-1)^{\alpha(i_{1},\ldots,i_{k})} [i_{1}, \ldots, i_{k}]
$$
\\
\ \\
\noindent We use $\Sigma_{n}$ to define $\overline{m}_{n}$ for $n \geq 1$:
$$
\overline{m}_{n} :=  \pi \circ \Sigma_{n} \circ \big(\iota \otimes \I^{\otimes\,(n-1)}\big)
$$
For $n=1$ we use the boundary map $\overline{m}_{1}$. Then $\{\overline{m}_{i}\}_{i=1}^{\infty}$ equips $\overline{M}$ with the structure of a right $A_{\infty}$-module. \\
\ \\
\noindent The morphisms which induce the homotopy equivalence are similarly defined. For $n=1$ we will use $\pi_{1} = \pi$ and $\omega_{1} = \iota$, while for $n > 1$ we use
$$
\pi_{n} := (-1)^{n}\left(\pi \circ \Sigma_{n} \circ \big(H \otimes |\I|^{\otimes\,(n-1)}\big)\right)
$$
$$
\omega_{n} := H \circ \Sigma_{n} \circ \big(\iota \otimes \I^{\otimes\,(n-1)}\big)
$$
The additional $H$ means that these are maps $\pi_{n} : M \otimes \I^{\otimes (n-1)} \longrightarrow \overline{M}[n-1]$ and $\omega_{n}: \overline{M} \otimes \I^{\otimes (n-1)} \longrightarrow M[n-1]$. As defined, these morphisms satisfy the relations $\Pi \circ \Omega = I_{\overline{M}}$ and $\Omega \circ \Pi \simeq_{\Lambda} I_{M}$ where $\lambda_{1} = H$ and
$$
\lambda_{n} := (-1)^{n}\left( H \circ \Sigma_{n} \circ \big(H \otimes |\I|^{\otimes\,(n-1)}\big)\right)
$$
and all homotopy equivalences occur in the category of (right) $A_{\infty}$-modules. 
\ \\
\noindent We note that even when $m_{i} \equiv 0$ for $i > 2$, a homotopy equivalence as described in \ref{prop:typeAcancel} can have higher order action terms. Indeed, the new module structure is given by
$$
\overline{m}_{n} = (-1)^{\epsilon} \pi [2,2,\ldots,2] (\iota \otimes \I^{n})
$$
where there are exactly $n-1$ $2$'s inside the square brackets and $\epsilon = 0$ if $n \equiv 1,2$ modulo 4, and $\epsilon = 1$ if $n \equiv 3,4$ modulo $4$. Thus, in all cases
$$
\overline{m}_{2} = \pi \circ m_{2} \circ (\iota \otimes \I)
$$
just comes from appropriately adjusting $m_{2}$. The effect of $\pi$, however, is substantial when doing calculations. With this observation and \ref{prop:typeAcancel}, we can, by directly analyzing the diagrams before and after a Reidemeister move, see that the $A_{\infty}$-module structure is preserved up to homotopy equivalence. This affords us the difficult part of

\begin{thm}
Let $\lefty{\mathcal{T}}$ be an inside tangle with boundary $P_{2n}$. 
\begin{enumerate}
\item Let $\lefty{T}$ be a diagram for $\lefty{\mathcal{T}}$ in $\leftHalf$. If $\mathfrak{o}_{1}$ and $\mathfrak{o}_{2}$ are two orderings of $\cross{\lefty{T}}$ then $\leftComplex{T, \mathfrak{o}_{1}}$ and $\leftComplex{T,\mathfrak{o}_{2}}$ are isomorphic type $A$ structures. 
\item If $\lefty{T_{1}}$ and $\lefty{T_{2}}$ are two diagrams for $\lefty{\mathcal{T}}$, then $\leftComplex{\lefty{T_{1}}}$ and $\leftComplex{\lefty{T}_{2}}$ are homotopy equivalent type $A$ structure. 
\end{enumerate}
\end{thm}

\begin{cor}
The homotopy type of the type A structure $\leftComplex{\lefty{T}}$, for any diagram $T$ of an inside tangle $\lefty{\mathcal{T}}$, is a tangle invariant.
\end{cor}

\noindent We will not prove these theorems here, as the proofs are modifications of those for the type D structure for an outside tangle found in \cite{typeD}. In addition, there are easier ways to prove these results once we have generalized the gluing theory in section \ref{sec:pairing}. Instead we content ourselves with computing some examples using \ref{prop:typeAcancel} which will illustrate the argument.\\
\ \\
\noindent{\bf How we will use this:} Suppose we have a chain complex $\{C_{i}\,|\,i\in\Z\}$ with explicit generators for each {\em free} chain group $C_{i}$. 
If the generators of $C_{i}$ are $\{x_{1},\ldots, x_{n}\}$ and those for $C_{i-1}$ are $\{y_{1},\ldots, y_{m}\}$ we can find a homotopy $H$ as in proposition \ref{prop:typeAcancel} by searching through the images $\partial\,x_{i} = \sum a_{i}^{j}y_{j}$ to find one where $a^{i}_{j} = u$ is a unit, for some $j$. We will reorder the generators so that this occurs for $i=j=1$. We can then construct a new chain complex on where the other chain groups and boundary maps are taken to be the same, but $C_{i}'$ is spanned by $x_{2}', \ldots, x_{n}'$ and $C_{i-1}'$ is spanned by $y_{2}', \cdots, y_{m}'$. We let $p(x_{i}) = x_{i}'$ for $i > 1$ and $p(x_{1})=0$, and likewise for the $y_{j}$. Otherwise $p$ is the identity.  The new boundary $\partial_{i}': C_{i'} \rightarrow C_{i-1}'$ is given by $\partial'\,x_{i}' = \big(p\circ\partial\big)\big(x_{i} - a^{1}_{i}u^{-1}\,x_{1}\big)$. If we let $\iota(x_{i}') = x_{i} - a^{1}_{i}u^{-1}\,x_{1}$
and $H(y_{1})=-u^{-1}x_{1}$, $H(\beta) = 0$ otherwise, we are in the situation envisioned in proposition \ref{prop:typeAcancel}. The map $\pi$ is the quotient map found by quotienting out the subcomplex generated by $\{x_{1},\partial\,x_{1}\}$. The formulas involving $p$ are a specific presentation of this quotient complex for the specific basis. $\partial'$ is computed by calculating $\pi \circ \partial$ in this presentation. \\
\ \\
\noindent We now compute $\overline{m}_{2}(x_{i}' \otimes e)$. Since $\overline{m}_{2} = \pi \circ m_{2} \circ (\iota \otimes \I)$ we first compute
$$m_{2}((x_{i} - a^{1}_{i}u^{-1}\,x_{1}) \otimes e) = m_{2}(x_{i}\otimes e) - a^{1}_{i}u^{-1}m_{2}(x_{1}\otimes e)$$
and then compute $\pi$. In particular, suppose $\langle \partial x_{j}, y_{1} \rangle = a^{1}_{j} = 0$, but $m_{2}(x_{j} \otimes e) = a\,y_{1} + Y$, then 
$\overline{m}_{2}(x_{j}' \otimes e) =  \pi(a\,y_{1} + Y)$ $= a(u^{-1} \sum_{j >1} a_{1}^{j}y_{j}) + Y$. This is the same process as for adjusting the
boundary maps above. \\
\ \\
\noindent However, now suppose $\langle \partial x_{j}, y_{1} \rangle =  a^{1}_{j} = 0$, but $m_{2}(x_{j} \otimes e_{1}) = a\,y_{1} + Y$ and $m_{2}(x_{1} \otimes e_{2}) = \Omega$. Then $\overline{m}_{3}(x_{j}' \otimes e_{1} \otimes e_{2}) = \pi(m_{3}(x_{j} \otimes e_{1} \otimes e_{2}) - \pi(m_{2}(H \circ m_{2}(x_{j} \otimes e_{1})\otimes e_{2})$. We concentrate upon $\pi(m_{2}(H \circ m_{2}(x_{j} \otimes e_{1})\otimes e_{2}) = \pi(m_{2}(H(a\,y_{1} + Y)))$ $=\pi(m_{2}(-u^{-1}a\,x_{1} \otimes e_{2}) = -u^{-1}a\Omega$. We thus pick up a higher order action.\\
\ \\

\section{Examples of the type A structure}\label{sec:examplesA}
 
\noindent{\bf Example I (Reidemeister tangles):} The three tangles below appear in the local description of the Reidemeister moves. We will analyze each in turn.\\
\ \\
$$
\inlinediag[0.5]{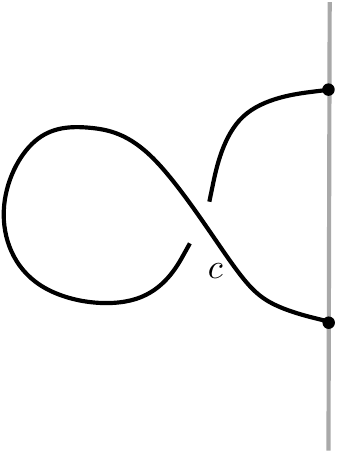} \hspace{1in} \inlinediag[0.5]{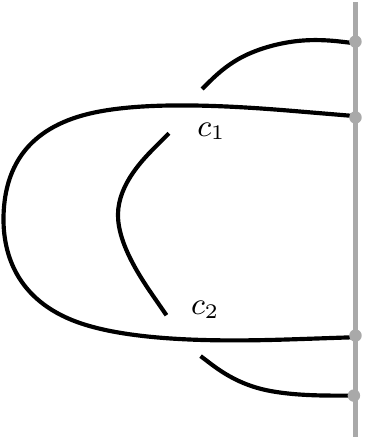} \hspace{1in} \inlinediag[0.5]{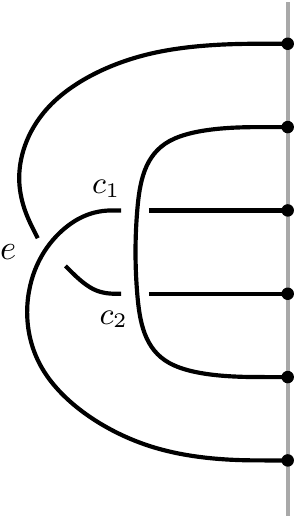}
$$
\noindent{\em (a) First move:} For the RI move we have a tangle diagram  $R_{I}$, over $P_{2}$ with one crossing. There are two resolutions, corresponding to $\rho = 0$ and $\rho = 1$, and the unique matching on $P_{2}$. Since the crossing is right handed, the $0$ resolution has a single free circle. Writing the decoration on the cleaved circle first, we can think of the decorated resolutions as $z_{++}$, $z_{+-}$, $z_{-+}$, and $z_{--}$ which occur in bigrading $(0, 5/2)$, $(0,1/2)$, $(0,1/2)$ and $(0,-1/2)$. For the $1$ resolution there is only the cleaved circle, so we get two state $t_{+}$ in grading $(1,5/2)$ and $t_{-}$ in $(1,3/2)$. We can compute $d_{APS}$ as $d_{APS}(z_{++}) = t_{+}$ and $d_{APS}(z_{-+}) = t_{-}$. $\mathcal{B}\Gamma_{1}$ has two non-idempotent elements $\lefty{e_{C}}$ and $\righty{e_{C}}$. The actions of these elements are
$z_{+\ast} \righty{e_{C}} = z_{-\ast}$ and $t_{+} \righty{e_{C}} = t_{-}$, since $\righty{e}_{C}$ only changes the sign on the cleaved circle. On the other hand, $z_{+-} \lefty{e_{C}} = t_{-}$ is the only non-trivial action for $\lefty{e}_{C}$. Since $t_{+}$ is not in the image of the action, or $d_{APS}$ except for $z_{++}$ we can cancel both to result in the homotopy equivalent structure with generators $z_{-+}, z_{+-}, z_{--}$ and $t_{-}$. We now cancel $t_{-}$ through the image of $z_{-+}$. To compute the new action on $z_{+-}$ and $z_{--}$ we consider $\iota$, which in this case is just inclusion, followed by $m_{2}(z_{\ast-} \otimes e)$ followed by projection. For $z_{+-}$ there is the non-trivial action $m_{2}(z_{+-} \otimes \righty{e_{C}}) = z_{--}$, which projects to an action as well. However, while $m_{2}(z_{+-} \otimes \lefty{e_{C}}) = t_{-}$, the projection will kill this image.  Furthermore, all the higher actions vanish since $m_{2}$ acts trivially on $z_{--}$, and any computation of $\overline{m}_{n}$ for $n > 2$ starts with $H \circ m_{2}(z_{+-} \otimes \lefty{e_{C}}) = z_{-+}$, but the action on $z_{-+}$ is trivial for all non-idempotents. The idempotent will fix $z_{-+}$, but this will be killed under $\pi$, or $H$, and the computation cannot proceed. Thus, $\leftComplex{R_{I}}$ is isomorphic to $\alpha_{+} = z_{+-}$  in grading $(0,1/2)$ and $\alpha_{-}=z_{--}$ in $(0,-1/2)$ with
$d_{APS}\equiv 0$ and the only non-trivial action term being $\alpha_{+} \cdot \righty{e_{C}} = \alpha_{-}$. This is isomorphic to $\leftComplex{U_{2}}$ where $U_{2}$ is the planar matching on $P_{2}$ found from untwisting the crossing.\\
\ \\
\noindent{\em (b) Second move:} For the RII move we analyze the tangle below, $R_{II}$ over $P_{4}$, with two opposite crossings. Thus $n_{+}=1$ and $n_{-}=1$ for every choice of orientation. We label the crossings from top to bottom. Now consider the states corresponding to the $01$ resolution. There is a free circle in this resolution, and we can divide the states into $S^{+}_{01}$ and $S^{-}_{01}$ based on the decoration of the circle (we do this regardless of the matching $\righty{m}$ used to construct the state). $d_{APS}$ maps $S^{+}_{01}$ isomorphically to $S_{11}$ and $S_{00}$ isomorphically to $S^{-}_{01}$. The action $m_{2}(\xi \otimes e)$ for $\xi$ in $S^{+}_{01}$ has image in $S^{+}_{01}$ since it will not change the decoration on the $+$ free circle, and merging the $+$ free circle does not change the boundary of the state. Consequently if we cancel along the isomorphism from $S^{+}_{01}$ to $S_{11}$ the image of $H$ is in $S^{+}_{01}$ and the image of $\iota$ on $\nu \in S_{10}$ is a sum $\nu + \nu'$ where $\nu' \in S^{+}_{01}$. Thus $\pi \circ m_{2} \circ (\iota \otimes \I)$ will have image equal to the part of $m_{2}(\nu \otimes e)$ in $S_{10}$, since $\pi$ will kill $S^{+}_{01}$. The only element $e$ for which the image may not be in $S_{10}$ is $\lefty{e_{\gamma}}$ for the unique class of bridges $\gamma$ in the boundary of any element in $S_{10}$. Its action would have image in $S_{11}$, but does not contribute to $\underline{m}_{n}$ for $n>1$ since $H:S_{11} \rightarrow S_{01}^{+}$, and thus any additional actions stay in $S_{01}^{+}$, which will be killed by $\pi$. \\
\ \\
\noindent The effect of the cancellation, therefore, is to reduce our module to $S_{10} \oplus S_{00} \oplus S_{01}^{-}$ with action
defined by restricting the image of $m_{2}$ to the remaining summands. Now a similar argument applies to the isomorphism found by the image $d_{APS}|S_{00}$ in $S_{01}^{-}$. Now, however, no element from $S_{10}$ can have a term in its action or boundary within $S_{01}$, so these will be unchanged. After the cancellation we obtain all the states in $S_{10}$ having trivialized the $\lefty{\gamma}$ action, but otherwise left the action unchanged. This is the same type $A$ structure as for the matching of the top point in $P_{4}$ with the bottom, and the second with the third. Thus it is isomorphic to the structure obtained after removing the crossings with the RII move. Being in $S_{10}$ means the states have no free circle, and just receive grading based on the cleaved circles. Furthermore, they are shifted by $(1,1) + (-1,1-2\cdot1) = (0,0)$ when we account for the resolution and the crossings. Thus, as a bigraded type $A$ structure the RII tangle is homotopy equivalent to the planar matching obtained from the RII move.\\
\ \\
\noindent{\em (c) Third move (sketch):} Let $R_{b}$ be the diagram before the move and $R_{a}$ be the diagram after:
$$
\inlinediag[0.5]{RIIIL} \hspace{1in} \inlinediag[0.5]{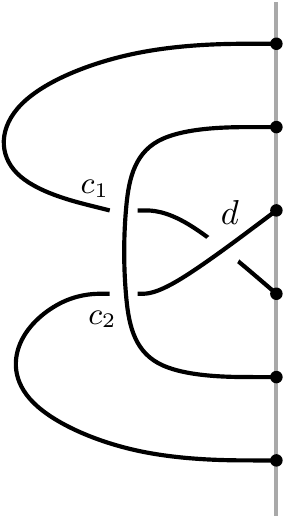}
$$
\noindent In each if we $0$ resolve the second crossing from the top we obtain a diagram with an RII move. As usual the states with a $1$ resolution here give rise to the same type A structure. It is enough then to see what happens in the $0$ resolved sub-module. As with the RII move we can use $d_{APS}$ to leave $S_{100}$ with its action intact, including the action of $\lefty{e_{\gamma_{1}}}$ which has image in $S_{110}$. However, $d_{APS}$ now maps $S^{+}_{001}$ to $S_{101}\oplus S_{011}$, isomorphically to each factor, given by minus the Khovanov maps. We let $\xi$ be the state in $S_{101}$ then $\xi'$ is the corresponding state, found by planar isotopy, in $S_{011}$. The effect of $\pi$ is to identify $\xi$ with $-\xi'$. Now, let $\nu$ be a state in $S_{100}$ and let $\nu' = H \circ d_{APS}(\nu)$ in $S^{+}_{001}$. Then $\iota(\nu) = \nu + \nu'$. The action $m_{2}((\nu  + \nu')\otimes \lefty{e_{\gamma_{2}}}) = m_{2}(\nu\otimes \lefty{e_{\gamma_{2}}})$ since the action of $\lefty{e_{\gamma_{2}}}$ on $S^{+}_{001}$ is trivial ($\gamma_{2}$ is used in the calculation of $d_{APS}$ for these states). If $m_{2}(\nu\otimes \lefty{e_{\gamma_{2}}})$ is non-zero in $S_{101}$, then the effect of $\pi$ is to identify it with $-m_{2}(\nu\otimes \lefty{e_{\gamma_{2}}})'$. \\
\ \\
\noindent If we repeat this argument with $R_{a}$, with the same crossing ordering, we get $S_{001}$ being the planar matching diagram and $S^{+}_{100}$ being used in the cancellation process. For $\nu$ in $S_{001}$ the effect of $\lefty{e_{\gamma_{2}}}$ is the same as before, but as it takes image in $S_{011}$ it occurs with a minus sign. On the other hand, the image of $\lefty{e_{\gamma_{1}}}$ will be in $S_{101}$, occurring with a minus sign, due to the crossing ordering, and thus will be identified with $-(-\eta)$ where $\eta$ is the image $m_{2}(\nu\otimes \lefty{e_{\gamma_{1}}})$ from $R_{b}$ in the previous paragraph. As such the actions of the bridges will be the same, and the APS-complexes will be the same. It is straightforward to see that the higher actions all vanish.\\
\ \\      
\noindent{\bf Example II (Hopf Tangle):} For the Hopf tangle over $P_{2}$
$$
\inlinediag[0.5]{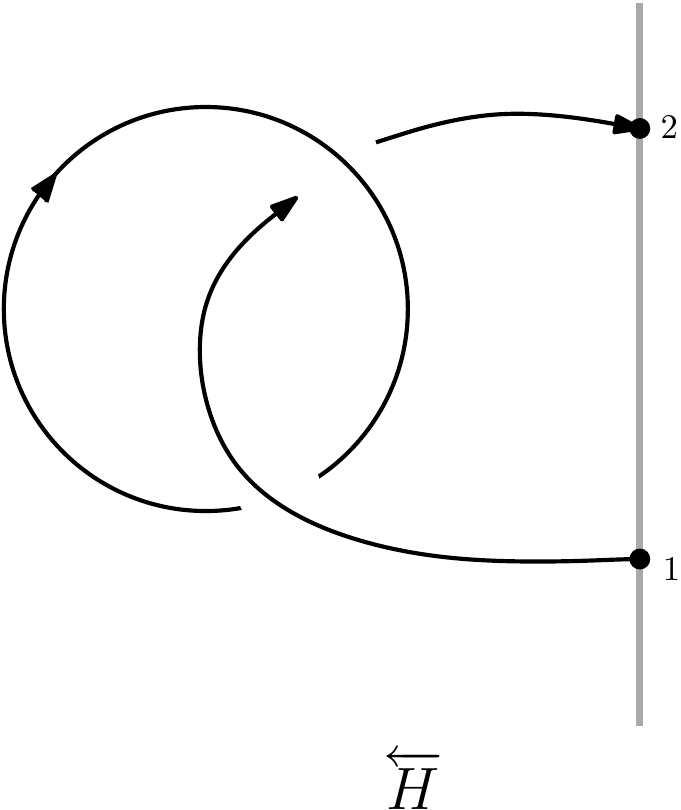}
$$
\noindent we will enumerate the crossings as shown, and write states with the decoration of the cleaved circle first. For the moment we will ignore orientations. There are four states $s^{00}_{\pm,\pm}$ in homological grading $0$ and quantum gradings $\pm 1/2 \pm 1$. There are two states $s^{10}_{\pm}$ for the $10$ resolution, and two states $s^{01}_{\pm}$ for the $01$. These occur in the bigradings $(1,1\pm 1/2)$. Finally, there are four states $s^{11}_{\pm\pm}$ with bigrading $(2,2\pm 1/2 \pm 1)$. For these states:
\begin{enumerate}
\item $d_{APS}$ is computed as 
\begin{align*}
s^{00}_{++} &\rightarrow s^{10}_{+} + s^{01}_{+}\\
s^{00}_{-+} &\rightarrow s^{10}_{-} + s^{01}_{-}\\
s^{10}_{+}  &\rightarrow s^{11}_{+-}\\
s^{10}_{-}  &\rightarrow s^{10}_{--}\\
s^{01}_{+}  &\rightarrow -s^{11}_{+-}\\
s^{01}_{-}  &\rightarrow -s^{11}_{--}\\
\end{align*} 
\item The action of $\righty{e_{C}}$ simply takes $s^{\ast}_{+\ast} \rightarrow s^{\ast}_{-\ast}$ where $\ast$ matches anything in those spots.
\item The action of $\lefty{e_{C}}$ is given by
\begin{align*}
s^{00}_{+-} &\rightarrow s^{10}_{-} + s^{01}_{-}\\
s^{10}_{+}  &\rightarrow s^{11}_{-+}\\
s^{01}_{+}  &\rightarrow -s^{11}_{-+}\\
\end{align*} 
\end{enumerate} 
If we cancel $s^{00}_{++}$ with $s^{10}_{+}$ we will have no effect except to remove these generators, as $s^{10}_{+}$ does not occur in the image of $d_{APS}$ or the action for any other state. Once we have done that, we can cancel $s^{01}_{+}$ with $-s^{11}_{+-}$ with no other effect, since $s^{11}_{+-}$ only occurs in the image of a previously canceled state. $s^{11}_{--}$ will then appear only in the image $d_{APS}(s^{01}_{-})$ and $d_{APS}(s^{10}_{-})$ (since $s^{11}_{+-}$ has been canceled, otherwise we would also need to include it in the image of $\righty{e_{C}}$). As, $d_{APS}(s^{01}_{-}) = -s^{11}_{--}$ we can cancel these without affecting the rest of the maps. Finally we can cancel $s^{00}_{-+}$ with $s^{10}_{-}$. $s^{10}_{-}$ occurs as $s^{00}_{+-}\cdot \lefty{e_{C}}$, but there are no other terms to consider, so the effect of the cancellation (through the projection $\pi$) is to cancel this portion of the action of $\lefty{e}_{C}$. \\
\ \\
\noindent Following these steps results in $s^{00}_{+-}$ and $s^{00}_{--}$ in bigrading $(0,-1/2)$ and $(0,-3/2)$, and $s^{11}_{++}$ in bigrading $(2,7/2)$ and $s^{-+}$ in bigrading $(2,5/2)$. The residual action is that of $\righty{e_{C}}$, which takes $s^{00}_{+-}$ to $s^{00}_{--}$ and $s^{11}_{++}$ to $s^{11}_{-+}$. \\
\ \\
\noindent The Hopf tangle will either have two positive or two negative crossings, depending upon the orientation of the components. If there are two positive crossings we will shift the bigrading up $(0,2)$. Otherwise, for negative crossings, we add $(-2, -4)$ to each bigrading. \\
\ \\
\noindent Consequently, for the {\em positive Hopf tangles} we will have $\F_{(0,3/2)} \stackrel{\righty{e_{C}}}{\longrightarrow} \F_{(0,1/2)}$ and $\F_{(2,11/2)} \stackrel{\righty{e_{C}}}{\longrightarrow} \F_{(2,9/2)}$. \\
\ \\
\noindent For {\em negative} Hopf tangles we will have $\F_{(-2,-9/2)} \stackrel{\righty{e_{C}}}{\longrightarrow} \F_{(-2,-11/2)}$ and $\F_{(0,-1/2)} \stackrel{\righty{e_{C}}}{\longrightarrow} \F_{(0,-3/2)}$.
%\\
%\ \\
%\noindent{\bf Example III: A trefoil} We consider the following right handed trefoil tangle. Here there are three crossings, and $n_{+} = 3$. We will shift by $(0,3)$ after we compute the type A structure. 

\section{Gluing inside and outside tangles}\label{sec:pairing}

\noindent Let $\lefty{\mathcal{T}_{1}}$ be an inside tangle for $P_{2n}$ and $\righty{\mathcal{T}_{2}}$ be an outside tangle. We let $\mathcal{T} = \lefty{\mathcal{T}_{1}}\#\righty{\mathcal{T}_{2}}$ be the link in $\R^{3}$ obtained by gluing $\lefty{\R^{3}}$ to $\righty{\R^{3}}$ and thereby gluing $\lefty{\mathcal{T}_{1}}$ to  $\righty{\mathcal{T}_{2}}$ along $P_{2n}$. Likewise, if $\lefty{T_{1}}$ is a diagram for $\lefty{\mathcal{T}_{1}}$ in $\leftHalf$ and $\righty{T_{2}}$ is a diagram for $\righty{\mathcal{T}_{2}}$ in $\rightHalf$, we can glue these diagrams along $P_{2n}$ to obtain a diagram $T$ for $\mathcal{T}$.\\
\ \\
\noindent In \cite{typeD} we showed how to associate a bigraded type D structure to $\righty{T_{2}}$ whose homotopy type is an isotopy invariant of $\righty{\mathcal{T}_{2}}$. In particular, we constructed a bigrading preserving map
$$
\righty{\delta_{T}} : \rightComplex{\righty{T}} \longrightarrow \mathcal{B}\bridgeGraph{n} \otimes_{\mathcal{I}} \rightComplex{\righty{T}} [(-1,0)]
$$
which satisfies the type D structure equation
$$
(\mu_{\mathcal{B}\bridgeGraph{n}} \otimes \I)\,(\I \otimes \righty{\delta_{T}}) \, \righty{\delta_{T}}  +  (d_{\Gamma_{n}} \otimes |\I|)\,\righty{\delta_{T}} = 0
$$
where $\mu_{\mathcal{B}\bridgeGraph{n}}: \mathcal{B}\bridgeGraph{n} \otimes \mathcal{B}\bridgeGraph{n} \rightarrow \mathcal{B}\bridgeGraph{n}$ is the multiplication map on $\mathcal{B}\bridgeGraph{n}$.
\ \\ 
\begin{defn}
By $\leftComplex{T_{1}} \boxtimes \rightComplex{T_{2}}$ we mean the bigraded module
$$
\leftComplex{T_{1}} \otimes_{\mathcal{I}_{n}} \rightComplex{T_{2}}
$$
equipped with the map
$$
\partial^{\boxtimes}(x \otimes y) =  d_{APS}(x) \otimes |y| + (m_{2, T_{1}} \otimes \I)(x \otimes \righty{\delta_{T_{2}}}(y))
$$ 
\end{defn} 

\begin{prop}
$\partial^{\boxtimes}$ is a $(1,0)$ differential map on $\leftComplex{T_{1}} \boxtimes \rightComplex{T_{2}}$. 
\end{prop}

\noindent{\bf Proof:} First, we rewrite $\partial^{\boxtimes}$ as an operator: 
$$
\partial^{\boxtimes} = d_{APS} \otimes |\I| + (m_{2,T_{1}} \otimes \I)(\I \otimes \righty{\delta_{T_{2}}})
$$
We note that $\I \otimes \righty{\delta_{T_{2}}}$ is a $(1,0)$ map $\leftComplex{T_{1}} \otimes_{\mathcal{I}} \rightComplex{T_{2}} \rightarrow  \leftComplex{T_{1}} \otimes_{\mathcal{I}} \mathcal{B}\Gamma_{n} \otimes_{\mathcal{I}} \rightComplex{T_{2}}$ , while $m_{2,T_{1}} \otimes \I$ preserves the bigrading as a map $\leftComplex{T_{1}} \otimes_{\mathcal{I}} \mathcal{B}\Gamma_{n} \otimes_{\mathcal{I}} \rightComplex{T_{2}} \rightarrow \leftComplex{T_{1}} \otimes_{\mathcal{I}} \rightComplex{T_{2}}$. In addition, $d_{APS}$ is a $(1,0)$ map. Hence $\partial^{\boxtimes}$ is a $(1,0)$ map. We now verify that $\partial^{\boxtimes}$ is a differential. The rest of the result follows from section \ref{sec:boxtimes} in the appendix, and that when $m_{i, T_{1}} = 0$ for $i > 2$, $m_{1,T_{1}} = d_{APS}$, implies that $\partial^{\boxtimes}$ above coincides with the definition in the appendix. $\Diamond$\\
\ \\
\noindent By $\complex{T}$ we will mean the usual bigraded Khovanov complex over $\Z$, equipped with its invariant bigrading. 

\begin{prop}
$\complex{T}  \cong (\leftComplex{T_{1}} \boxtimes \rightComplex{T_{2}}, \partial^{\boxtimes})$.
\end{prop}

\noindent{\bf Important Comment:} We have not required that the orientations on $T_{1}$ and $T_{2}$ match along $P_{n}$. If they do, $\complex{T}$ is exactly the Khovanov complex from \cite{Khov}, as described in \cite{Bar1}. However, the statement still holds even if the orientations do not match. The Khovanov complex in the latter case is for a link with a finite number of orientation changes, constructed in the same manner as before. Now, however, it has an invariant bigrading only as long as isotopies do not take a strand across a point where the orientation changes. In the latter case there is a bigrading shift of $\pm (1,3)$ due to the conversion of a negative crossing to a positive crossing, or vice-versa.\\ 
\ \\
\noindent{\bf Proof:} We start by identifying the generators of $\leftComplex{T_{1}} \otimes_{\mathcal{I}} \rightComplex{T_{2}}$ with the generators of $\complex{T}$. For $(r_{1},s_{1}) \otimes_{\mathcal{I}} (r_{2},s_{2})\neq 0$ we need that $I_{\partial(r_{1},s_{1})} \cdot (r_{2},s_{2}) \neq 0$ since
$(r_{1},s_{1}) \cdots I_{\partial(r_{1},s_{1})} =(r_{1},s_{1})$. However, only $I_{\partial(r_{2},s_{2})} \cdot (r_{2},s_{2}) \neq 0$, so $\partial(r_{1},s_{1}) = \partial(s_{2},r_{2}) = (L,\sigma)$. If $r_{1} = (\rho_{1}, \righty{m_{1}})$ and $r_{2} = (\lefty{m_{2}},\rho_{2})$  we use $\righty{m}_{1} = \righty{L}$ to identify $\righty{m_{1}}$ with the {\em arcs} in $\rho_{2}(T_{2})$, and likewise we can identify $\lefty{m_{2}} = \lefty{L}$ with the arcs in $\rho_{1}(T_{1})$. Furthermore, $s_{1}$ and $s_{2}$ to $\sigma$, so we can take $\rho_{1}(T_{1})\# \rho_{2}(T_{2})$ with $s_{1} \# s_{2}$ to get a resolution diagram for $T$ where every circle is unambiguously decorated with $\pm$. \\
\ \\
\noindent Furthermore, we can reverse the construction. If $\rho$ is a resolution of $T$, we let $\rho_{1}$ be $\rho$ restricted to those crossings in $\leftHalf \cap T = T_{1}$ and $\rho_{2}$ be $\rho$ restricted to $\rightHalf \cap T = T_{2}$. Furthermore, the arcs in $\rho_{2}(T_{2})$ form an (outside) planar matching $\righty{m_1}$, and we define $r_{1} = (\rho_{1}, \righty{m_{1}}$. Likewise the arcs of $\rho_{1}(T_{1})$ define an (inside) planar matching $\lefty{m_{2}}$ and we let $r_{2} = (\lefty{m_2}, \rho_{2})$. A generator of $\complex{T}$ is a pair $(\rho,s)$ where $s$ is a decoration of $\circles{\rho(T)}$. By restriction $s$ defines decorations, $s_{1}$, $s_{2}$ on $r_{1}(T_{1})$ and $r_{2}(T_{2})$ with $\partial(r_{1},s_{1}) = \partial(r_{2},s_{2})$. It is straightforward to see that $(r_{1},s_{1}) \otimes_{\mathcal{I}} (r_{2},s_{2}) = (\rho,s)$, so that this is the inverse of the previous map. \\ 
\ \\
\noindent Furthermore, the bigrading of $(r_{1},s_{1}) \otimes_{\mathcal{I}} (r_{2},s_{2})$ coming from the tensor product is identical to that of $(\rho,s)$ from the construction of $\complex{T}$. The bigrading of $(r_{1},s_{1}) \otimes_{\mathcal{I}} (r_{2},s_{2})$ is the sum $(h(r_{1})-n_{-}(T_{1}), h(r_{1}) + q(r_{1},s_{1}) + 1/2 \iota(\partial(r_{1},s_{1})) + n_{+}(T_{1}) - 2n_{-}(T_{1}))$ $+ (h(r_{2})-n_{-}(T_{2}), h(r_{2}) + q(r_{2},s_{2}) + 1/2 \iota(\partial(r_{2},s_{2})) + n_{+}(T_{2}) - 2n_{-}(T_{2}))$. However $h(r_{1}) + h(r_{2})$ is the number of $1$ resolutions in $\rho_{1}(T_{1})$ added to the number in $\rho_{2}(T_{2})$, which equals the total number in $\rho(T)$. Likewise, since they are counts over crossings, $n_{+}(T_{1}) + n_{+}(T_{2}) = n_{+}(T)$ and $n_{-}(T_{1}) + n_{-}(T_{2}) = n_{-}(T)$. Finally, $\iota(\partial(r_{1},s_{1})) + \iota(\partial(r_{2},s_{2})) = 2\iota(L,\sigma)$ so the second entry in the bigrading equals the sum of the decorations on the free circles in $r_{1}(T_{1})$ plus the sum of the decorations on the circles in $(L,\sigma)$ plus the sum of the decorations on the free circles in $r_{2}(T_{2})$. In $\rho(T)$ this is just the quantum grading for the usual Khovanov generator. Thus the bigrading of 
$(r_{1},s_{1}) \otimes_{\mathcal{I}} (r_{2},s_{2})$ is $(h(\rho) - n_{-}(T), h(\rho) + q(\rho,s) + n_{+}(T) - 2n_{-}(T))$ which is the bidgrading of $(r,s)$ in $\complex{T}$. The tensor product identifies $\complex{T}$ with $\leftComplex{T_{1}} \otimes_{\mathcal{I}} \rightComplex{T_{2}}$ as bigraded modules over $\Z$.\\
\ \\
\noindent To see that $\partial^{\boxtimes}$ is the $(1,0)$ Khovanov differential $\partial_{KH}$  under this isomorphism, we must first specify the order of the crossings to be used in calculating the signs in $\partial_{KH}$. The chain isomorphism type of $\complex{T}$ is unaffected by this choice of ordering, \cite{Khov}. If $\mathfrak{o}_{i}$ is the ordering of the crossings in $T_{i}$ and, then $\mathfrak{o}_{1} ||\mathfrak{o}_{2}$ is an ordering of the crossings for $T$, which we now fix. In short, all the crossings of the inside tangle $T_{1}$ come before all the crossings of $T_{2}$, and in the same order as in $T_{1}$. \\
\ \\
\noindent We compute $\partial^{\boxtimes}$ in stages. First $(d_{APS} \otimes |\I|)\big[(r_{1},s_{1}) \otimes_{\mathcal{I}} (r_{2},s_{2})\big]$ is a sum
over the crossings of $T_{1}$. For each crossing, $c$, we get either $0$ or $(-1)^{m}(-1)^{h(r_{2}})(r',s') \otimes_{\mathcal{I}} (r_{2},s_{2})$
where $m$ is the number of $1$ resultions in $(r_{1},s_{1})$ following $c$, $h(r_{2})$ is the total number of $1$ resolutions in $r_{2}$, and $(r',s')$ is as specified previously, which has $\partial(r',s') = \partial (r_{1},s_{1})$. Consequently, $m+h(r_{2})$ is the number of $1$ resolutions of $T$ following $c$ in our fixed order, and $(r',s') \otimes_{\mathcal{I}} (r_{2},s_{2})$ is a generator of $\complex{T}$. Following the definition of $d_{APS}$ this is precisely a term in $\partial_{KH}(r,s)$. In fact, the sum of these is precisely the terms in $\partial_{KH}(r,s)$ which have the same decorated cleaved link, and occur from a crossing change in $T \cap \leftHalf$. Those terms in $\partial_{KH}(r,s)$ which have the same decorated cleaved link, and occur from a crossing change in $\rightHalf \cap T$ correspond to terms in $(m_{2,T_{1}} \otimes \I)(\I \otimes \righty{\delta_{T_{2}}})$ applied to $(r_{1},s_{1}) \otimes_{\mathcal{I}} (r_{2},s_{2})$. In the definition of $\righty{\delta_{T_{2}}}(r_{2},s_{2})$ there is a term $I_{\partial(r_{2},s_{2})} \otimes d_{APS}(r_{2},s_{2})$. Since $\partial(r_{2},s_{2}) = \partial(r_{1},s_{1})$ we conclude that $(m_{2,T_{1}} \otimes \I)((r_{1},s_{1}) \otimes I_{\partial(r_{2},s_{2})} \otimes d_{APS}(r_{2},s_{2})) = (r_{1},s_{1} \otimes_{\mathcal{I}} d_{APS}(r_{2},s_{2})$. Note that the signs are also correct for $\partial_{KH}$ since the sign of a term in $I_{\partial(r_{2},s_{2})} \otimes d_{APS}(r_{2},s_{2})$ is $(-1)^{m}$ where $m$ is the number of $1$ resolutions in $r_{2}$ following the crossing that yields the term. As this crossing follows all those of $T_{1}$, the same sign is used in $\partial_{KH}$.\\
\ \\
\noindent This leaves the terms of $\partial_{KH}\big[(r_{1},s_{1}) \otimes_{\mathcal{I}} (r_{2},s_{2})\big]$ which change the decorated, cleaved link. We divide them into two groups, based on whether the crossing change giving the term occurs in $\leftHalf$ or $\rightHalf$. We start with those in $\rightHalf$. Each such crossing gives an active arc which is either in $\rightBridges{r_{2}}$ or $\dec{r_{2},s_{2}}$. In the first case, $\righty{\delta_{T_{2}}}(r_{2},s_{2})$ will have a term, or two terms,  $(-1)^{m}(\righty{e}_{\gamma} \otimes (r',s'))$ which corresponds to the crossing change. as before, $m$ is the number of $1$ resolved crossings following the crossing. This is the same sign as in $\partial_{KH}$, and the decorations on the decorated, cleaved link also follow the pattern for Khovanov homology. In $(m_{2,T_{1}} \otimes \I)(\I \otimes \righty{\delta_{T_{2}}})$ we get the term $(-1)^{m}(m_{2,T_{1}} \otimes \I)\big[(r_{1},s_{1}) \otimes \righty{e}_{\gamma} \otimes (r',s')\big]$. From the definition of $m_{2,T_{1}}$ the action of $\righty{e}_{\gamma}$ on $\leftComplex{T_{1}}$ is just to change the decorated, cleaved link to have the same boundary as $(r',s')$ (which occurs purely in $\righty{m_{1}}$). Consequently, we obtain the tensor product of compatible pairs, and we replicate the term in $\partial_{KH}$. The case of an arc in $\dec{r_{2},s_{2}}$ is similar, except only the decoration on one circle changes, and not the underlying cleaved link. This is the effect of $\righty{e_{C}}$, for that circle, on $\leftComplex{T_{1}}$. \\
\ \\
\noindent This leaves the terms of $\partial_{KH}$ which come from crossing changes in $\leftHalf$ that change the decorated, cleaved link. Let $c$ be such a crossing, and $\gamma$ be the active arc. Suppose $\gamma$ has image in $\leftBridges{\partial(r_{1},s_{1})}$, which we will denote by $\gamma'$. There is then a term in $\righty{\delta_{T_{2}}}(r_{2},s_{2})$ of the form $(-1)^{h(r_{2})}(\lefty{e_{\gamma'}} \otimes (r'_{2},s'_{2}))$ where $(r'_{2},s'_{2})$ is the result of $\gamma'$ surgery on $r(T_{2}) \cap \leftHalf$ which reflects the decoration changes necessary for the Khovanov differential. In $(-1)^{h(r_{2})}(m_{2,T_{1}} \otimes I)((r_{1},s_{1}) \otimes  \lefty{e_{\gamma'}} \otimes (r'_{2},s'_{2}))$ we get a sum over all the terms in $\partial_{KH}$ which correspond to $\gamma'$ and the decoration changes for $\lefty{e_{\gamma'}}$, but with sign $(-1)^{h(r_{2})}(-1)^{m}$ where $m$ is the number of $1$ resolved crossings following that for $\gamma$ (not $\gamma'$) in ordering on the crossings of $T_{1}$. However, $h(r_{2}) + m$ is the number of $1$ resolved crossings following that for $\gamma$ in the ordering on $T$. Thus, the sign is the same as that for $\partial_{KH}$. If $\gamma \in \dec{r_{1},s_{1}}$ then the argument is the same except that the term in $\partial_{KH}$ comes from the action of $(-1)^{h(r_{2})}(\lefty{e_{C}} \otimes (r_{2},s_{2,C})$ where $C$ is the cleaved circle whose decoration changes. Note that a crossing change can occur in multiple terms, but that with the decoration changes included, each crossing and decoration change occurs in precisely one way above. Thus we recover all the terms of $\partial_{KH}(r,s)$ with the correct signs from the ordering of crossings.  $\Diamond$\\
\ \\
\noindent The advantage of using $\leftComplex{T_{1}} \boxtimes \rightComplex{T_{2}}$ arises from the ability to separately simplify $\leftComplex{T_{1}}$ and $\rightComplex{T_{2}}$ without changing the homotopy type of $\leftComplex{T_{1}} \boxtimes \rightComplex{T_{2}}$. We show this in the appendix through a series of propositions which replicate, for our sign conventions, results from \cite{Bor1}. In particular, propositions \ref{prop:Asimplify} and the corollary to \ref{prop:Dsimplify} imply the following result.

\begin{prop}
Suppose $(N,\delta)$ is  homotopy equivalent, as a  type $D$ structure over $\mathcal{B}\Gamma_{n}$, to $\rightComplex{T_{2}}$, and $(M,\{m_{i}\})$ is homotopy equivalent to $\leftComplex{T_{1}}$, as a type $A$ structure. Then $(M,\{m_{i}\}) \boxtimes (N,\delta) \simeq \leftComplex{T_{1}} \boxtimes \rightComplex{T_{2}} \simeq \complex{T_{1}\#T_{2}}$
\end{prop}

\noindent In the preceding proposition, we assume that the homotopy equivalences preserve the quantum grading.\\
\ \\
\noindent We have seen in section \ref{sec:invariance} how to affect such a homotopy equivalence by simplifying the chain complex $(\leftComplex{T_{1}}, d_{APS})$. A similar result holds for the type D structure on $\leftComplex{T_{2}}$: simplifications of the chain complex with differential $d_{APS}$ results in a homotopy equivalent type $D$ structure on the simplified complex. Over a field, $\mathbb{F}$, such simplifications show that $(\leftComplex{T_{1}} \otimes \mathbb{F}, d_{APS}) \simeq H_{\ast, \mathbb{F}}( \leftComplex{T_{1}})$ where the homology is taken with respect to $d_{1}$ and similarly for $(\rightComplex{T_{2}} \otimes \mathbb{F}, d_{APS})$. These homologies are determined by the tangle homology of  Asaeda, Przytycki, and  Sikora.  Consequently,

\begin{cor}
There is a type $A$ structure on $H_{\ast, \mathbb{F}}( \leftComplex{T_{1}})$ and a type $D$ structure on $H_{\ast, \mathbb{F}}( \rightComplex{T_{2}})$ for which $$\complex{T} \simeq H_{\ast, \mathbb{F}}( \leftComplex{T_{1}}) \boxtimes H_{\ast, \mathbb{F}}( \rightComplex{T_{2}})$$ 
\end{cor}

\noindent For example, this result applies to the rational coefficient theory and the theory over $\Z/2\Z$. 
\section{Examples of pairing type A structures and type D structures}\label{sec:expairing}
\ \\
\noindent{\bf Example I:}({\em Reidemeister Invariance of Khovanov Homology}) Suppose that $L$ and $L'$ are two link diagrams for an oriented link in $S^{3}$. Furthermore, suppose they differ by a Reidemeister move. If $D^{2} \subset \R^{2}$ is the local region in which the move occurs, we can use $\partial D^{2}$ and the orientation on $\R^{2}$ to think of $R = L \cap D^{2}$ as an inside tangle, and $\righty{L}$ as an outside tangle. Then $\complex{L} \cong \leftComplex{R} \boxtimes \rightComplex{\righty{L}}$. If we let $R' = L' \cap D^{2}$ then $\complex{L'} \cong \leftComplex{R'} \boxtimes \rightComplex{\righty{L}}$.  In section \ref{sec:examplesA} we compute the type $A$ structure for three of the tangles involved in the Reidemeister moves. In each we saw that the the type $A$ structure was homotopy equivalent to the structure obtained for the tangle after applying the Reidemeister move. Due to the results in the appendix in section \ref{sec:boxtimes} this implies that 
$$\complex{L} \cong \leftComplex{R} \boxtimes \rightComplex{\righty{L}} \simeq \leftComplex{R'} \boxtimes \rightComplex{\righty{L}} \cong \complex{L'} $$
This gives a new perspective on the locality arguments for invariance in various forms of Khovanov homology.\\
\ \\
\noindent{\bf Example II:} In \cite{typeD} we computed the type $D$ structure $$\rightComplex{\righty{T_{L}}}$$ for the following tangle $\righty{T_{L}}$, based on the left handed trefoil,$T_{L}$,
$$
\inlinediag[0.3]{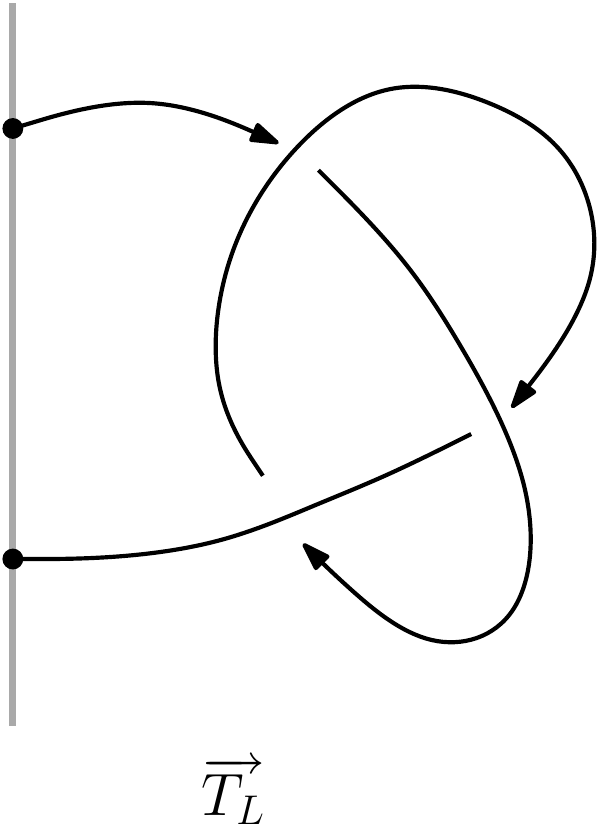}
$$
This structure is the map $\righty{\delta}$ where
\begin{equation}
\begin{split}
\righty{\delta}(s_{(-3,-15/2)}^{+}) &= 2\,\righty{e_{C}} \otimes s_{(-2,-13/2)}^{-} + \lefty{e_{C}} \otimes s_{(-3,-17/2)}^{-} \\
\righty{\delta}(s_{(-2,-11/2)}^{+}) &= - \lefty{e_{C}} \otimes s_{(-2,-13/2)}^{-}\\
\righty{\delta}(s_{(0,-3/2)}^{+})\ &= - \lefty{e_{C}} \otimes s_{(0,-5/2)}^{-}\\
\end{split}
\end{equation}
where the superscript indicates the decoration on the cleaved circle $C$ in the corresponding resolutions, and the subscript is the bigrading. We use the pairing theorem to compute several connect sums.\\
\ \\
\noindent {\em (i) With the unknot:} We can think of the unknot $U$ as a cleaved circle on $P_{2}$ where $\lefty{U} = U \cap \leftHalf$ and $\righty{U} = U \cap \rightHalf$ are the unique planar matchings. Then the left handed trefoil is the connect sum $U \# T_{L}$ which we can think of as gluing $\lefty{U}$ with the tangle above. The type $A$ structure $\leftComplex{\lefty{U}}$ is isomorphic to $\Z\,f_{(0,1/2)}\oplus \Z\,f_{(0,-1/2)}$, which is the idempotent decomposition for $I_{C^{+}}$ and $I_{C^{-}}$ (see the example in section 2). The action of $\lefty{e_{C}}$ is trivial since there are no crossings in the standard diagram. On the other hand $\righty{e_{C}}$ takes $f_{(0,1/2)}$ to $f_{(0,-1/2)}$. Since $d_{APS} \equiv 0$ for $\lefty{U}$, we need only compute $(m_{2} \otimes \I)(\I \otimes \righty{\delta})$. Using the idempotents we see that there are six generators. Furthermore, the only terms in $(m_{2} \otimes \I)(\I \otimes \righty{\delta})$ come from $\righty{e_{C}}$. This gives the following chain complex for $\leftComplex{\lefty{U}} \boxtimes \rightComplex{\righty{T_{L}}}$ 
$$
\begin{array}{lll}
f_{(0,1/2)} \otimes s_{(-3,-15/2)}^{+} &\stackrel{\cdot 2}{\longrightarrow} & f_{(0,-1/2)} \otimes s_{(-2,-13/2)}^{-} \\ 
f_{(0,1/2)} \otimes s_{(-2,-11/2)}^{+} & & \\
f_{(0,1/2)} \otimes s_{(0,-3/2)}^{+} & & \\
f_{(0,-1/2)} \otimes s_{(0,-5/2)}^{-} & & \\ 
f_{(0,-1/2)} \otimes s_{(-3,-17/2)}^{-} & & \\ 
\end{array}
$$
whose homology consists of a $\Z$-summand in bigradings $(-2,-5)$, $(0,-1)$, $(0,-3)$, $(-3,-9)$ and a $\Z/2\Z$ summand in bigrading $(-2,-7)$. This is the Khovanov homology of the left handed trefoil. \\
\ \\
\noindent {\em (ii) with the positive Hopf tangle in section \ref{sec:examplesA}:} To compute the Khovanov homology of the connect sum of the left handed trefoil with the Hopf link with $+1$ linking number. Recall that for the Hopf tangle with positive crossings we obtained for $\leftComplex{H}$ the following action $r_{(0,3/2)} \stackrel{\righty{e_{C}}}{\longrightarrow} r_{(0,1/2)}$ and $r_{(2,11/2)} \stackrel{\righty{e_{C}}}{\longrightarrow} r_{(2,9/2)}$ where the first entry in each corresponds to the $+$ decoration. Consequently, $\leftComplex{H} \otimes_{\mathcal{I}} \rightComplex{T_{L}}$ has the generators
$$
\begin{array}{llll}
r_{(0,3/2)}\otimes s_{(-3,-15/2)}^{+} &  r_{(2,11/2)}\otimes s_{(-3,-15/2)}^{+} &  r_{(0,3/2)} \otimes s_{(-2,-11/2)}^{+} &  r_{(2,11/2)} \otimes s_{(-2,-11/2)}^{+}\\ 
r_{(0,3/2)} \otimes s_{(0,-3/2)}^{+} & r_{(2,11/2)} \otimes s_{(0,-3/2)}^{+} & r_{(0,1/2)} \otimes s_{(-2,-13/2)}^{-} & r_{(2,9/2)} \otimes s_{(-2,-13/2)}^{-}\\
r_{(0,1/2)} \otimes s_{(-3,-17/2)}^{-} & r_{(2,9/2)} \otimes s_{(-3,-17/2)}^{-} & r_{(0,1/2)} \otimes s_{(0,-5/2)}^{-} & r_{(2,9/2)} \otimes s_{(0,-5/2)}^{-}\\
\end{array}
$$
Since $d_{APS} \equiv 0$ after the simplification, we need only compute $(m_{2} \otimes \I)(\I \otimes \righty{\delta})$ on these twelve generators. Since the action $m_{2}$ is trivial except for on $\righty{e_{C}}$, we can ignore all the terms in $\righty{\delta}$ except for those with $\righty{e_{C}}$. This leaves the following as the only non-trivial maps in $\partial^{\boxtimes}$:
$$
\begin{array}{l}
r_{(0,3/2)}\otimes s_{(-3,-15/2)}^{+} \longrightarrow 2\cdot(r_{(0,1/2)} \otimes s_{(-2,-13/2)}^{-})\\
r_{(2,11/2)}\otimes s_{(-3,-15/2)}^{+} \longrightarrow 2\cdot(r_{(2,9/2)} \otimes s_{(-2,-13/2)}^{-})
\end{array}
$$
Taking the homology of this new complex gives two $\Z/2\Z$ summands in bigradings $(-2,-6)$ and $(0,-2)$, and $8$ $\Z$ summands for the remaining generators in the corresponding bigrading:
$$
\begin{array}{ll}
r_{(0,3/2)} \otimes s_{(-2,-11/2)}^{+} \rightarrow (-2,-4) &  r_{(2,11/2)} \otimes s_{(-2,-11/2)}^{+} \rightarrow (0,0)\\ 
r_{(0,3/2)} \otimes s_{(0,-3/2)}^{+} \rightarrow (0,0) & r_{(2,11/2)} \otimes s_{(0,-3/2)}^{+}  \rightarrow (2,4)\\
r_{(0,1/2)} \otimes s_{(-3,-17/2)}^{-} \rightarrow (-3,-8) & r_{(2,9/2)} \otimes s_{(-3,-17/2)}^{-} \rightarrow (-1,-4) \\
r_{(0,1/2)} \otimes s_{(0,-5/2)}^{-} \rightarrow (0,-2) & r_{(2,9/2)} \otimes s_{(0,-5/2)}^{-} \rightarrow (2,2)\\
\end{array} 
$$
which is isomorphic to the Khovanov homology of the connect sum of the positive Hopf link with the left handed trefoil. \\
\ \\
\noindent {\em (iii) With a right-handed trefoil:} We give some of the details of the computation in the introduction. We consider the tangle $\lefty{T_{R}}$ found by removing an arc from the right-handed trefoil. We can compute $\leftComplex{\lefty{T_{R}}}$ directly. The result is a bigraded module spanned by $t^{+}_{(0,5/2)}$, $t^{+}_{(2,13/2)}$, $t^{+}_{(3,17/2)}$, $t^{-}_{(0,3/2)}$, $t^{-}_{(2,11/2)}$, $t^{-}_{(3,15/2)}$, where the superscript identifies the corresponding idempotent. The action of $\righty{e_{C}}$ is given by
$t^{+}_{(0,5/2)} \rightarrow t^{-}_{(0,3/2)}$, $t^{+}_{(2,13/2)} \rightarrow t^{-}_{(2,11/2)}$, $t^{+}_{(3,17/2)} \rightarrow t^{-}_{(3,15/2)}$. The action of $\lefty{e_{C}}$ is $t^{+}_{(2,13/2)} \rightarrow 2\cdot t^{-}_{(3,15/2)}$. \\
\ \\
\noindent Consequently, the module $\leftComplex{\lefty{T_{R}}} \otimes_{\mathcal{I}} \rightComplex{\righty{T_{L}}}$ has eighteen generators. Those, along with their images under $\partial^{\boxtimes}$ are shown in the following list. 
$$
\begin{array}{lll}
t^{+}_{(0,5/2)} \otimes s_{(-3,-15/2)}^{+} & \longrightarrow & 2\,t^{-}_{(0,3/2)}\otimes s_{(-2,-13/2)}^{-}  \\ 
t^{+}_{(0,5/2)} \otimes s_{(-2,-11/2)}^{+} & \longrightarrow & 0\\
t^{+}_{(0,5/2)} \otimes s_{(0,-3/2)}^{+} & \longrightarrow & 0 \\
t^{+}_{(2,13/2)} \otimes s_{(-3,-15/2)}^{+} & \longrightarrow & 2\,t^{-}_{(2,11/2)} \otimes s_{(-2,-13/2)}^{-} + 2\,t^{-}_{(3,15/2)} \otimes s_{(-3,-17/2)}^{-}\\
t^{+}_{(2,13/2)} \otimes s_{(-2,-11/2)}^{+} & \longrightarrow & -2\,t^{-}_{(3,15/2)} \otimes s_{(-2,-13/2)}^{-}\\
t^{+}_{(2,13/2)} \otimes s_{(0,-3/2)}^{+} & \longrightarrow & -2\,t^{-}_{(3,15/2)} \otimes s_{(0,-5/2)}^{-}\\
t^{+}_{(3,17/2)} \otimes s_{(-3,-15/2)}^{+} & \longrightarrow & 2\,t^{-}_{(3,15/2)}\otimes s_{(-2,-13/2)}^{-} \\
t^{+}_{(3,17/2)} \otimes s_{(-2,-11/2)}^{+}  & \longrightarrow & 0 \\
t^{+}_{(3,17/2)} \otimes s_{(0,-3/2)}^{+} & \longrightarrow & 0 \\
t^{-}_{(0,3/2)} \otimes s_{(-3,-17/2)}^{-} & \longrightarrow & 0\\
t^{-}_{(0,3/2)} \otimes s_{(-2,-13/2)}^{-} & \longrightarrow & 0\\
t^{-}_{(0,3/2)} \otimes s_{(0,-5/2)}^{-} & \longrightarrow & 0\\
t^{-}_{(2,11/2)} \otimes s_{(-3,-17/2)}^{-} & \longrightarrow & 0\\
t^{-}_{(2,11/2)} \otimes s_{(-2,-13/2)}^{-} & \longrightarrow & 0\\
t^{-}_{(2,11/2)} \otimes s_{(0,-5/2)}^{-} & \longrightarrow & 0\\
t^{-}_{(3,15/2)} \otimes s_{(-3,-17/2)}^{-} & \longrightarrow & 0\\
t^{-}_{(3,15/2)} \otimes s_{(-2,-13/2)}^{-} & \longrightarrow & 0\\
t^{-}_{(3,15/2)} \otimes s_{(0,-5/2)}^{-} & \longrightarrow & 0\\
\end{array}
$$    
From this we see immediately that there are $\Z$ summands for each of $t^{+}_{(0,5/2)} \otimes s_{(-2,-11/2)}^{+}$ in bigrading $(-2,-3)$, $t^{+}_{(0,5/2)} \otimes s_{(0,-3/2)}^{+}$ in $(0,1)$, $t^{+}_{(3,17/2)} \otimes s_{(-2,-11/2)}^{+}$ in $(1,3)$, $t^{+}_{(3,17/2)} \otimes s_{(0,-3/2)}^{+}$ in $(3,7)$, $t^{-}_{(0,3/2)} \otimes s_{(-3,-17/2)}^{-}$ in $(-3,-7)$, $t^{-}_{(0,3/2)} \otimes s_{(0,-5/2)}^{-}$ in $(0,-1)$, $t^{-}_{(2,11/2)} \otimes s_{(-3,-17/2)}^{-}$
in $(-1,-3)$, $t^{-}_{(2,11/2)} \otimes s_{(0,-5/2)}^{-}$ in $(2,3)$. \\
\ \\
\noindent The remaining generators occur in the non-zero rows for $\partial^{\boxtimes}$. We will have a $\Z/2\Z$-summand for $t^{-}_{(3,15/2)} \otimes s_{(0,-5/2)}^{-}$ in $(3,5)$ and for $t^{-}_{(0,3/2)}\otimes s_{(-2,-13/2)}^{-}$ in $(-2,-5)$. That $\partial^{\boxtimes}(t^{+}_{(2,13/2)} \otimes s_{(-3,-15/2)}^{+})$ equals $2\,t^{-}_{(2,11/2)} \otimes s_{(-2,-13/2)}^{-} + 2\,t^{-}_{(3,15/2)} \otimes s_{(-3,-17/2)}^{-}$ gives a $\Z \oplus \Z/2\Z$ in $(0,-1)$. In addition, that $\partial^{\boxtimes}(t^{+}_{(2,13/2)} \otimes s_{(-2,-11/2)}^{+})= -2\,t^{-}_{(3,15/2)} \otimes s_{(-2,-13/2)}^{-}$ and $\partial^{\boxtimes}(t^{+}_{(3,17/2)} \otimes s_{(-3,-15/2)}^{+})= 2\,t^{-}_{(3,15/2)}\otimes s_{(-2,-13/2)}^{-} $ means that $(t^{+}_{(2,13/2)} \otimes s_{(-2,-11/2)}^{+}) + (t^{+}_{(3,17/2)} \otimes s_{(-3,-15/2)}^{+})$ generates a $\Z$ summand in homology in bigrading $(0,1)$, while $t^{-}_{(3,15/2)}\otimes s_{(-2,-13/2)}^{-}$ generates a $\Z/2\Z$ summand in $(1,1)$.  \\
\ \\
\noindent Consequently, the Khovanov homology of this connected sum has free part
$$
\Z_{(-3,7)} \oplus \Z_{(-2,-3)} \oplus \Z_{(-1,-3)} \oplus \Z^{2}_{(0,-1)} \oplus \Z^{2}_{(0,1)} \oplus \Z_{(1,3)} \oplus \Z_{(2,3)} \oplus \Z_{(3,7)}
$$
while the torsion part is
$$
\big(\Z/2\Z)_{(-2,-5)} \oplus \big(\Z/2\Z)_{(0,-1)} \oplus \big(\Z/2\Z)_{(1,1)} \oplus \big(\Z/2\Z)_{(3,5)}
$$
This agrees with the Khovanov homology of the knot as computed by Bar-Natan and Greene's JavaKH program. Note that we have correctly computed the torsion terms. In Khovanov's original paper the connect sum gives rise to a long exact sequence, which in general is not enough to compute the homology due to the usual ambiguity in long exact sequences. However, our approach will compute the torsion correctly, and provide a modular approach so that previous computations may be reused. 

\appendix
\section{Graded modules and conventions}
\noindent Let $M$ be a $\Z$-graded module over a ring $R$ and let $M_{i}$ be the module of elements in grading $i \in \Z$. For a homogeneous element $m \in M$, $|m|$ will denote the grading of $m$: if $m \in M_{i}$ then $|m| = i$.\\
\ \\
\noindent A module map $f: M \rightarrow M'$ has {\em order} $r$ if the composition $M_{i} \hookrightarrow M \stackrel{f}{\longrightarrow} M'$ has image in $M_{i+r}$ for each $i \in \Z$. \\
\ \\
\noindent {\em Degree shift convention:} If $M$ is a $\Z$-graded module, $M[n]$ is the graded module with $(M[n])_{i} = M_{i-n}$, i.e. the module found by shifting the homogeneous elements of $M$ up $n$ levels. If $m \in M$, the corresponding element in $M[n]$ will be denoted $m[n]$. Thus $|m[n]| = |m| + n$. \\
\ \\
\noindent An order $r$ map $f: M \rightarrow M'$ induces order $0$ maps $M \rightarrow M'[-r]$ and $M[r] \rightarrow M'$, along with maps of different orders $M[n] \rightarrow M[s]$. These will also be denoted by $f$, except where confusion could arise.\\ 

\noindent The identity on $M$ will be denoted $\I_{M}$. We will also have need of a graded version of the identity:

\begin{defn}
$|\I_{M}| : M \rightarrow M$ is the $0$-order map defined by
setting 
$$
|\I_{M}|(m) = (-1)^{|m|} m
$$
for homogeneous $m \in M$ and linearly extending to $M$. $|\I_{M}|^{j}$ is the $j$-fold composition of $|\I_{M}|$. 
\end{defn}

\noindent If $m \in M_{i}$ then $|\I_{M}|^{j}(m) = (-1)^{ij} m$. Consequently, $|\I_{M}|^{j} \circ |\I_{M}|^{k} =  |\I_{M}|^{j+k}$ while $\big(|\I_{M}|^{j}\big)^{k} = |\I_{M}|^{jk}$.  \\
\ \\
In addition, shifting changes the sign:
$$|\I_{M[n]}| = (-1)^{n} |\I_{M}|$$ and $|\I_{M[n]}|^{j} = (-1)^{jn} |\I_{M}|$.\\

\subsection{Tensor Algebras}

\noindent We fix a $\Z$-graded $R$-module $A$. As usual,  
 $$\mathcal{T}^{\ast}(A) = \bigoplus_{i = 0}^{\infty} A^{\otimes n}$$
where $A^{\otimes 0} = R$ and for $n > 0$, $A^{\otimes n} = A \otimes_{R} A \otimes_{R} \cdots \otimes_{R} A$ using exactly $n$ factors. $A^{\otimes n}$ is graded using the standard rule $|a_{1}\otimes \cdots \otimes a_{n}| = \sum |a_{i}|$.\\
\ \\
\noindent Furthermore, $\mathcal{T}^{\ast}(A)$ has a filtration 
$$
R \subset \mathcal{T}^{1}(A) \subset \cdots \subset \mathcal{T}^{k}(A) \subset \cdots
$$
where $\mathcal{T}^{k}(A) = \bigoplus_{i = 0}^{k} A^{\otimes n}$. \\
\ \\ 
\noindent By $\I_{A}^{\otimes n}$ we will mean the identity on $A^{\otimes n}$  thought of as the map $\I_{A}\!\otimes\!\I_{A}\!\otimes\!\cdot\!\otimes\!\I_{A}$. In general, we will only use the subscript when we need to distinguish $A$; by default, $\I^{\otimes n}$ will be the identity on $A^{\otimes n}$. Furthermore, by $|\I|^{j \otimes n}$ we will mean the map  $|\I|^{j}\otimes \cdots \otimes |\I|^{j}$ on $A^{\otimes n}$. 

\begin{defn}
For any $\Z$-graded module, $\mathcal{T}^{\ast}_{A}(M)$ is the $\Z$-graded $R$-module $M \otimes_{R} \mathcal{T}^{\ast}(A)$ filtered by the submodules $M \otimes \mathcal{T}^{k}(A)$ for $k = 0,1,2, \ldots$.
\end{defn}

\begin{defn}
Let $\mathfrak{T}_{A}$ be the category whose objects are the $R$-modules $\mathcal{T}^{\ast}_{A}(M)$ for each $\Z$-graded module $M$, and whose morphisms, $\mathcal{T}_{A}(M,M')$, are filtered $R$-module maps $\Phi : \mathcal{T}^{\ast}_{A}(M) \rightarrow \mathcal{T}^{\ast}_{A}(M')$. 
\end{defn}

\begin{defn}
Let $\Phi \in \mathcal{T}_{A}(M,M')$. For $1, j \in \N$, the $ij^{th}$ {\em component} of $\Phi$ is the map
$$
\Phi_{ij}:  M \otimes A^{\otimes(i\!-\!1)} \hookrightarrow \mathcal{T}^{\ast}_{A}(M) \stackrel{\Phi^{\ast}}{\longrightarrow} \mathcal{T}^{\ast}(M') \longrightarrow M' \otimes A^{\otimes(j\!-\!1)}
$$
\end{defn}
\noindent Since $\Phi$ is filtered, $\Phi_{ij} = 0$ unless $j \leq i$.

\subsection{The $\infty$-sub-category of $\mathcal{T}_{A}^{\ast}$} 

\begin{prop}\label{prop:cat}
Let $\mathcal{C}_{A}(M,M') \subset \mathcal{T}_{A}(M,M')$ be those module maps $\Phi: \mathcal{T}^{\ast}_{A}(M) \rightarrow \mathcal{T}^{\ast}_{A}(M')$ such that $\Phi$ has order $r$ for some $r \in \Z$, and 
\begin{equation}\label{eqn:criterion1}
\Phi_{nm}= \Phi_{n-m+1,1} \otimes |\I|^{(n\!+\!m\!+\!r)\otimes (m\!-\!1)}
\end{equation}
for every $n,m \in \N$ with $1 \leq m \leq n$. Then $\mathcal{C}_{A}(M,M')$ are the sets of morphisms for a full subcategory of $\mathcal{T}_{A}$
\end{prop}

\noindent{\bf Proof:} First, we verify that $\I_{\mathcal{T}^{\ast}_{A}(M)} \in \mathcal{C}_{A}(M,M)$. $I_{nm}$ is non-zero only if $n=m$. When $n=m$ the right side of \ref{eqn:criterion1} equals $\I_{11} \otimes |\I|^{(n\!+\!n\!+\!0)\otimes (n\!-\!1)}$. However, $I_{11} = \I_{M}$, and $|\I|^{(n\!+\!n\!+\!0)\otimes (n\!-\!1)} = \I^{\otimes (l\!-\!1)}$ since an even entry in the first place in the exponent of $|\I|$ will not change the sign. Thus, $I_{nn}= \I_{M} \otimes \I^{\otimes (n\!-\!1)}$, which is the identity on $M \otimes A^{\otimes(n-1)}$. On the other hand, if $n > m$, then $I_{nm}\equiv 0$ and $I_{n-m+1,1} \equiv 0$ as well. Thus, $\I_{\mathcal{T}^{\ast}_{A}(M)} \in \mathcal{C}_{A}(M,M)$ for every $M$. \\
\ \\
\noindent We now need to verify that composition of morphisms with the property in \ref{eqn:criterion1} will still have this property. Suppose $\Phi \in \mathcal{C}_{A}(M,M')$ has order $r$ and $\Psi \in \mathcal{C}_{A}(M',M'')$ has order $s$, and components of each satisfy \ref{eqn:criterion1}, then the order $(r+s)$ morphism $\Psi \circ \Phi$ has components
\begin{equation}
\begin{split}
(\Psi\circ\Phi)_{nm} &= \sum_{m \leq k \leq n} \Psi_{km} \circ \Phi_{nk}\\
&= \sum_{m \leq k \leq n}\big(\Psi_{k-m+1,1} \otimes |\I|^{(k\!+\!m\!+\!s)\otimes (m\!-\!1)}\big) \circ \big(\Phi_{n-k+1,1} \otimes |\I|^{(n\!+\!k\!+\!r)\otimes (k\!-\!1)}\big)\\
&= \sum_{m \leq k \leq n}\Psi_{k-m+1,1}(\Phi_{n-k+1,1} \otimes |\I|^{(n\!+\!k\!+\!r)\otimes (k\!-\!m)}\big)\otimes |\I|^{(n\!+\!m\!+\!r+\!s)\otimes (m\!-\!1)}\\
\end{split}
\end{equation}
If we let $i = k - m + 1$ and $j = n - k + 1$, then  $n + k \equiv j + 1$ modulo 2, as $k$ changes from $m$ to $n$, 
$i$ changes from $1$ to $n-m+1$, so we can rewrite the previous result as
\begin{equation}
(\Psi\circ\Phi)_{nm} = \big(\sum_{i+j=n-m+2} \Psi_{i,1}(\Phi_{j,1} \otimes |\I|^{(j\!+\!r\!+\!1)\otimes(i\!-\!1)})\big)
\otimes \otimes |\I|^{(n\!+\!m\!+\!(r+s))\otimes (m\!-\!1)}\big)\\
\end{equation}
On the other hand, 
\begin{equation}\label{eqn:identity}
\begin{split}
(\Psi\circ\Phi)_{n-m+1,1} &= \sum_{1 \leq i \leq n-m+1} \Psi_{i,1} \circ \Phi_{n-m+1,i}\\
&= \sum_{1 \leq i \leq n-m+1} \Psi_{i,1}(\Phi_{n-m-i+2,1} \otimes  |\I|^{(n\!+\!m\!+\!1\!+\!i\!+\!r)\otimes (i\!-\!1)}\big)\\
&= \sum_{i+j=n-m+2} \Psi_{i,1}(\Phi_{j,1} \otimes |\I|^{(j\!+\!r\!+\!1)\otimes (i\!-\!1)}\big)\\
\end{split}
\end{equation}
when we let $j=n-m+2-i$. Thus $\Psi \circ \Phi$ satisfies \ref{eqn:criterion1}.$\Diamond$ \\

\begin{defn}
For $\Phi \in \mathcal{C}(M,M')$ of order $r$ the {\em core} of $\Phi$ is the set of order $r$ module maps $\Phi^{\ast} = \{\,\phi_{k}\,|\,n \in \N\,\}$ where $\phi_{k} = \Phi_{k1} : M \otimes A^{k-1} \rightarrow M'$. Given a set of order $r$ module maps $R = \{\,\rho_{k}\,|\,n \in \N\,\}$ with $\rho_{k}: M \otimes A^{k-1} \rightarrow M'$ the {\em extension} of $R$ is the map $\overline{R} \in \mathcal{C}(M,M')$ with with components
\begin{equation}\label{eqn:extension}
\overline{R}_{nm}= \rho_{n-m+1} \otimes |\I|^{(n\!+\!m\!+\!r)\otimes (m\!-\!1)}
\end{equation}
\end{defn}

\noindent The argument in proposition \ref{prop:cat} shows that these are inverses: for $\Phi \in \mathcal{C}(M,M')$, $\overline{\Phi^{\ast}} = \Phi$ while for $R = \{\,\rho_{k}\,|\,n \in \N\,\}$, $\big(\overline{R}\big)^{\ast}$ equals $R$. Consequently, we can describe $\mathcal{C}_{A}$ completely in terms of a composition on the cores $\Phi^{\ast}$ which directly reflects the usual module map composition for filtered maps on $\mathcal{T}^{\ast}_{A}(M)$. This allows us to pull the operations of $\mathcal{C}_{A}$ back to the category of $R$-modules. 

\begin{defn}
$\mathcal{C}^{\ast}_{A}$ is the category whose objects are $\Z$-graded $R$-modules, and whose morphisms $\Phi^{\ast}: M {\ast\!\!\!\!\rightarrow} M'$ are sets $\Phi^{\ast}=\{\,\phi_{i}\,|\,i \in \N\,\}$ of $R$-module maps $\phi_{i}: M \otimes A^{\otimes (i\!-\!1)} \rightarrow M'$ such that  every $\phi_{j}$ has order $r$ for some $r \in \Z$.  The identity $\I^{\ast}_{M} : M {\ast\!\!\!\!\rightarrow} M$ is the set of $0$-order module maps with $(\I^{\ast}_{M})_{1} = \I_{M}$ and $(\I^{\ast}_{M_{i}} = 0$ for $i > 1$. The composition of an order $r$ morphism $\Phi^{\ast}: M {\ast\!\!\!\!\rightarrow} M'$ with an order $s$ morphism $\Psi^{\ast}: M' {\ast\!\!\!\!\rightarrow} M''$ is the set of order $r+s$ module maps given by
$$\big( \Psi^{\ast}\!\ast\! \Phi^{\ast} \big)_{k} = \sum_{i+j=k+1} \psi_{i}(\phi_{j} \otimes |\I|^{(j\!+\!r\!+\!1)\otimes(i\!-\!1)}) $$
for $k = 1,2,\ldots$.
\end{defn}

\noindent Proposition \ref{prop:cat} implies

\begin{prop}
There is a functor $\mathcal{F}: \mathcal{C}^{\ast}_{A} \rightarrow \mathcal{C}_{A}$ which takes $M  \rightarrow \mathcal{T}^{\ast}_{A}(M)$ and $\Phi^{\ast}: M {\ast\!\!\!\!\rightarrow} M'$ to 
its extension $\Phi: \mathcal{T}^{\ast}_{A}(M) \rightarrow \mathcal{T}^{\ast}_{A}(M')$. 
\end{prop}

\noindent We will generally work in $\mathcal{C}_{A}$, and then pull back our results to $\mathcal{C}^{\ast}_{A}$. 

\subsection{$\infty$-structures}

\noindent Since $A$ is a $\Z$-graded $R$-module, we may take $M = A$ above. Then $\mathcal{T}_{A}^{\ast}(A)=$ $A \otimes \mathcal{T}^{\ast}(A) \cong $ $\bigoplus_{n=1}^{\infty} A^{\otimes n}$. Let $P: \mathcal{T}_{A}(A) \rightarrow \mathcal{T}_{A}^{\ast}(A)$ be an order $r$ map in $\mathcal{T}_{A}(A,A)$. Then we can form a new map $P + (\I_{A} \otimes P)$ in $\mathcal{T}_{A}(A,A)$. This is evidently still filtered. 

\begin{defn}
An $\infty$-algebra structure on $A$ is an order $1$ map $D \in \mathcal{T}(A,A)$ such that 
\begin{enumerate}
\item $D \circ D = 0$, and 
\item $D +  (\I \otimes D)$ is in $\mathcal{C}(A,A)$
\end{enumerate}  
\end{defn} 

For an $\infty$-algebra structure $D$, we will let $\mu = D +  \I \otimes D$. Then the core of $\mu$ is a collection of maps $\mu^{\ast} = \{\mu_{i}: A\otimes \mathcal{A}^{\otimes (i-1)} \rightarrow A\}$ in $\mathcal{C}^{\ast}_{A}(A,A)$. When we have a prescribe $\mu$ in mind, we will write $D_{\mu}$ for the corresponding structure. 

\begin{defn}
A right $\infty$-module $M$ over $(A,D_{\mu})$ is a $\Z$-graded module and an order $1$ morphism $D_{M} \in \mathcal{T}(M,M)$ such that 
\begin{enumerate}
\item $D_{M} \circ D_{M} = 0$, and 
\item $D_{M} + \I \otimes D_{\mu}$ is in $\mathcal{C}(M,M)$.
\end{enumerate}
\end{defn}

\noindent Notice that $A$ with the map $D_{\mu}$ is a right module over $(A,D_{\mu})$. A right $\infty$-module over $(A,D_{\mu})$ is a chain complex with an additional requirement placed on its boundary map.  We can similarly adapt the notion of chain map and chain homotopy to this context.

\begin{defn}
An $\infty$-module map between right $\infty$-modules $(M, D_{M})$ and $(M', D_{M'})$ (over $(A,D_{\mu})$) is an order $0$ morphism $\Psi \in \mathcal{C}(M,M')$ such that $\Psi \circ D_{M} = D_{M'} \circ \Psi$
\end{defn}

\begin{defn}
An $\infty$-homotopy between $\infty$-module maps $\Phi$ and $\Psi$, each mapping $(M,D_{M})$ to $(M', D_{M'})$, is an order $-1$ map $H \in \mathcal{C}(M,M')$ such that $\Phi - \Psi = H \circ D_{M} + D_{M'} \circ H$. 
\end{defn}

\noindent Since chain complexes form a category, and the maps are drawn from the morphisms of $\mathcal{C}(M,M')$, we obtain a category of right $\infty$-modules. Furthermore, we can quotient by chain homotopies to obtain a notion of chain homotopy equivalence. 

\subsection{$\infty$-structures in terms of the core category}

\noindent Let $(A,D_{\mu})$ be an $\infty$-algebra. The following identity is an immediate consequence of the definition
$$
D_{\mu} = \mu - (\I \otimes D_{\mu})
$$
From this identity we obtain 
$$
D_{\mu} = \mu - \I \otimes \mu + \I\otimes \I \otimes \mu + \cdots = \sum_{l=0}^{\infty} (-1)^{l}\big(\I^{\otimes l} \otimes \mu\big)
$$
Note that the sum is actually finite on any summand $A^{\otimes n}$. \\
\ \\
\noindent If we wish to write out the relations for $\infty$-algebras, modules, morphisms, etc. in terms of their cores we encounter the difficulty that $D_{\mu}$ is not itself in $\mathcal{C}_{A}(A,A)$; furthermore composing with it is not likely to be in $\mathcal{C}_{A}(A,A)$ either. However, we a graded commutator with $I \otimes D_{\mu}$ will be an extension. Before we prove this, we must understand commutators  with $|\I|^{j}$:

\begin{prop}
Let $R \in \mathcal{C}_{A}(A,A)$ have order $r$, then $|\I_{A}|^{j} \circ R_{k,1} = (-1)^{rj}R_{k,1} \circ |\I|^{j\otimes k}$
\end{prop}

\begin{prop}
Let $\Phi \in \mathcal{C}_{A}(M,M')$ have order $r$. Then $\Phi(\I\otimes D)-(-1)^{r}(\I \otimes D)\Phi$ is in $\mathcal{C}_{A}(M,M')$ and has core $\big\{\,\big(\Phi(\I \otimes D)\big)_{n,1}\,|\,n\in\N\,\big\}$. 
\end{prop}

\noindent{\bf Proof:} We compute the components of  $(\I \otimes D)\Phi$: 
\begin{equation}
\begin{split}
\big((\I \otimes &D)\Phi\big)_{lm} = \sum_{m \leq k \leq l} (\I \otimes D)_{km} \circ \Phi_{lk}\\
%&= \sum_{m \leq k \leq l}(\I \otimes D)_{km} \big(\Phi^{\ast}_{l-k+1} \otimes |\I|^{(l\!+\!k\!+\!r)\otimes (k\!-\!1)}\big)\\
&= \sum_{m \leq k \leq l}\big(\sum_{s=0}^{\infty}(-1)^{s} \I \otimes \I^{s} \otimes \mu_{k-s-1,m-s-1}\big)\circ \big(\Phi_{l-k+1,1} \otimes |\I|^{(l\!+\!k\!+\!r)\otimes (k\!-\!1)}\big)\\
&= \sum_{m \leq k \leq l}\big(\sum_{s=0}^{\infty}(-1)^{s} \I \otimes \I^{s} \otimes \mu_{k-m+1,1} \otimes |\I|^{(k\!+\!m\!+\!1)\otimes(m\!-\!s\!-\!2)}\big)\circ \big(\Phi_{l-k+1,1} \otimes |\I|^{(l\!+\!k\!+\!r)\otimes (k\!-\!1)}\big)\\
&= \sum_{m \leq k \leq l}\big(\sum_{s=0}^{\infty}(-1)^{s\!+\!l\!+\!k\!+\!r}\big(\Phi_{l-k+1,1} \otimes |\I|^{(l\!+\!k\!+\!r)\otimes (m\!-\!1)}\big)( \I \otimes \I^{l-k+s} \otimes \mu_{k-m+1,1} \otimes |\I|^{(k\!+\!m\!+\!1)\otimes(m\!-\!s\!-\!2)}\big)\big)\\
&= \sum_{m \leq k \leq l}\big(\sum_{s'=l-k}^{\infty}(-1)^{s'\!+\!r}\big(\Phi_{l-k+1,1} \otimes |\I|^{(l\!+\!k\!+\!r)\otimes (m\!-\!1)}\big)( \I \otimes \I^{s'} \otimes \mu_{k-m+1,1} \otimes |\I|^{(k\!+\!m\!+\!1)\otimes(m\!+\!l\!-\!s'\!-\!k\!-\!2)}\big)\big)\\
&= (-1)^{r} \sum_{m \leq k \leq l}\big(\Phi_{l-k+1,1} \otimes |\I|^{(l\!+\!k\!+\!r)\otimes (m\!-\!1)}\big)\big(\sum_{s'=l-k}^{\infty}(-1)^{s'} \I \otimes \I^{s'} \otimes \mu_{k-m+1,1} \otimes |\I|^{(k\!+\!m\!+\!1)\otimes(m\!+\!l\!-\!s'\!-\!k\!-\!2)}\big)\\
%&= (-1)^{r} \sum_{m \leq k \leq l}\Phi_{l+m-k,m}\big(\sum_{s'=l-k}^{\infty}(-1)^{s'} \I \otimes \I^{s'} \otimes \mu_{l-s'-1,l+m-k-s'-1}\big)\\
&= - (-1)^{r} \sum_{m \leq k \leq l} \big(\Phi_{l-k+1,1} \otimes |\I|^{(l\!+\!k\!+\!r)\otimes (m\!-\!1)}\big)\left(\sum_{s'=0}^{l-k-1}(-1)^{s'} \I \otimes \I^{s'} \otimes \mu_{k-m+1,1} \otimes |\I|^{(k\!+\!m\!+\!1)\otimes(m\!+\!l\!-\!s'\!-\!k\!-\!2)}\right)\\
&\hspace{0.5in} (-1)^{r} \sum_{m \leq k \leq l}\Phi_{l+m-k,m}\big(\I \otimes D\big)_{l,l+m-k} \\
&= -(-1)^{r}\left[ \sum_{m \leq k \leq l}\Phi_{l-k+1,1}\big(\sum_{s'=0}^{l-k-1}(-1)^{s'}(\I \otimes \I^{s'} \otimes \mu_{k-m+1,1} \otimes |\I|^{(k\!+\!m\!+\!1)\otimes(l\!-\!s'\!-\!k\!-\!1)}) \big) \otimes |\I|^{(l\!+\!m\!+\!r\!+\!1)\otimes (m\!-\!1)} \right]\\
&\hspace{0.5in} + (-1)^{r} \big(\Phi(\I\otimes D)\big)_{l,m} \\
\end{split}
\end{equation}
However, $\Phi_{l-k+1,1}$ will consume the $M$ factor as well as the first $l-k$ factors of $A$. Since $s'$ only has range up to $l-k-1$,the $\mu$ term must feed into an argument of $\Phi_{l-k+1,1}$. By the identity \ref{eqn:extension}, $(-1)^{r}(\I\otimes D)\Phi - \Phi(\I\otimes D)$ is the extension of $\{\rho_{n}\}$ where
\begin{equation}
\begin{split}
\rho_{n} &= -\sum_{1 \leq k \leq n}\Phi_{n-k+1,1}\big(\sum_{s=0}^{n-k-1}(-1)^{s}(\I \otimes \I^{s} \otimes \mu_{k,1} \otimes |\I|^{k \otimes(n\!-\!s\!-\!k\!-\!1)}) \big)\\
&= \sum_{1 \leq k \leq n}\Phi_{n-k+1,1}\big(\sum_{s=0}^{n-k-1}(-1)^{s}(\I \otimes \I^{s} \otimes \mu_{n-s-1,n-s-k}) \big)
\end{split}
\end{equation}
Thus $\rho_{n}=(-\Phi(\I \otimes D))_{n1}$. $\Diamond$\\
\ \\
\noindent To write out the requirement that $(A,D_{\mu})$ be an $\infty$-algebra in terms of the $\mu_{i}$, we first note
that $0 = D_{\mu} \circ D_{\mu} = \mu \circ \mu - (\I \otimes D_{\mu})\mu - \mu(\I \otimes D_{\mu})$. Using preceding proposition, we see that $0 = \mu \circ \mu - \overline{\{\big(\mu(\I\otimes D_{\mu})\big)_{n1}\}}$. However, $\mu \circ \mu = \overline{\mu^{\ast} \ast \mu^{\ast}}$, so for each $n \in \N$, $(\mu \ast \mu)_{n} - \big(\mu(\I\otimes D_{\mu})\big)_{n,1} = 0$. Unpacking the definitions above produces the following:

$$
\sum_{\footnotesize\begin{array}{c} i\!+\!j\!=\!n\!+\!1 \\ 1 \leq l \leq i \end{array}\normalsize} \mu_{i}(\mu_{j} \otimes |\I|^{j \otimes (i-l)} ) - \sum_{\footnotesize\begin{array}{c} i\!+\!j\!=\!n\!+\!1 \\ 0 \leq l \leq i-2 \end{array}\normalsize} (-1)^{l}\mu_{i}(\I \otimes \I^{\otimes l} \otimes \mu_{j} \otimes |\I|^{j \otimes (i-l-2)} ) = 0  
$$ 

\noindent which is equivalent to 

$$
\sum_{\footnotesize\begin{array}{c} i\!+\!j\!=\!n\!+\!1 \\ 1 \leq l \leq i \end{array}\normalsize} (-1)^{l+1}\mu_{i}(\I^{\otimes (l-1)} \otimes \mu_{j} \otimes |\I|^{j \otimes (i-l)} ) = 0  
$$ 

\begin{defn}
An $A_{\infty}$-algebra structure on a $\Z$-graded $R$-module $A$ is an $\infty$-algebra structure $D_{\mu}$ on $A[-1]$
\end{defn}

\noindent If $\mu^{\ast}=\{\mu_{i}\}$, then $\mu_{i}: (A[-1])^{\otimes i} \rightarrow A[-1]$ being order $1$ means $\mu_{i} : \big((A[-1])^{\otimes i}\big)_{k} \rightarrow (A[-1])_{k+1}$ or $(A^{\otimes i})_{k+i} \rightarrow A_{k+2}$. If we let $k' = k+i$ then $\mu_{i}: A^{\otimes i}_{k'} \rightarrow A_{k'-i+2} = A[i-2]_{k'}$. Thus, in terms of $A$ with its original grading, each $\mu_{n}$ needs to be an order $2-n$ map $A^{\otimes n} \rightarrow A$. Alternatively, $\mu_{n}$ is a grading preserving map $A^{\otimes n} \rightarrow  A[n-2]$.\\
\ \\
\noindent The preceding relation for an $\infty$-structure on $A[-1]$ is
$$
\sum_{\footnotesize\begin{array}{c} i\!+\!j\!=\!n\!+\!1 \\ 1 \leq l \leq i \end{array}\normalsize} (-1)^{l+1}\mu_{i}(\I^{\otimes (l-1)} \otimes \mu_{j} \otimes |\I|_{A[-1]}^{j \otimes (i-l)} ) = 0  
$$ 

\noindent Since $|\I|_{A[-1]} = - |\I|_{A}$ in terms of the grading on $A$ we obtain

$$
\sum_{\footnotesize\begin{array}{c} i\!+\!j\!=\!n\!+\!1 \\ 1 \leq l \leq i \end{array}\normalsize} (-1)^{j(i-l) + l+1}\mu_{i}(\I^{\otimes (l-1)} \otimes \mu_{j} \otimes |\I|_{A}^{j \otimes (i-l)} ) = 0  
$$ 

\noindent However, $j(i-l) + l + 1 \equiv ji + lj + l + 1 \equiv j(i+1)+(l+1)(j+1)$ so

$$
\sum_{\footnotesize\begin{array}{c} i\!+\!j\!=\!n\!+\!1 \\ 1 \leq l \leq i \end{array}\normalsize} (-1)^{j(i+1) + (j+1)(l+1)}\mu_{i}(\I^{\otimes (l-1)} \otimes \mu_{j} \otimes |\I|_{A}^{j \otimes (i-l)} ) = 0  
$$ 

\noindent We therefore see that our definition of an $A_{\infty}$-algebra is equivalent to the more standard

\begin{defn}
An $A_{\infty}$-algebra $A$ over a ring $R$ is a $\Z$-graded $R$-module $A$ equipped with maps $\mu_{n}:A^{\otimes n} \rightarrow A[n-2]$ for each $n \in \N$, which satisfy the relation
$$
0 = \sum_{\footnotesize\begin{array}{c} i\!+\!j\!=\!n\!+\!1 \\ l\!\in\!\{1,\!\ldots\!,\!i\!\} \end{array}\normalsize} (-1)^{j(i+1)+(j+1)(l+1)}\mu_{i}(\I^{\otimes (l-1)} \otimes \mu_{j} \otimes |\I|^{j \otimes (i-l)} )  
$$
\end{defn} 

\noindent Similarly, if we define 

\begin{defn}
A (right) $A_{\infty}$-module structure on a $\Z$-graded $R$-module $M$, over an $A_{\infty}$ algebra $(A,\mu)$, is an right $\infty$-module structure $D_{M[-1]}$ on $M[-1]$ over $(A[-1],D_{\mu})$.
\end{defn}

\noindent Following the argument above, suppose the core of $D + (\I \otimes D)$ is a set of order $1$ maps $m_{i} : M[-1] \otimes (A[-1])^{\otimes (i-1)} \rightarrow M[-1]$ with extension $\overline{m}$.  Then $D \circ D \equiv 0$ is equivalent to $\overline{m}\circ \overline{m} - (\I_{A[-1]} \otimes D_{\mu})\overline{m} - \overline{m}(\I_{A[-1]} \otimes D_{\mu}) \equiv 0$. Pushing this identity back to one involving the maps $m_{i} : M \otimes A^{\otimes(i-1)} \rightarrow M[i-2]$ yields 

\begin{defn}[\cite{Bor1}]
A right $A_{\infty}$-module $M$ over an $A_{\infty}$-algebra $A$ is a set of maps $\{m_{i}\}_{i\in \N}$ with
$m_{i}: M \otimes A^{\otimes (i-1)} \rightarrow M[i-2]$, and satisfying the following relations for each $n \geq 1$:
\begin{equation}
\begin{split}
0 = \sum_{i+j=n+1} (-1)^{j(i+1)}m_{i}(m_{j}& \otimes |\I|^{j \otimes (i-1))})\\
& + \sum_{i+j=n+1,k>0}(-1)^{k(j+1)+j(i+1)} m_{i}(\I^{\otimes\,k} \otimes \mu_{j} \otimes |\I|^{j \otimes (i-k-1)})
\end{split}
\end{equation}
$M$ is said to be strictly unital if for any $\xi \in M$,  $m_{2}(\xi \otimes \I_{A}) = \xi$, but for $n > 1$, $m_{n}(\xi \otimes a_{1} \otimes a_{2} \otimes \cdots \otimes a_{n-1}) = 0$ if any $a_{i} = \I_{A}$.
\end{defn} 

\noindent The definition for an $\infty$-morphism unpacks similarly:

\begin{defn}
An $A_{\infty}$ morphism $\Psi$ from $M$ to $M'$ over $(A,\mu)$, is an $\infty$-morphism from $(M[-1],D_{M[-1]})$ to $(M'[-1], D_{M'[-1]})$ over $(A[-1],D_{\mu})$. 
\end{defn}

\noindent The same argument as above allows us to write this requirement in terms of the core maps for $\Psi$, conceived of as order $0$ module maps $\psi_{i}: M[-1] \otimes (A[-1])^{\otimes\,(i-1)} \longrightarrow M'[-1]$. The requirement that $\Psi \circ D_{M[-1]} = D_{M'[-1]} \circ \Psi$ becomes $\Psi \circ \overline{m} - \overline{m'}\circ\Psi = \Psi \circ (\I \otimes D_{\mu}) - (\I \otimes D_{\mu}) \circ \Psi$. By our proposition, $\Psi \circ (\I \otimes D_{\mu}) - (\I \otimes D_{\mu}) \circ \Psi$ is the extension of $\big\{(\Psi(\I \otimes D_{\mu}))_{n1}\big\}$, since $\Psi$ is order $0$. Writing out his relation in terms of the cores, and adjusting $|\I|_{A[-1]} = - |\I|_{A}$ as above yields the standard definition:

\begin{defn}[\cite{Bor1}]
An $A_{\infty}$-morphism $\Psi$ of right $A$-modules $M$ and $M'$ is a set of maps
$\psi_{i}: M \otimes A^{\otimes\,(i-1)} \longrightarrow M'[i-1]$ for $i \in \N$,  satisfying
\begin{equation}
\begin{split}
\sum_{i+j=n+1}(-1)^{(i+1)(j+1)}m'_{i}(\psi_{j}&\otimes |\I|^{(j+1)\otimes(i-1)}) =\\
&= \sum_{i+j=n+1}(-1)^{j(i+1)}\psi_{i}(m_{j} \otimes |\I|^{j \otimes(i-1)})\\
& \hspace{.5in} + \sum_{i+j=n+1,k>0} (-1)^{j(i+1)+k(j+1)}\psi_{i}(\I^{\otimes\,k} \otimes \mu_{j} \otimes |\I|^{j \otimes (i-k-1)})
\end{split}
\end{equation} 
$\Psi$ is {\em strictly unital} if $\psi_{i}(\xi \otimes a_{1} \otimes \cdots \otimes a_{i-1}) = 0$ when $a_{j} = \I_{A}$ for some $j$ and $i > 1$. The {\em identity} morphism $I_{M}$ is the collection of maps $i_{1}(\xi) = \xi$, $i_{j} = 0$ for $j > 1$
\end{defn}

\noindent Likewise, if we have two morphisms of $A_{\infty}$-modules $\Phi: (M'[-1],D_{M'[-1]} \rightarrow (M''[-1],D_{M''[-1]})$ and $\Psi: (M[-1],D_{M[-1]} \rightarrow (M'[-1],D_{M'[-1]})$ over $(A,\mu)$, when we take their composition $\Phi \circ \Psi$, we can write it it in terms of the cores of $\Phi$ and $\Psi$, and then adjust the signed identities to be on $A$. This process gives 

\begin{defn}[\cite{Bor1}]
Let $\Psi$ be an $A_{\infty}$-morphism from $M$ to $M'$, and let $\Phi$ be an $A_{\infty}$-morphism from $M'$ to $M''$. The composition $\Phi \ast \Psi$ is the morphism whose component maps for $n \geq 1$ are
$$
(\Phi \ast \Psi)_{n}^{1} = \sum_{i+j=n+1}(-1)^{(i+1)(j+1)}\phi_{i}(\psi_{j} \otimes |\I|^{(j+1) \otimes (i-1)})
$$
\end{defn}

\noindent This is almost the composition defined in $\mathcal{C}^{\ast}_{A[-1]}$, but in transferring to $A$, we use
$|\I|_{A[-1]}^{(j+1) \otimes (i-1)} = (-1)^{(j+1)(i-1)}|\I|_{A}^{(j+1) \otimes (i-1)}$. This accounts for the additional
sign. 

\begin{defn}
Two $A_{\infty}$ morphisms $\Psi, \Phi$ from $M$ to $M'$ over $(A,\mu)$ are homotopic if they are homotopic as $\infty$-morphisms from $(M[-1],D_{M[-1]})$ to $(M'[-1], D_{M'[-1]})$ over $(A[-1],D_{\mu})$. 
\end{defn}

\noindent If $H$ is a homotopy, it is order $1$. Writing out the conditions in terms of its core, using the commutator proposition, and adjusting the signs in using $|\I|_{A}$ instead of $|\I|_{A[-1]}$ produces an equivalent definition. 

\begin{defn}[\cite{Bor1}]
Let $\Psi, \Phi$ be $A_{\infty}$-morphisms from $M$ to $M'$. $\Psi$ and $\Phi$ are homotopic if there is a set of maps
$\{h_{i}\}$ with $h_{i} : M \otimes A^{\otimes\,(i-1)} \longrightarrow M'[i]$ such that
\begin{equation}
\begin{split}
\psi_{i} - \phi_{i} = \sum_{i+j=n+1}(-1)^{(i+1)j}m'_{i}&(h_{j}\otimes |\I|^{j \otimes (i-1)})\\
& + \sum_{i+j=n+1} (-1)^{(i+1)j}h_{i}(m_{j} \otimes |\I|^{j \otimes (i-1)}) \\
& \hspace{0.5in} + \sum_{i+j=n+1,k>0}(-1)^{k(j+1) + j(i+1)}h_{i}(\I^{\otimes\,k} \otimes \mu_{j} \otimes |\I|^{j \otimes (i-k-1)})
\end{split}
\end{equation}
and for $i > 1$, $h_{i}(\xi \otimes a_{1} \otimes \cdots \otimes a_{i-1}) = 0$ when $a_{j} = \I_{A}$ for some $j$.
\end{defn} 

\noindent In short, all the notions of an $A_{\infty}$-object, $O$ come from the same notion for an $\infty$-object applied to $O[-1]$, and then adjusting the signs on $|\I|_{O[-1]}$ to get a formula without grading shifts.  

\subsection{Incorporating a factor on the right}

\noindent Let $A$ and $N$ be $\Z$-graded $R$-modules (as above). We can lift maps $\mathcal{T}^{\ast}(A) \rightarrow \mathcal{T}^{\ast}(A)$ to maps which take account of $N$: 

\begin{defn}
Let $\Psi : \mathcal{T}^{\ast}(A) \rightarrow \mathcal{T}^{\ast}(A)$ be an order $r$ module map. $\Psi_{N}$ is the map
$\mathcal{T}^{\ast}(A) \otimes N \rightarrow \mathcal{T}^{\ast}(A) \otimes N$ with component maps $A^{\otimes n} \otimes N \rightarrow A^{\otimes m} \otimes N$ given by
$$
\big(\Psi_{N}\big)_{n,m}  = \Psi_{n,m} \otimes |\I|_{N}^{n-m+r}
$$
\end{defn}

%\begin{defn}
% is the map whose components are given by
%$$
%\big(D_{\mu,N}\big)_{n,m}  = \big(D_{\mu}\big)_{n,m} \otimes |\I|_{N}^{n-m+1}
%$$
%for $n,m \geq 1$, and $0$ otherwise.   
%\end{defn}

\noindent We can also extend maps with domain $N$:

\begin{defn}
Let $\phi: N \rightarrow \mathcal{T}^{\ast}(A) \otimes N'$ be a degree $r$ map with projections $\phi_{i} : N \rightarrow A^{\otimes i} \otimes N'$. The extension of $\phi$ is the degree $r$ map $\overline{\phi}: \mathcal{T}^{\ast}(A) \otimes N \rightarrow \mathcal{T}^{\ast}(A) \otimes N'$ with component $\overline{\phi}_{nm} : A^{\otimes n}\otimes N \rightarrow A^{\otimes m}\otimes N'$ given by
$$
\overline{\phi}_{nm} = (-1)^{nr}\big(\I_{A}^{n} \otimes \phi_{m-n}\big)
$$
for $n \leq m$ and $0$ otherwise. 
\end{defn}

\begin{prop}
$\overline{\phi}$ is the extension of $\phi$ if and only if $\overline{\phi} = \phi \oplus (-1)^{r}(\I_{A} \otimes \overline{\phi})$ under the isomorphism $\mathcal{T}^{\ast}(A) \otimes N \cong N \oplus A \otimes \mathcal{T}^{\ast}(A) \otimes N$. 
\end{prop}

\noindent {\bf Proof:} For $m \geq n > 0$ we have $(\I_{A} \otimes \overline{\phi})_{nm} = (\I_{A} \otimes \overline{\phi}_{n-1,m-1})$ $=(-1)^{(n-1)r}(\I_{A} \otimes \I_{A}^{n-1} \otimes \phi_{m-n})$ $=(-1)^{r} \overline{\phi}_{nm}$.  If $n = 0$ then $(\I_{A} \otimes \overline{\phi})_{nm} = 0$ but $\overline{\phi}_{0m} = \phi_{m}$. $\Diamond$.  \\
\ \\
\noindent {\bf Examples:}\\
\noindent (1) We think of $\psi : N \rightarrow A \otimes N$ as a map $N \rightarrow \mathcal{T}^{\ast}(A) \otimes N$ by setting $\psi_{i} = 0$ except for $\psi_{1} = \psi$. In this case, $\overline{\psi}$ only has non-zero entries $\overline{\psi}_{n,n+1} = (-1)^{nr} \big(\I_{A}^{\otimes n}\otimes \psi\big)$.\\
\noindent (2) $\I_{N} : N \rightarrow N$ can be considered as a degree $0$ map $\iota: N \rightarrow \mathcal{T}^{\ast}(A) \otimes N$ by setting $\iota_{i} = 0$ except for $\iota_{0} = \I_{N}$. In this case, $\overline{\iota}_{nn} = \I_{A}^{n} \otimes \I_{N}$ while $\overline{\iota}_{nm} = 0$ for $n \neq m$. Thus $\overline{\iota} = \I_{\mathcal{T}^{\ast}(A)\!\otimes\!N}$.\\
\ \\

\noindent We now fix an $\infty$-structure $D_{\mu}$ on $A$. Let $\mu$ be the corresponding map on $\mathcal{T}^{\ast}(A)$ with core maps $\mu^{\ast}=\{\mu_{i}\}$ and extension $\mu_{N}: \mathcal{T}^{\ast}(A) \otimes N \rightarrow \mathcal{T}^{\ast}(A) \otimes N$. The ``core'' of  $\mu_N$, is the map $\mu^{\ast}_{N}:\mathcal{T}^{\ast}(A) \otimes N \rightarrow A \otimes N$ found by extending each $\mu_{i}$ to $A^{\otimes n}\otimes N \rightarrow A \otimes N$:
$$
\bigoplus_{n=1}^{\infty} (\mu_{n} \otimes |\I_{N}|^{n})
$$ 
A similar set of identities obtain for these maps, and the extension of $D_{\mu}$.

\begin{prop}
Let $D_{\mu,N} : \mathcal{T}^{\ast}(A) \otimes N \rightarrow \mathcal{T}^{\ast}(A) \otimes N$ be the extension of $D_{\mu}$. Then  $D_{\mu,N} + \I \otimes D_{\mu,N} = \mu_{N}$, and $D_{\mu,N} \circ D_{\mu,N} = 0$
\end{prop}

\noindent{\bf Proof:} We know that $(D_{\mu})_{n,m} = \mu_{n,m} - (\I \otimes D_{\mu})_{n,m}$. On the other hand, $(\I \otimes D_{\mu})_{n,m} = (\I \otimes (D_{\mu})_{n-1,m-1})$. Thus $(D_{\mu})_{n,m} \otimes |\I|_{N}^{n-m+1} = \mu_{n,m} \otimes |\I|_{N}^{n-m+1} - (\I \otimes (D_{\mu})_{n-1,m-1} \otimes |\I|_{N}^{n-m+1})$. Consequently, $\big(D_{\mu,N}\big)_{n,m} = \mu_{n,m} \otimes |\I|_{N}^{n-m+1} - \big(\I \otimes (D_{\mu,N})_{n-1,m-1}\big)$. Thus $D_{\mu,N} + \I \otimes D_{\mu,N} = \mu_{N}$. $\Diamond$\\
\ \\
\noindent As a consequence of the proposition,
$$
D_{\mu,N} = \sum_{l=0}^{\infty}(-1)^{l}(\I^{\otimes l} \otimes \mu_{N})
$$  

\subsection{Type D-structures}

\begin{defn}
A type $D$-structure on $N$ over $(A,D_{\mu})$ is an order $0$ map $\Delta:  N \rightarrow \mathcal{T}^{\ast}(A) \otimes N$ such that
\begin{enumerate}
\item $\Delta_{0} = \I_{N}$
\item $(\I_{A}^{\otimes n} \otimes \Delta_{m})\Delta_{n} = \Delta_{m+n}$, and
\item $D_{\mu,N} \circ \Delta = 0$
\end{enumerate}
\end{defn}

\begin{defn}
A type $D$-structure $\Delta$ on $N$ (over $(A,D_{\mu})$) is {\em bounded} if there is an $N \in \N$ such that $\Delta_{n} \equiv 0$ whenever $n \geq N$.  
\end{defn}

\noindent From now on we will all type $D$ structures in this paper will be bounded, unless otherwise stated.\\
\ \\
\noindent If we let $\delta = \Delta_{1}$ then $(\I_{A}^{\otimes n} \otimes \Delta_{m})\Delta_{n} = \Delta_{m+n}$ implies that
\begin{equation}
\begin{split}
\Delta_{0} &= \I_{N} \\
\Delta_{1} &= \delta \\
\Delta_{n} &= (\I^{\otimes (n-1)} \otimes \delta)\Delta_{n-1}
\end{split}
\end{equation}
We will denote type $D$ structures by this core map: $(N, \delta)$ where $\delta: N \rightarrow A \otimes N$. 
%If $N$ is supported in finitely many gradings then any type $D$ structure from will be bounded.\\
\ \\
\noindent Note that we may also extend $\delta: N \rightarrow A \otimes N$ as the map $\overline{\delta}: \mathcal{T}^{\ast}(A) \otimes N \rightarrow \mathcal{T}^{\ast}(A) \otimes N$ with $\overline{\delta}_{n,n+1} = \I_{A}^{\otimes n}\otimes \delta$, and that we can similarly extend $\Delta$. 

\begin{prop}
Let $\Delta$ be the map for $\delta$, then $\Delta$ satisfies the following identities:\\
\noindent (1) $(\overline{\Delta}\Delta)_{n} = (n+1)\Delta_{n}$,\\
\noindent (2) $\Delta = \I_{N} \oplus (\I_{A} \otimes \Delta)\delta$, and\\ 
\noindent (3) $\Delta = \I_{N} \oplus \overline{\delta}\Delta $
%Thus, $\overline{\Delta} = \I_{\mathcal{T}^{\ast}(A)\!\otimes\!N} + \sum_{i > 0} \overline{\delta}^{\circ i}$.
\end{prop}

\noindent{\bf Proof:} Item (1) follows from noting that $(\I^{k}\otimes \Delta_{l})$ is $\overline{\Delta}_{k,k+l}$ and $(\I^{k} \otimes \Delta_{l})\Delta_{k} = \Delta_{l+k}$. In the composition, $l$ and $k$ are independent, so we obtain $\Delta_{n}$ in each of the $(n+1)$ ways we can write $n+1=l+k$ with $l,k \geq 0$. For item (2), note that $(\I_{A} \otimes \Delta)_{n}\delta = (\I_{A} \otimes \Delta_{n-1})\delta$ $= (\I^{\otimes (n-1)} \otimes \delta)(\I \otimes \Delta_{n-2})\delta$ $= (\I^{\otimes (n-1)} \otimes \delta) (\I_{A} \otimes \Delta)_{n-1}\delta$. Thus the components of $(\I_{A} \otimes \Delta)\delta$ follow the same definition as $\Delta$. Furthermore, $(\I_{A} \otimes \Delta)_{1}\delta = (\I_{A} \otimes \I_{N})\delta=\delta$, but $(\I_{A} \otimes \Delta)_{0} = 0$ since there must be at least one $A$ factor. Item (2) follows after adjusting the $0^{th}$ level to compensate. For item (3), we compute $\big(\overline{\delta}\Delta\big)_{0m} = \overline{\delta}_{m-1,m} \circ \Delta_{m-1}$ for $m \geq 1$.  Thus, this component equals $(\I^{\otimes (m\!-\!1)}\otimes \delta)\Delta_{m-1} = \Delta_{m}$ for $m \geq 1$. However, $\Delta_{00} = \I_{N}$.$\Diamond$\\
\ \\
\noindent The definition above uses the map $D_{\mu,N}$, but we can use the identities to replace this condition with one depending solely on the core map $\mu^{\ast}_{N}$.

\begin{prop}
$(N, \delta)$ being a type $D$ structure for $(A,D_{\mu})$ is equivalent to either
$\mu_{N}\Delta = 0$, or
$$
\mu^{\ast}_{N}\Delta = \sum_{n=1}^{\infty} (\mu_{n} \otimes |\I_{N}|^{n})\Delta_{n} \equiv 0
$$
\end{prop}

\noindent{\bf Proof:} First, we note that $D_{\mu,N}\Delta = \mu_{N}\Delta - (\I_{A}\otimes D_{\mu,N})\Delta$. If we replace the second $\Delta$ on the right with $\Delta = \I_{N} + (\I_{A}\otimes \Delta)\delta$ we will get $D_{\mu,N}\Delta = \mu_{N}\Delta - (\I_{A}\otimes D_{\mu,N})(\I_{A}\otimes \Delta)\delta$ since $(\I_{A}\otimes \Delta)$ is zero on $N \subset \mathcal{T}^{\ast}(A) \otimes N$. However, $(\I_{A}\otimes D_{\mu,N})(\I_{A}\otimes \Delta)=(\I_{A}\otimes D_{\mu,N}\Delta))$. So
$$\mu_{N}\Delta = D_{\mu,N}\Delta + (\I_{A}\otimes D_{\mu,N}\Delta)\delta$$. \\
\ \\
\noindent Thus, when $D_{\mu,N}\Delta = 0$, then $\mu_{N}\Delta = 0$. On the other hand, if $\mu_{N}\Delta = 0$ then $D_{\mu,N}\Delta = -(\I_{A}\otimes D_{\mu,N}\Delta)\delta$. Iterating this relation, yields
$D_{\mu,N}\Delta = (\I_{A}\otimes (\I_{A}\otimes D_{\mu,N}\Delta)\delta)\delta $ $= (\I_{A}^{\otimes 2} \otimes D_{\mu,N}\Delta)(\I_{A} \otimes \delta)\delta $ $= (\I_{A}^{\otimes 2} \otimes D_{\mu,N}\Delta) \Delta_{2}$. By induction, we can show that $D_{\mu,N}\Delta = (-1)^{n} (\I_{A}^{\otimes n} \otimes D_{\mu,N}\Delta) \Delta_{n}$. Since $\delta$ is assumed to be bounded, $\Delta_{n} \equiv 0$ for $n$ large enough. Thus $D_{\mu,N}\Delta = 0$ when $\mu\Delta = 0$. \\
\ \\
\noindent To complete the argument, we show that the identity in the proposition is equivalent to $\mu_{N}\Delta = 0$ is equivalent to $\mu_{N}^{\ast}\Delta = 0$. It follows from the definition of $\mu$ that 
$$(\mu\Delta)_{0n} = \sum_{j-i=n-1}\big(\mu_{i} \otimes |\I_{A}|^{i \otimes(\!j\!-\!i\!)} \otimes |\I_{N}|^{i}\big)\Delta_{j}$$
for $n \geq 1$. But $\Delta_{j} = \overline{\delta}^{j-i}\Delta_{i}$ for $j \geq 1$. Since $\delta$ is order $0$, we then have
$$(\mu\Delta)_{0n} = \overline{\delta}^{n-1}\sum_{j-i=n-1}\big(\mu_{i} \otimes |\I_{N}|^{i}\big)\Delta_{i}$$
since the later application of $\overline{\delta}$ produces factors on the {\em right} of the tensor products of $A$. Rewriting the sum to be in terms of $i = j - n + 1$, and noting that any $j \geq n$ is possible, we get
$$
(\mu\Delta)_{0n} = \overline{\delta}^{n-1}\sum_{i=0}^{\infty}\big(\mu_{i} \otimes |\I_{N}|^{i}\big)\Delta_{i}
$$
from which the statement in the proposition follows directly. $\Diamond$\\
\ \\
\noindent We now consider maps between type $D$-structures and their compositions. 

\begin{defn}
Let $(N,\delta)$ and $(N',\delta')$ be bounded type $D$-structures for $(A,D_{\mu})$. An order $r$ type $D$ map $\psi : (N,\delta) \circ\!\!\!\rightarrow (N',\delta')$ is an order $r$ map of graded modules $\psi: N \longrightarrow A \otimes N'$.
\end{defn}

\begin{defn}
Let $\psi_{i} : (N_{i},\delta_{i}) \circ\!\!\!\rightarrow (N_{i+1},\delta_{i+1})$, $i=1, \ldots, n$ be order $r_{i}$ type $D$ maps. Define $M_{n}(\psi_{n},\ldots,\psi_{1})$ to be the order $1+\sum r_{i}$ map given by
$$
M_{n}(\psi_{n},\ldots,\psi_{1}) = \mu_{N_{n+1}}^{\ast}\overline{\Delta}_{n+1}\overline{\psi}_{n}\overline{\Delta}_{n}\cdots \overline{\Delta}_{2}\overline{\psi_{1}}\Delta_{1}
$$
\end{defn}
%(-1)^{n-1}
\noindent The basic proposition relating the $M_{n}$ compositions to the $\infty$-algebra $(A,\mu)$ is 

\begin{prop}\label{prop:workhorse}
Let $\Delta_{i}$, $i=1,\ldots,n+1$ be type $D$-structures for $(A,D_{\mu})$ and let $\psi_{i}$ be a degree $r_{i}$ maps from $(N_{i},\Delta_{i})$ to $(N_{i+1},\Delta_{i+1})$. Then
\begin{equation}
\begin{split}
D_{\mu,N_{n+1}}&\overline{\Delta}_{n+1}\overline{\psi}_{n}\overline{\Delta}_{n}\cdots \overline{\Delta}_{2}\overline{\psi_{1}}\Delta_{1} \\
&= \sum_{\footnotesize\begin{array}{c} 1\!\leq i\!\leq\!n\\ 0\!\leq\!l\!\leq\!n\!-\!i\end{array}\normalsize}\!(-1)^{\sum_{p=n-L+2}^{n}\mathrm{deg}\ \psi_{p}}\overline{\Delta}_{n\!+\!1}\overline{\psi_{n}}\overline{\Delta}_{n}\cdots \overline{\Delta}_{n\!-\!l\!+\!1}\overline{M_{i}(\psi_{n\!-\!l},\ldots, \psi_{n\!-\!l\!-\!i\!+\!1})} \overline{\Delta}_{n\!-\!l\!-\!i\!+\!1}\cdots \overline{\Delta}_{2}\overline{\psi_{1}}\Delta_{1}
\end{split}
\end{equation}
\end{prop}
\noindent{\bf Proof:} We consider the image of $\xi \in N_{1}$ under the map found by alternating the $\psi_{i}$ and the $\Delta_{i}$: 
$$
\overline{\Delta}_{n+1}\overline{\psi}_{n}\overline{\Delta}_{n}\cdots \overline{\Delta}_{2}\overline{\psi_{1}}\Delta_{1}
$$
The image of $\xi$ is then a sum of terms of the following form
\begin{equation}
\epsilon_{k_{1}\!,\!k_{2}\!,\!\ldots\!,\!k_{n+1}} a^{1}_{1}\!\otimes\!\cdots\!\otimes\!a^{1}_{k_{1}}\!\otimes\!\gamma_{1}\!\otimes\!a^{2}_{1}\!\otimes\!\cdots\!a^{n}_{k_{n}}\!\otimes\!\gamma_{n}\!\otimes\!a_{1}^{n+1}\!\otimes\!\cdots\!a^{n+1}_{k_{n+1}}\otimes \xi'
\end{equation}
where 1) $a^{i}_{j} \in A$, 2) each $\gamma_{i} \in A$ marks the factor coming from a $\psi_{i}$, and 3) $\xi'$ is some element of $N_{n+1}$. The sign in front equals 
$$
\epsilon_{k_{1}\!,\!k_{2}\!,\!\ldots\!,\!k_{n}} = (-1)^{k_{1}r_{1} + (k_{1} + 1 + k_{2})r_{2} + \ldots + (k_{1} + 1 + \cdots + 1 + k_{n})r_{n}}
$$
These signs come from the signs in $\overline{\psi}_{i}$ for each of $i=1,\ldots, n$: $(k_{1} + 1 + \cdots + 1 + k_{i})r_{i}$, comes from the number of factors preceding $\psi_{i}$, including the $i-1$ factors arising from $\psi_{j}$ with $j < i$, times the degree of $\psi_{i}$. \\
\ \\
\noindent To this we will apply the map:
$$
D_{\mu,N_{n+1}}=\sum_{\footnotesize\begin{array}{c} 1\!\leq i\!\leq\!M\\ 0\!\leq\!l\!\leq\!M\!-\!i\end{array}\normalsize} (-1)^{l}\big(\I^{\otimes l} \otimes \mu_{i} \otimes |\I|^{i \otimes (M\!-\!l\!-\!i)} \otimes |\I_{N_{n+1}}|^{i})= 0  
$$
with $M = k_{1} + 1 + \cdots + 1 + k_{n} + 1 + k_{n+1}$. To simplify the computation let $L$ be the number of $\gamma_{k}$ factors which are after the closing parenthesis for $\mu_{i}$ and $I$ be the number of such factors inside $\mu_{i}$. Finally, let $\chi_{s_{1},\ldots,s_{n}} = {s_{1} + 1 + s_{2} + \ldots + 1 + s_{n}}$. We will fix the value of $I$ for a minute. After applying $D_{\mu,N_{n+1}}$ we  obtain terms of the form

\begin{equation}
\begin{split}
\epsilon_{k_{1}\!,\!k_{2}\!,\!\ldots\!,\!k_{n}}&\cdot(-1)^{\chi_{k_{1}\!,\!\ldots\!,\!k_{n-L-I}\!,\!s}}a^{1}_{1}\!\otimes\!\cdots\!\otimes\!a^{n-L-I+1}_{s}\!\otimes \\ 
&\mu_{i}\big(a^{n-L-I+1}_{s+1}\!\otimes\!\cdots\!\otimes\!\gamma_{n-L-I+1}\!\otimes\!\cdots\!\otimes\!\gamma_{n-L}\!\otimes\!a^{n-L+1}_{1}\!\otimes\!\cdots\!\otimes\!a^{n-L+1}_{s'}\big)\!\otimes\\ 
&\hspace{.5in}|a^{n-L+1}_{s'}|^{i}\!\otimes\!\cdots\!\otimes\!|\gamma_{n}|^{i}\!\otimes\!|a_{1}^{n+1}|^{i}\!\otimes\!\cdots\!|a^{n+1}_{k_{n+1}}|^{i}\otimes |\xi'|^{i}
\end{split}
\label{eqn:elt}
\end{equation}
where the additional sign comes from $(-1)^{l}$ in the definition of $D_{\mu,N}$.\\
\ \\
\noindent We now do some sign accounting. First, $\epsilon_{k_{1}\!,\!k_{2}\!,\!\ldots\!,\!k_{n}} = \prod_{t=1}^{n+1} (-1)^{\chi_{k_{1},\ldots, k_{t}}\cdot r_{t}}$. Consequently, we can use this sign to replace each $\gamma_{u}$ with $(-1)^{\chi_{k_{1},\ldots,k_{u}}r_{u}}\gamma_{u}$. Note that this is the sign which would be used in an application of $\overline{\psi_{u}}$ in the product above. For $\gamma_{n-L},\ldots, \gamma_{n-L-I+1}$, however, we rewrite $\chi_{k_{1},\ldots,k_{u}}r_{u}$ as
$\chi_{k_{1}\!,\!\ldots\!,\!k_{n-L-I}\!,\!s}r_{u} + p_{u}r_{u}$. Then $p_{u}$ is the number of factors inside $\mu_{i}$ which precede $\gamma_{u}$. We can then bring the $\chi$-sign from the front into the factors, and rewrite the portion which uses $\mu_{i}$ as
$$(-1)^{(1+ \sum_{s=1}^{I} r_{n-L-I+s})\chi_{k_{1}\!,\!\ldots\!,\!k_{n-L-I}\!,\!s}}\mu_{i}\big(a^{n-L-I+1}_{s+1}\!\otimes\!\cdots\!\otimes\!(-1)^{p_{n-L-I+1}r_{n-L-I+1}}\gamma_{n-L-I+1}\!\otimes\!\cdots\!\otimes\!(-1)^{p_{n-L}r_{n-L}}\gamma_{n-L}\!\otimes\!a^{n-L+1}_{1}\!\otimes\!\cdots\!\otimes\!a^{n-L+1}_{s'}\big)
$$
The sign in front is the same as the sign introduced in extending to get $\overline{M_{I}}(\psi_{n-L}, \ldots, \psi_{n-L-I+1})$, a degree $1+ \sum_{s=1}^{I} r_{n-L-I+s}$ map, after skipping $\chi_{k_{1}\!,\!\ldots\!,\!k_{n-L-I}\!,\!s}$ preceding $A$-factors.  Each $p_{u}$ is the number of factors preceding the application of $\psi_{u}$ in $M_{I}(\psi_{n-L-I+1}, \ldots, \psi_{n-L})$ before extending. There is another sign which is also added when we change to $\overline{M}_{I}$: in $\overline{M}_{I}$ we use $\mu^{\ast}_{N}$ not $\mu^{\ast}$. The action on $N$ factor introduces another sign: that in $|\xi'|^{i}$ versus $\xi'$ where $\xi'$ is the term in the $N$ factor coming right after the application of $\mu_{i}$. This sign is $(-1)^{i\mathrm{deg}(\tilde{\xi})}$\\ %To get $M_{I}$ we still need a $(-1)^{I-1}$ as well. \\
\ \\
\noindent Last we consider the terms on the third line. We note that the sign introduced is $-1$ raised to $i(\sum \mathrm{deg} a^{l}_{j} + \sum \mathrm{deg}\gamma_{t} + \mathrm{deg}(\xi'))$ from the signed identity terms, times $(-1)$ raised to the sum of $\chi_{k_{1}\!,\!\ldots\!,\!k_{p}}r_{p}$ for $n \geq p \geq n - L + 1$. In the $\infty$-relation for $\psi_{j}$ we apply $\overline{\psi}_{r}$ after we have used $\mu_{i}$ to contract $i$ factors to $1$ factor. Thus the exponent we need differs from $\chi_{k_{1}\!,\!\ldots\!,\!k_{p}}r_{p}$ by $(i-1)r_{p}$. This occurs for each of the $n-(n-L + 1) + 1 = L$ factors after $\mu_{i}$. Thus, the sign is different by $(-1)^{\sum_{p=n-L+2}^{n} (i-1)r_{p}}$. In addition, $\sum \mathrm{deg} a^{l}_{j} + \sum \mathrm{deg}\gamma_{t} + \mathrm{deg}(\xi')$ is $\mathrm{deg}(\tilde{\xi}) + \sum_{p=n-L+2}^{n} r_{p}$ where $\tilde{\xi}$ is the result in the $N$ factor immediately after applying $\overline{M}_{i}$. This introduces another $i\sum_{p=n-L+2}^{n} r_{p}$ in the exponent. Consequently, the sign remaining after combining is $(-1)^{i\mathrm{deg}(\tilde{\xi}) + \sum_{p=n-L+2}^{n}r_{p}}$. Combining with the sign above, we are left with $(-1)^{\sum_{p=n-L+2}^{n}r_{p}}$.   \\
\ \\
\noindent Thus  $(-1)^{l}\big(\I^{\otimes l} \otimes \mu_{j} \otimes |\I|^{j \otimes (i-l)}\big)$ applied to each term in 
$$
\big(\overline{\Delta}_{n+1}\overline{\psi}_{n}\overline{\Delta}_{n}\cdots \overline{\Delta}_{2}\overline{\psi_{1}}\Delta_{1}\big)(\xi)
$$
is the same as a term in 
$$
(-1)^{\sum_{p=n-L+2}^{n}r_{p}}
\big(\overline{\Delta}_{n+1}\overline{\psi}_{n}\overline{\Delta}_{n}\cdots\overline{\Delta}_{n-L+2}\overline{M_{J}}(\psi_{n-L+1},\ldots, \psi_{n-L-J+2})\overline{\Delta}_{n-L-J+2}\cdots \overline{\Delta}_{2}\overline{\psi_{1}}\Delta_{1}\big)(\xi)
$$
If we add over all the terms we obtain the desired identity.\\
\ \\
\noindent Note: We are interpreting $I=0$ as the case where $\mu_{i}$ is applied solely to $A$-factors which come from $\overline{\Delta}$'s. In the final summation these will all cancel since $\Delta_{i}$ is a type $D$-structure. $\Diamond$ \\
\ \\
\noindent We note that for type $D$ morphisms, $\sum_{p=n-L+2}^{n}r_{p} = L-1$ and the signs will mimic the $\infty$-relations used above. 

\begin{prop}
The compositions $M_{n}, n \in \N$ satisfy the following $\infty$-relations:
$$
\sum_{\footnotesize\begin{array}{c} i\!+\!j\!=\!n\!+\!1 \\ 1 \leq l \leq i \end{array}\normalsize} (-1)^{\sum_{p=n-L+2}^{n}\mathrm{deg}\ \psi_{p}}M_{i}(\psi_{n},\ldots, \psi_{n-l+2}, M_{j}(\psi_{n-l+1},\ldots, \psi_{n-l-j+2}),\psi_{n-l-j+1},\ldots,\psi_{1}) = 0  
$$ 
\end{prop}
\ \\
\noindent{\bf Proof:} We compose $\mu^{\ast}_{N_{n+1}}$ to $D_{\mu,N_{n+1}}\overline{\Delta}_{n+1}\overline{\psi}_{n}\overline{\Delta}_{n}\cdots \overline{\Delta}_{2}\overline{\psi_{1}}\Delta_{1}$. Since $(A,\mu)$ is an $\infty$-algebra we know that
$ \mu^{\ast}_{N_{n+1}}D_{\mu,N_{n+1}} = 0$. On the other hand, using \ref{prop:workhorse}, we see that this
implies
$$
0 = \mu^{\ast}_{N_{n+1}}\left(\sum_{\footnotesize\begin{array}{c} 1\!\leq i\!\leq\!n\\ 0\!\leq\!l\!\leq\!n\!-\!i\end{array}\normalsize}\!(-1)^{\sum_{p=n-L+2}^{n}\mathrm{deg}\ \psi_{p}}\overline{\Delta}_{n\!+\!1}\overline{\psi_{n}}\overline{\Delta}_{n}\cdots \overline{\Delta}_{n\!-\!l\!+\!1}\overline{M_{i}(\psi_{n\!-\!l},\ldots, \psi_{n\!-\!l\!-\!i\!+\!1})} \overline{\Delta}_{n\!-\!l\!-\!i\!+\!1}\cdots \overline{\Delta}_{2}\overline{\psi_{1}}\Delta_{1}\right)
$$
Moving $\mu^{\ast}_{N_{n+1}}$ inside the summation, and then using the definition of $M_{n}$, we obtain the $\infty$-relations we desired.$\Diamond$.\\
\ \\
\noindent We now concentrate on $M_{1}$ and $M_{2}$. Note that $M_{2}(\psi_{2},\psi_{1})$ has degree $1 + r_{1} + r_{2}$. If we limit $\psi_{i}$ to have degree $-1$, then $M_{2}(\psi_{2},\psi_{1})$ will also have degree $-1$. Thus $M_{2}$ defines a product on the degree $-1$ maps. Indeed, the $\infty$-relation on $-1$ maps has a simpler form:
$$
\sum_{\footnotesize\begin{array}{c} i\!+\!j\!=\!n\!+\!1 \\ 1 \leq l \leq i \end{array}\normalsize} (-1)^{l-1}M_{i}(\psi_{n},\ldots, \psi_{n-l+2}, M_{j}(\psi_{n-l+1},\ldots, \psi_{n-l-j+2}),\psi_{n-l-j+1},\ldots,\psi_{1}) = 0  
$$ 

\noindent From the $\infty$-relations we see that $M_{1}$ is a boundary map. We will call a map $\psi$ with $M_{1}(\psi) = 0$ a {\em closed} map. We define

\begin{defn}
A type $D$ morphism $\psi : (N,\delta) \circ\!\!\!\rightarrow (N',\delta')$ is a closed order $-1$ module map $\psi: N \longrightarrow A \otimes N'$
\end{defn}

\begin{prop}
A degree $-1$ map $\psi: N \longrightarrow A \otimes N'$  is closed if and only if
$$
D_{\mu,N'} \circ \overline{\Delta}' \circ \overline{\psi} \circ \Delta \equiv 0
$$
\end{prop}

\noindent{\bf Proof:} By the proposition \ref{prop:workhorse}, we know
$$
D_{\mu,N'} \circ \overline{\Delta}' \circ \overline{\psi} \circ \Delta = \overline{\Delta}'\overline{M_{1}(\psi)}\Delta
$$
If $M_{1}(\psi) = 0$, then $D_{\mu,N'} \circ \overline{\Delta}' \circ \overline{\psi} \circ \Delta = 0$ since $\overline{M_{1}(\psi)}= 0$. On the other hand, if the left hand side is $0$, we get that $\overline{\Delta}'\overline{M_{1}(\psi)}\Delta = 0$. This map has image in $\mathcal{T}^{\ast}(A)\otimes N'$. If we look at the image in $A \otimes N'$, we see that it equals $ \overline{\Delta}_{11}'\overline{M_{1}(\psi)}_{01}\Delta_{00} = M_{1}(\psi)$. thus, when  $D_{\mu,N'} \circ \overline{\Delta}' \circ \overline{\psi} \circ \Delta = 0$ we have $M_{1}(\psi) = 0$. $\Diamond$. \\
\ \\
We will now restrict ourselves to degree $-1$ maps of type $D$ structures. The $\infty$-identity for $n=3$ reduces to
$$
M_{1}(M_{2}(\psi_{2},\psi_{1})) + M_{2}(M_{1}(\psi_{2}),\psi_{1}) - M_{2}(\psi_{2},M_{1}(\psi_{1})) = 0
$$
We see from this identity that $M_{2}$ will take closed maps to closed maps, thereby defining a product on the kernel of $M_{1}$. Furthermore, $M_{1}$ is a (signed, right) differential for the composition $-M_{2}$. We formalize this as

\begin{prop}
If $\psi : (N,\delta) \circ\!\!\!\rightarrow (N',\delta')$ and $\phi : (N',\delta') \circ\!\!\!\rightarrow (N'',\delta'')$ are two type $D$-morphisms, then $M_{2}(\phi, \psi) : (N,\delta) \circ\!\!\!\rightarrow (N'',\delta'')$ is a type $D$ morphism. The composition $\phi \ast \psi$ is the type $D$-morphism $-M_{2}(\phi,\psi)$.
\end{prop}

\noindent We require that $A$ be (strictly) unital with identity $1_{A} \in A_{-1}$\footnote{Recall that we will let $A$ be $A'[-1]$ for some $A'$, thus $(A'[-1])_{-1} = A'_{-1+1} = A'_{0}$}. The identity is a two-sided identity for $\mu_{2}$ (which will map $A_{-1}\otimes A_{-1} \rightarrow A_{-2+1}$), but its presence as any argument in the application of another core map $\mu_{i}$ will mean the image is $0$ . 

\begin{prop}
Let $\I_{(N,\delta)}: (N,\delta) \circ\!\!\!\rightarrow (N,\delta)$ be the map $N \rightarrow A \otimes N$ defined by $x \rightarrow 1_{A} \otimes x$. Then $\I_{(N,\delta)}$ is a type $D$ morphism with $M_{2}(\I_{(N',\delta')},\psi) = \psi$ and $M_{2}(\phi,\I_{(N,\delta)})= \phi$. Furthermore, the presence of $I_{(N_{i},\delta_{i})}$ as an argument in $M_{n}$, $n \geq 3$ results in $0$. 
\end{prop} 

\noindent{\bf Proof:} First, $I_{(N,\delta)}$ has degree $-1$ since $1_{A}$ is in $A_{-1}$. Second, we show that
$I_{(N,\delta)}$ is a morphism, i.e. that is it closed for $M_{1}$: 
 
$$
\mu^{\ast}_{N}\overline{\Delta}\overline{I}_{(N,\delta)}\Delta =  0
$$

\noindent Note that the image in $A^{\otimes n} \otimes N$ will be non-zero only for $n \geq 1$. If it is non-zero, then its image will be linear combinations of terms with a $1_{A}$ in some factor of $A^{\otimes n}$, due to the presence of $I_{(N,\delta)}$. This factor will be fed into a core map $\mu_{i}$ in $\mu_{N}^{\ast}$. When $i =1$ or $i > 2$, the image will then be zero. The only potentially non-zero terms of the composition applied to $\xi$ are those with $\mu_{2}$. If $\delta(\xi) = \sum c_{i}\otimes x_{i}$ we have
\begin{equation}
\begin{split}
(\mu_{2} \otimes I_{N})(-\I_{A} &\otimes I_{(N,\delta)})\delta(\xi) +  (\mu_{2} \otimes I_{N})(\I_{A}\otimes \delta)I_{(N,\delta)}(\xi)\\
&= (\mu_{2} \otimes I_{N})(-\I_{A} \otimes I_{(N,\delta)})(\sum c_{i}\otimes x_{i}) +  (\mu_{2} \otimes I_{N})(\I_{A}\otimes \delta)(1_{A} \otimes \xi)\\
&= (\mu_{2} \otimes I_{N})(-\sum c_{i} \otimes 1_{A} \otimes x_{i}) +  \mu_{2}\big(1_{A} \otimes (\sum c_{i}\otimes x_{i}\big)\\
&= (\mu_{2} \otimes I_{N})\big(\sum (1_{A} \otimes c_{i} - c_{i} \otimes 1_{A}) \otimes x_{i}\big)\\
&= 0\\
\end{split}
\end{equation}

\noindent To verify that $I_{(N,\delta)}$ composes as the identity on both sides, we compute 
$$
\mu_{N'}^{\ast}\circ \overline{\Delta}'\circ \overline{\I_{(N,\delta)}} \circ \overline{\Delta}' \circ \overline{\psi} \circ \Delta
$$
AS above, the strict unitality of the maps $\mu_{i}$ mean that the only terms in this composition which are non-zero will be those which feed two factors of $A$ into $\mu_{N'}$. These must come from $\overline{I_{(N,\delta)}}$ and $\overline{\psi}$. Thus the entire composition collapses to a sum of terms $(\mu_{2} \otimes \I_{N})(-c_{i} \otimes 1_{A} \otimes x'_{i}) = -\mu_{2}(c_{i},1_{A})\otimes x'_{i} =- c_{i} \otimes x'_{i}$ where $\psi(\xi) = \sum c_{i} \otimes x'_{i}$. So $\psi \ast I_{(N,\delta)} = - M_{2}(I_{(N,\delta),\psi}) = \psi$. A similar argument shows that $\I_{(N,\delta)}$ acts as an identity on the left (recalling the order reversal in the product). \\
\ \\
\noindent To see than $\I_{(N,\delta)}$ in an argument of $M_{n}$ for $n > 2$ we note that since there are $n$ morphisms the $\mu_{i}$ map applied will have $i \geq n > 2$. Furthermore, at least one argument in that $\mu_{i}$ will come from $\I_{(N,\delta)}$ and thus be $1_{A}$. Since the $\mu_{i}$ form a strictly unital
$\infty$-algebra, this means that the result must be $0$. $\Diamond$\\
\ \\
\noindent However, the product $-M_{2}$ is not associative. Instead, again from the $\infty$-relations, $M_{2}$ satisfies the generalized associativity relation

\begin{equation}
\begin{split}
M_{2}(M_{2}(\psi_{3}&,\psi_{2}),\psi_{1}) - M_{2}(\psi_{3},M_{2}(\psi_{2},\psi_{1})) = \\
&-M_{1}(M_{3}(\psi_{3},\psi_{2},\psi_{1})) - M_{3}(M_{1}(\psi_{3}),\psi_{2},\psi_{1}) + M_{3}(\psi_{3},M_{1}(\psi_{2}),\psi_{1}) -  M_{3}(\psi_{3},\psi_{2},M_{1}(\psi_{1}))\\
\end{split}
\end{equation}

%\begin{defn}
%If $\psi : (N,\Delta) \circ\!\!\!\rightarrow (N',\Delta')$ and $\phi : (N',\Delta') \circ\!\!\!\rightarrow (N'',\Delta'')$ are two morphisms, their composition $\phi \ast \psi : (N,\Delta) \circ\!\!\!\rightarrow (N'',\Delta'')$ is the map
%\begin{equation}
%\begin{split}
%\phi \ast \psi = \mu_{N''}^{\ast}\circ \overline{\Delta}''\circ \overline{\phi} \circ \overline{\Delta}' \circ \overline{\psi} \circ \Delta
%\end{split}
%\end{equation}
%\end{defn}

\noindent We can simplify this relation by quotienting by the image of $M_{1}$. To this end we declare equivalent any two $-1$ morphisms $\psi$ and $\phi$  if there is a degree $-2$ map $H: N \rightarrow A \otimes N'$ with 
$$
\psi - \phi  = M_{1}(H) = \mu_{N'}^{\ast}\overline{\Delta}'\overline{H}\Delta  
$$
We call such morphisms {\em homotopic}, following the terminology in \cite{Bor1}. However, equivalent maps represent the same homology class under $M_{1}$.\\
\ \\
\noindent  That $M_{1}$ is a differential for $M_{2}$ implies that $M_{2}$ defines a composition on the equivalence classes under homotopy. The generalized associativity relation implies that the composition $M_{2}$ is associative once restricted to equivalence classes.  As usual, once we have the $\infty$ structure above, we obtain an $\infty$ structure on the homology: the set of closed morphisms after modding out by homotopy. From the arguments above, and the $\infty$-relations we obtain the following

\begin{prop}
Let $\mathcal{D}$ be the collection of $D$-structures $(N,\delta)$ over $(A,D_{\mu})$. Let $\mathrm{\textsc{Mor}}((N,\delta),(N',\delta'))$ be the homotopy equivalence classes of the set of closed degree $-1$ type $D$ maps. Then $\mathcal{D}$ with these morphism sets forms a category if we take
\begin{enumerate}
\item the composition $\mathrm{\textsc{Mor}}((N,\delta),(N',\delta')) \otimes_{R} \mathrm{\textsc{Mor}}((N',\delta'),(N'',\delta'')) \rightarrow \mathrm{\textsc{Mor}}((N,\delta),(N'',\delta''))$ to be induced from $(\psi,\psi) \rightarrow -M_{2}(\psi,\phi)$, and
\item the identity morphism at $(N,\delta)$ to be $\I_{(N,\delta)}$ 
\end{enumerate}  
\end{prop}

\subsection{For $A$ a DGA}  We are interested in the following case: $A'$ is such that $A=A'[-1]$ has an $\infty$-structure with $\mu_{i} = 0$ for $i \geq 3$. This makes $A'$ into a differential graded algebra. In this case $M_{i} \equiv 0$ for $i \geq 3$ since these require the use of $\mu_{n}$ for $n \geq i$ due to the number of $A'[-1]$ factors involved. Examining the $\infty$-relation we see that $M_{2}$ defines an associative composition on type $D$ morphisms, before quotienting by homotopy. Furthermore,  $\I_{(N,\delta)}$ is still the identity map. Thus, in this case, type $D$ structures with type $D$ morphisms form a category {\em before quotienting by the homotopy relation}. Modding out by homotopy is then quotienting this category by an ideal. \\
\ \\  
\noindent We can write out the conditions for being a type $D$ structure, a type $D$-morphism, and composition of type $D$ structure in this setting. We note that these only require grading shifts and not changes in sign. First, a type $D$ structure is an order $0$ map $\delta: N \rightarrow A'[-1] \otimes N \cong (A' \otimes N)[-1]$ satisfying 
$$
\mu^{\ast}_{N}\Delta = \sum_{n=1}^{\infty} (\mu_{n} \otimes |\I_{N}|^{n})\Delta_{n} \equiv 0
$$
Since $\mu_{n} = 0$ for $n > 2$ we can simplify this to 
$$
(\mu_{2} \otimes \I_{N})\Delta_{2} + (\mu_{1} \otimes |\I_{N}|)\Delta_{1} = 0
$$
Since $\Delta_{2} = (\I_{A'[-1]} \otimes \delta)\delta$ and $\Delta_{1} = \delta$ we obtain the relation
$$
(\mu_{2} \otimes \I_{N}) (\I_{A'} \otimes \delta)\delta + (\mu_{1} \otimes |\I_{N}|)\delta = 0
$$
By a similar argument, we see that a morphism of type $D$ structure will be an order $-1$ map $\psi: N \rightarrow (A'[-1] \otimes N')$. Thus, $\psi$ maps $N_{k}$ to $(A'\otimes N)[-1]_{k-1} \cong (A'\otimes N)_{k}$. Thus we can take a type $D$ morphism to be an order $0$ map $N \rightarrow A' \otimes N'$ which satisfies $\mu_{N'}^{\ast}\overline{\Delta'}\overline{\psi}\Delta = 0$. This simplifies to 
$$
(\mu_{2} \otimes \I_{N})(\overline{\Delta'}\overline{\psi}\Delta)_{02} + (\mu_{1} \otimes |\I_{N}|)(\overline{\Delta'}\overline{\psi}\Delta)_{01} = 0
$$
$\overline{\psi}$ will increase the number of $A'$ factors by one, so $(\overline{\Delta'}\overline{\psi}\Delta)_{01} = (\I_{A'} \otimes \I_{N'})\psi\I_{N}$
since we must use  $\Delta'_{0}$ and $\Delta_{0}$ or else have too many factors. On the other hand, in the first term we may use either $\Delta_{1}$ or $\Delta'_{1}$, but not both. Then 
$$
(\overline{\Delta'}\overline{\psi}\Delta)_{02} = (\I_{A'[-1]} \otimes \delta')\psi - (\I_{A'[-1]} \otimes \psi)\delta
$$
Under our isomorphisms, this becomes
$$
(\mu_{2} \otimes \I_{N})(\I_{A'[-1]} \otimes \delta')\psi - (\mu_{2} \otimes \I_{N})(\I_{A'[-1]} \otimes \psi)\delta + (\mu_{1} \otimes |\I_{N}|)\psi = 0
$$
The composition of two morphisms $\psi : (N,\delta) \circ\!\!\!\rightarrow (N',\delta')$ and $\phi : (N',\delta') \circ\!\!\!\rightarrow (N'',\delta'')$ 
can be computed from $-M_{2}(\phi,\psi) = -\mu_{N''}^{\ast}\overline{\Delta}''\overline{\phi}\overline{\Delta}'\overline{\psi}\Delta$. Both $\overline{\psi}$ and $\overline{\phi}$ introduce $A'[-1]$-factors. Hence, the contributions of $\overline{\Delta}''$, $\overline{\Delta}'$ and $\Delta$ must either be the identity on the respective modules, or introduce $A'[-1]$-factors which force $\mu_{N''}^{\ast}$ to evaluate to $0$ due to $\mu_{i} = 0$ for $i > 2$. Thus,
$$
-M_{2}(\phi,\psi) = -(\mu_{2}\otimes \I_{N''})(-\I_{A'} \otimes \phi)\psi = (\mu_{2}\otimes \I_{N''})(\I_{A'} \otimes \phi)\psi
$$
Furthermore, a homotopy $H: N \rightarrow (A'[-1] \otimes N')$ is a degree $-2$ map, and thus can be thought of as a map $N_{k} \rightarrow (A'\otimes N'[-1])_{k-2} \cong (A' \otimes N'[+1])_{k}$. It is thus an order $0$ map $N \rightarrow (A' \otimes N')[+1]$. Furthermore, as above, we can compute
$M_{1}(H) =  (\mu_{2} \otimes \I_{N})(\overline{\Delta'}\overline{H}\Delta)_{02} + (\mu_{1} \otimes |\I_{N}|)(\overline{\Delta'}\overline{H}\Delta)_{01}$
which simplifies using
$$
(\overline{\Delta'}\overline{H}\Delta)_{02} = (\I_{A'[-1]} \otimes \delta')H + (\I_{A'[-1]} \otimes H)\delta
$$
since $H$ has even order. Thus if $\psi$ and $\phi$ are homotopic type $D$-morphisms, with homotopy $H$, if
$$
\psi - \phi = (\mu_{2} \otimes \I_{N})(\I_{A'[-1]} \otimes \delta')H + (\mu_{2} \otimes \I_{N})(\I_{A'[-1]} \otimes H)\delta + (\mu_{1} \otimes |\I_{N}|)H
$$
or, after applying the shift isomorphisms
$$
\psi - \phi = (\mu_{2} \otimes \I_{N})(\I_{A'} \otimes \delta')H + (\mu_{2} \otimes \I_{N})(\I_{A'} \otimes H)\delta + (\mu_{1} \otimes |\I_{N}|)H
$$

\subsection{Pairing}
\noindent  Since $D_{M,N}\circ D_{M,N} = 0$ and $D_{M,N} = m_{N} - \I_{M} \otimes D_{\mu,N}$ we see that
$m_{N}m_{N} =  m_{N}\big(\I_{M} \otimes D_{\mu,N}) + \big(\I_{M} \otimes D_{\mu,N})m_{N}$. $\big(\I_{M} \otimes D_{\mu,N})m_{N}$ has image in $\oplus_{n > 0} M \otimes A^{n} \otimes N$ since there is always an $A$ factor remaining in the codomain of $D_{\mu,N}$. Thus, after restricting to have domain and codomain $M \otimes_{R} N$, $$m^{\ast}_{N}\circ m_{N} = m^{\ast}_{N}\big(\I_{M} \otimes D_{\mu,N})$$ This works for any right $\infty$-module $M$.\\
\ \\
\noindent Now suppose we have type $D$ structures $(N_{i},\Delta_{i})$ for $i=1,\ldots,n+1$. Let $\psi_{i} : (N_{i},\delta_{i}) \circ\!\!\!\rightarrow (N_{i+1},\delta_{i+1})$, $i=1, \ldots, n$ be order $r_{i}$ type $D$ maps. We will now apply both sides of $m^{\ast}_{N_{n+1}}\circ m_{N_{n+1}} = m^{\ast}_{N_{n+1}}\big(\I_{M} \otimes D_{\mu,N_{n+1}})$ to 
$$
\xi=\I_{M} \otimes \big(\overline{\Delta}_{n+1}\overline{\psi}_{n}\overline{\Delta}_{n}\cdots \overline{\Delta}_{2}\overline{\psi_{1}}\Delta_{1}\big)
$$
We let $$\Omega_{n}(\psi_{n},\ldots,\psi_{1}) = m_{N_{n+1}}^{\ast}(\I_{M} \otimes\overline{\Delta}_{n+1})(\I_{M} \otimes\overline{\psi}_{n})(\I_{M} \otimes\overline{\Delta}_{n})\cdots (\I_{M} \otimes \overline{\Delta}_{2})(\I_{M} \otimes \overline{\psi_{1}})(\I_{M} \otimes\Delta_{1})$$ for $n \geq 1$, and $\Omega_{0}=m_{N_{1}}^{\ast}\Delta_{1}$.\\
\ \\
\noindent By proposition \ref{prop:workhorse}, 
\begin{equation}
\begin{split}
 m^{\ast}_{N}\big(&\I_{M} \otimes D_{\mu,N})(\xi) = \\
&\sum_{\footnotesize\begin{array}{c} i\!+\!j\!=\!n\!+\!1 \\ 1 \leq l \leq i \end{array}\normalsize} (-1)^{\sum_{p=n-l+2}^{n}\mathrm{deg}\ \psi_{p}}\Omega_{i}(\psi_{n},\ldots, \psi_{n-l+2}, M_{j}(\psi_{n-l+1},\ldots, \psi_{n-l-j+2}),\psi_{n-l-j+1},\ldots,\psi_{1})
\end{split}
\end{equation}

\noindent Similar to the proof of proposition \ref{prop:workhorse} we can analyze $(m^{\ast}_{N}\circ m_{N})(\xi)$. There are two differences between this argument and that in the proof of \ref{prop:workhorse}. The first occurs in the signs: there we removed $j$ factors of $A$ by applying $\mu_{j}$ and replaced it with an new factor (the image) which resulted in a difference of sign of $(j-1)r_{k}$ for each $\psi_{k}$ occurring after the application of $\mu_{j}$. Here, however, applying $m_{j}$ removes $j-1$ factors and merges them into the $M$ factor out front. This also results in a change of $(j-1)\,r_{k}$. The second difference is that $m^{\ast}_{N}\Delta=\Omega_{0} = \partial^{\boxtimes}$ and not zero as before. Furthermore, when we apply $m_{j}$ followed by $m_{i}$ we obtain a composition of $\Omega_{i}$ and $\Omega_{j}$ since each has image in $M \otimes N$. Putting these observations together with
the proof of proposition \ref{prop:workhorse} we get 
$$(m^{\ast}_{N}\circ m_{N})(\xi) = \sum_{\footnotesize\begin{array}{c} i\!+\!j\!=\!n \end{array}\normalsize} (-1)^{\sum_{p=j+1}^{n}\mathrm{deg}\ \psi_{p}}\Omega_{i}(\psi_{n},\ldots, \psi_{j+1}) \Omega_{j}(\psi_{j},\ldots,\psi_{1}) $$
where both $i$ and $j$ on the left side can equal $0$. Consequently,
\begin{equation}
\begin{split}
&\sum_{\footnotesize\begin{array}{c} i\!+\!j\!=\!n \end{array}\normalsize} (-1)^{\sum_{p=j+1}^{n}\mathrm{deg}\ \psi_{p}}\Omega_{i}(\psi_{n},\ldots, \psi_{j+1}) \Omega_{j}(\psi_{j},\ldots,\psi_{1})\\ &=\sum_{\footnotesize\begin{array}{c} i\!+\!j\!=\!n\!+\!1 \\ 1 \leq l \leq i \end{array}\normalsize} (-1)^{\sum_{p=n-l+2}^{n}\mathrm{deg}\ \psi_{p}}\Omega_{i}(\psi_{n},\ldots, \psi_{n-l+2}, M_{j}(\psi_{n-l+1},\ldots, \psi_{n-l-j+2}),\psi_{n-l-j+1},\ldots,\psi_{1})\\
\end{split}
\end{equation}

\subsection{Pairing a left $\infty$-module and a type $D$-structure}\label{sec:boxtimes}

\noindent Let $(A,D_{\mu})$ be an $\infty$-algebra with $\mu^{\ast}= \{\mu_{i}\}$, and let $M$ be a right $\infty$-module over $(A,D_{\mu})$ with differential $D_{M}$. We let $m = D_{M} + \I_{M} \otimes D_{\mu}$, and $m^{\ast}=\{m_{i}\}$ be the corresponding core maps $m_{n} : M \otimes A^{\otimes(n-1)} \rightarrow M$. In addition, we let $(N,\delta)$ be a type $D$-structure over $(A,D_{\mu})$.

\begin{defn}
Define $M \boxtimes N$ to be the graded module $M \otimes_{R} N$, and $\partial^{\boxtimes} : M \boxtimes N \rightarrow (M \boxtimes N)[-1]$ to be the map
$$
\partial^{\boxtimes} = m^{\ast}_{N}\big(\I_{M} \otimes \Delta\big) = \sum_{k=0}^{\infty}\big(m_{k+1} \otimes  |\I_{N}|^{k+1}\big) \circ \big(\I_{M} \otimes \Delta_{k}\big)
$$
\end{defn}

\begin{thm}[\cite{Bor1}]
$(M \boxtimes N, \partial^{\boxtimes})$ is a chain complex
\end{thm}

\noindent{\bf Proof:} We note that $\partial^{\boxtimes} = \Omega_{0}$ for $N_{1} = N$. Taking the relation for $i=j=0$ we obtain $\Omega_{0}\Omega_{0} = 0$, since the right hand side contains no terms. This shows that $\partial^{\boxtimes}$ is a boundary map. $\Diamond$\\
\ \\
\begin{prop}\label{prop:Asimplify}
For each type $D$-structure $(N,\delta)$ over $(A,D_{\mu})$ there is a functor $\mathcal{F}_{(N,\delta)}$ from the category of right $\infty$-modules over $(A,D_{\mu})$ to
the category of chain complexes. $\mathcal{F}_{(N,\delta)}$ is defined by
\begin{equation}
\begin{split}
&\mathcal{F}_{(N,\delta)}(M, D_{M}) = (M \boxtimes N, \partial^{\boxtimes})\\
&\ \\
&\mathcal{F}_{(N,\delta)}(\Phi) = \Phi_{N}^{\ast}(\I_{M} \otimes \Delta)
\end{split}
\end{equation}
where $\Phi\in \mathcal{C}(M,M')$ is a morphism of right $\infty$-modules over $(A, D_{\mu})$. We will denote $\mathcal{F}_{(N,\delta)}{\Phi}$ by $\Phi \boxtimes \I_{N}$. Furthermore, if $\Phi$ and $\Psi$ are homotopic then $\Phi \boxtimes \I_{N}$ is chain homotopic to $\Psi \boxtimes \I_{N}$. Thus, $\mathcal{F}_{(N,\delta)}$ induces a functor from the homotopy category of right $\infty$-modules to the homotopy category of chain complexes. 
\end{prop}

\noindent Before proving this proposition, we introduce a useful lemma:

\begin{lemma}\label{lem:interchange}
Let $\Phi \in \mathcal{C}(M,M')$ have order $r$. Then
$\big(\I_{M'} \otimes \Delta\big)\circ\big(\Phi^{\ast}_{N}\big)\circ\big(\I_{M} \otimes \Delta\big)= \Phi_{N} \circ \big(\I_{M} \otimes \Delta\big) $
\end{lemma}

\noindent{\it Proof of lemma \ref{lem:interchange}:} The image of  $\big(\I_{M'} \otimes \Delta\big)\circ\big(\Phi^{\ast}_{N}\big)\circ\big(\I_{M} \otimes \Delta\big)$ in $M \otimes A^{\otimes n}\otimes N$ has the form 
$$\sum_{l} (\I_{M'} \otimes \Delta_{n})(\Phi^{\ast}_{l+1}\otimes |\I_{N}|^{l+r})(\I_{M} \otimes \Delta_{l})$$
Since $\big(\I_{M} \otimes \Delta_{n}\big)$ does not change the $M$ factor, and $\big(\Phi^{\ast}_{l+1}\otimes |\I_{N}|^{l+r}\big)$ only affects the $M$ factor and the $l$ available $A$ factors, we can rewrite 
$$(\I_{M'} \otimes \Delta_{n})(\Phi^{\ast}_{l+1}\otimes |\I_{N}|^{l+r}) = (\Phi^{\ast}_{l+1}\otimes \I_{A}^{\otimes n} \otimes \I_{N})(\I_{M} \otimes \I_{A}^{\otimes l} \otimes \Delta_{n})(\I_{M} \otimes \I_{A}^{\otimes l} \otimes |\I_{N}|^{l+r})$$
Furthermore, as $\Delta_{n}$ preserves grading, so 
$$\Delta_{n}(|\I_{N}|^{l+r}) = (|\I_{A}|^{(l\!+\!r)\otimes n}\otimes |\I_{N}|^{l\!+\!r})\Delta_{n}$$
Therefore,
$$(\I_{M'} \otimes \Delta_{n})(\Phi^{\ast}_{l+1}\otimes |\I_{N}|^{l+r}) = (\Phi^{\ast}_{l+1}\otimes \I_{A}^{\otimes n} \otimes \I_{N})(\I_{M} \otimes \I_{A}^{\otimes l} \otimes |\I_{A}|^{(l\!+\!r)\otimes n}\otimes |\I_{N}|^{l\!+\!1}))(\I_{M} \otimes \I_{A}^{\otimes l} \otimes \Delta_{n})$$
However, $(\Phi^{\ast}_{l+1}\otimes \I_{A}^{\otimes n} \otimes \I_{N})(\I_{M} \otimes \I_{A}^{\otimes l} \otimes |\I_{A}|^{(l\!+\!r)\otimes n}\otimes |\I_{N}|^{l\!+\!r})) = \Phi_{n+l+1,n+1} \otimes |\I_{N}|^{l\!+\!1}$ which is an entry in $\Phi_{N}$. If we pre-compose with $\I_{M} \otimes \Delta_{l}$ we can then replace $(\I_{M} \otimes \I_{A}^{\otimes l} \otimes \Delta_{n})(\I_{M} \otimes \Delta_{l})$ with $\I_{M} \otimes \Delta_{l+n}$. Consequently, the image in $M \otimes A^{\otimes n}\otimes N$ is the 
sum $\sum_{l} (\Phi_{n+l+1,n+1} \otimes |\I_{N}|^{l\!+\!r})(\I_{M} \otimes \Delta_{l+n})$, which is also the entry in $\Phi_{N}(\I_{M} \times \Delta)$ which maps $M \otimes N$ to $M \otimes A^{\otimes n} \otimes N$. $\Diamond$.\\
\ \\
\noindent{\bf Proof of proposition:} Let $\Phi \in \mathcal{C}(M,M')$ be an order $0$ morphism from $M$ to $M'$ with core $\Phi^{\ast}= \{\phi_{k}\}$. We define 
$$\phi \boxtimes \I_{N} = \Phi_{N}^{\ast}(\I_{M} \otimes \Delta)$$ 
To see that this is a chain map, we compute   $(\Phi \boxtimes \I_{N})\partial^{\boxtimes} = \Phi_{N}^{\ast}(\I_{M} \otimes \Delta)(m^{\ast}_{N})(\I_{M} \otimes \Delta)$. Using the previous lemma we can simplify this to $\Phi_{N}^{\ast}(m_{N})(\I_{M} \otimes \Delta)$. Since $\Phi$ is a morphism of right $\infty$-modules, we have $\Phi_{N}D_{M,N} = D_{M',N}\Phi_{N}$ or 
$\Phi_{N}\big(m_{N}- (\I_{M} \otimes D_{\mu,N})\big) = \big(m'_{N} - (\I_{M} \otimes D_{\mu,N})\big)\Phi_{N}$. If we look at those terms with image in $M \otimes N$ we obtain $\Phi^{\ast}_{N}m_{N}- \Phi^{\ast}_{N}(\I_{M} \otimes D_{\mu,N}) = (m'_{N})^{\ast}\Phi_{N}$. Thus 
$$
(\Phi \boxtimes \I_{N})\partial^{\boxtimes}= \big(\Phi^{\ast}_{N}(\I_{M} \otimes D_{\mu,N}) + (m'_{N})^{\ast}\Phi_{N}\big)(\I_{M} \otimes \Delta) = (m'_{N})^{\ast}\Phi_{N}(\I_{M} \otimes \Delta)
$$
On the other hand, $\partial^{\boxtimes}(\Phi \boxtimes \I_{N}) = (m^{\ast}_{N})(\I_{M} \otimes \Delta)(\Phi_{N}^{\ast})(\I_{M} \otimes \Delta)$, which reduces to $(m^{\ast}_{N})\Phi_{N}(\I_{M} \otimes \Delta)$, using the lemma above. Thus, $(\Phi \boxtimes \I_{N})$ is a chain map. \\
\ \\
\noindent Let $\Psi: M \rightarrow M'$ and $\Phi: M' \rightarrow M''$ be morphisms of right $\infty$ modules. Then $(\Phi \ast \Psi) \boxtimes \I_{N}$ is the map $ (\Phi \ast \Psi)_{N}^{\ast}(\I_{M} \otimes \Delta) = \Phi^{\ast}_{N}\Psi_{N}(\I_{M} \otimes \Delta)$. Since $\Psi_{N}(\I_{M} \otimes \Delta)= (\I_{M} \otimes \Delta)\Psi^{\ast}_{N}(\I_{M} \otimes \Delta)$, we see that $(\Phi \ast \Psi) \boxtimes \I_{N} = \Phi^{\ast}_{N}(\I_{M} \otimes \Delta)\Psi^{\ast}_{N}(\I_{M} \otimes \Delta) = (\Phi \boxtimes \I_{N})(\Psi \boxtimes \I_{N})$ as required. Furthermore, $\I^{\infty}_{M} \boxtimes \I_{N} = (\I^{\infty}_{M,N})^{\ast}(\I_{M} \otimes \Delta)$. Since $(\I^{\infty}_{M,N})^{\ast}_{k}$ will be non-zero only for $k=1$, we see that the only non-zero term is $ (\I^{\infty}_{M})^{\ast}_{1} \otimes |\I_{N}|^{0})(\I_{M} \otimes \Delta_{0}) = (\I_{M} \otimes \I_{N})(\I_{M} \otimes \I_{N})= \I_{M\otimes_{R} N}$. Our map preserves the identity morphisms. \\
\ \\
\noindent Finally, we verify that the functor preserves homotopy relations. Suppose $\Phi_{N} - \Psi_{N} = D_{M',N}H_{N} + H_{N}D_{M,N}$ for some homotopy map: an order $-1$ map in $\mathcal{C}(M,M')$. Since  $D_{M',N}H_{N} + H_{N}D_{M,N} = (m'_{N} + (\I_{M'}\otimes D_{\mu,N})H_{N} + H_{N}(m_{N} + (\I_{M}\otimes D_{\mu,N})$, the only terms with image in $M \otimes N$ will be those without a $D_{\mu,N}$ term. Thus, $\Phi_{N}^{\ast} - \Psi^{\ast}_{N} = (m')^{\ast}_{N}H_{N} + H^{\ast}_{N}m_{N}$.\\
\ \\
\noindent Now, let $\mathbb{H}= H_{N}^{\ast}(\I_{M} \otimes \Delta)$. Then, using lemma \ref{lem:interchange} above,
\begin{equation}
\begin{split}
\partial^{\boxtimes}\mathbb{H}& + \mathbb{H}\partial^{\boxtimes}\\
&=(m')^{\ast}_{N})(\I_{M} \otimes \Delta)H_{N}^{\ast}(\I_{M} \otimes \Delta)+ H_{N}^{\ast}(\I_{M} \otimes \Delta)(m^{\ast}_{N})(\I_{M} \otimes \Delta)\\
&= ((m')^{\ast}_{N})H_{N}(\I_{M} \otimes \Delta)+ H_{N}^{\ast}m_{N}(\I_{M} \otimes \Delta)\\
&= \big(\Phi_{N}^{\ast} - \Psi^{\ast}_{N}\big)(\I_{M} \otimes \Delta)\\
&= (\Phi \boxtimes \I_{N}) - (\Psi \boxtimes \I_{N})\\
\end{split}
\end{equation}
$\Diamond$\\
\ \\
\noindent Consequently, homotopy equivalent $\infty$-modules will result in chain equivalences in of the chain complexes. 

%\begin{prop}
%If $\psi: (N,\delta) \rightarrow (N',\delta')$ is a morphism of left $D$-structures over $A$, then there is a chain map $\I_{M} \boxtimes \psi : M \boxtimes N \longrightarrow M \boxtimes N'$. Furthermore, if $\psi$ is homotopic to $\phi$ as type $D$-morphisms then $\I_{M} \boxtimes \psi$ is chain homtopic to $\I_{M} \boxtimes \phi$. For each $M$, this is a functor from $\mathcal{D}$ to the homotopy category of chain complexes.
%\end{prop}

\begin{prop}\label{prop:Dsimplify}
For each right $\infty$-module $(M,D_{M})$ over $(A,D_{\mu})$ there is a functor $\mathcal{G}_{(M,D_{M})}$ from the category $\mathcal{D}$ of homotopy classes of type $D$-structures over $(A,D_{\mu})$ to
the homotopy category of chain complexes. $\mathcal{G}_{(M,D_{M})}$ is defined by
\begin{equation}
\begin{split}
&\mathcal{G}_{(M,D_{M})}(N, \delta) = (M \boxtimes N, \partial^{\boxtimes})\\
&\ \\
&\mathcal{G}_{(M,D_{M})}(\psi) = [m_{N}^{\ast}(\I_{M} \otimes \overline{\Delta}')(\I_{M} \otimes \overline{\psi})(\I_{M} \otimes \Delta)]
\end{split}
\end{equation}
where $\psi$ represents a homotopy class of morphisms of type $D$ structure over $(A, D_{\mu})$, and the image $\mathcal{G}_{(M,D_{M})}(\psi)$  is the homotopy class of the chain map inside the brackets. We will denote $\mathcal{G}_{(M,D_{M})}(\psi)$ by $\I_{M} \boxtimes \psi$. 
\end{prop}

\begin{cor}
When $(A,D_{\mu})$ has $\mu_{i} = 0$ for $i \geq 2$ then the functor $\mathcal{G}_{(M,D_{M})}$ can be extended to a functor on
the category whose objects are type $D$ structures over $(A,D_{\mu})$ and whose morphisms are all the type $D$ morphisms between two type $D$ structures. Homotopic morphisms will be taken by $\mathcal{G}_{(M,D_{M})}$ to chain homotopy equivalent chain maps.
\end{cor}

\noindent{\bf Proof:} Let $\psi:N \rightarrow A \otimes N'$ be an order $-1$ morphism of type $D$ structures. We define $\I_{M} \boxtimes \psi : M \boxtimes N \longrightarrow M \boxtimes N'$ to be $\Omega_{1}(\psi)$, or 
$$\I_{M} \boxtimes \psi = m_{N}^{\ast}(\I_{M} \otimes \overline{\Delta}')(\I_{M} \otimes \overline{\psi})(\I_{M} \otimes \Delta)$$
By taking $n=1$ in the pairing relation we obtain
$$
(-1)^{1}\Omega_{1}(\psi)\Omega_{0} + (-1)^{0}\Omega_{0}\Omega_{1}(\psi) = \Omega_{1}(M_{1}(\psi))
$$
so if $\psi$ is a type $D$ morphism, we obtain $(\I_{M} \boxtimes \psi)\partial^{\boxtimes} =  \partial^{\boxtimes}(\I_{M} \boxtimes \psi)$. Thus, $(\I_{M} \boxtimes \psi)$ is a chain map. \\
\ \\
\noindent Suppose that $H$ is homotopy of type $D$-morphisms $\psi$ and $\phi$: $\psi - \phi = M_{1}(H)$. If we apply the same identity to $H$ we obtain
$$
(-1)^{2}\Omega_{1}(H)\Omega_{0} + (-1)^{0}\Omega_{0}\Omega_{1}(H) = \Omega_{1}(M_{1}(H))
$$
or
$$
\Omega_{1}(H)\partial^{\boxtimes} + \partial^{\boxtimes}\Omega_{1}(H) = \Omega_{1}(\psi-\phi)
$$
If we let $$ \I_{M} \boxtimes H = \Omega_{1}(H) = m_{N}^{\ast}(\I_{M} \otimes \overline{\Delta}')(\I_{M} \otimes \overline{H})(\I_{M} \otimes \Delta)$$ then
$$
(\I_{N} \boxtimes \psi) - (\I_{N} \boxtimes \phi) = (\I_{M} \boxtimes H)\partial^{\boxtimes} + \partial^{\boxtimes}(\I_{M} \boxtimes H)
$$
so homtopic type $D$ morphisms will be taken to chain homtopic chain maps. Thus the functor takes morphisms in the homotopy category of type $D$ morphisms to morphisms in the homotopy category of chain complexes. \\
\ \\
\noindent We have seen that the map $\psi \mapsto \I_{M} \boxtimes \psi$ takes homotopy classes to homotopy classes. We now verify that it maps the identity correctly, and preserves compositions. The image of $\I_{(N,\delta)}$ is the map $m_{N}^{\ast}(\I_{M} \otimes \overline{\Delta}')(\I_{M} \otimes \overline{\I_{(N,\delta)}})(\I_{M} \otimes \Delta)$. This introduces a $1_{A}$ into each term. Since $M$ is strictly unital, the only term remaining will be that employing $m_{2}$. We are thus able to add only one $A$-factor, so we get $m_{N,2}^{\ast}(\I_{M} \otimes \I_{(N,\delta)}) = m_{N,2}^{\ast}(\I_{M} \otimes 1_{A} \otimes \I_{N})$ or $\I_{M} \otimes \I_{N}$.\\
\ \\
\noindent To verify that the functor preserves compositions {\em in the homotopy categories}, recall that the composition of two type $D$ morphisms $\psi$ and $\phi$ is given by 
$$
M_{2}(\phi,\psi)=\mu_{N''}^{\ast}\overline{\Delta}''\overline{\phi}\overline{\Delta}'\overline{\psi}\Delta
$$
As a consequence, we may use the pairing relations, when $\phi$ and $\psi$ are type $D$ morphisms and $n=2$, to get
\begin{equation}
\begin{split}
(-1)^{0}\partial^{\boxtimes}&\Omega_{2}(\phi,\psi) + (-1)^{1}\Omega_{1}(\phi)\Omega_{1}(\psi) + (-1)^{2} \Omega_{2}(\phi,\psi)\partial^{\boxtimes} \\
&= (-1)^{0}\Omega_{1}(M_{2}(\phi,\psi)) +(-1)^{0} \Omega_{2}(M_{1}(\phi),\psi) +(-1)^{1} \Omega_{2}(\phi,M_{1}(\psi))\\
\end{split}
\end{equation}
Since $M_{1}(\phi)=M_{1}(\psi) = 0$, this identity becomes
$$
\partial^{\boxtimes}\Omega_{2}(\phi,\psi) + \Omega_{2}(\phi,\psi)\partial^{\boxtimes} = \Omega_{1}(\phi)\Omega_{1}(\psi)
+\Omega_{1}(M_{2}(\phi,\psi)) 
$$
Thus the map $\Omega_{1}(\phi)\Omega_{1}(\psi)$ is chain homotopic to $\Omega_{1}(-M_{2}(\phi,\psi))$. Consequently, after modding out by homotopies:
$$
(\I \boxtimes \phi)(\I \boxtimes \psi) \simeq (\I \boxtimes (\phi \ast \psi))
$$
and we have verified that the map preserves compositions and thus is a functor. 


\begin{thebibliography}{99}

\bibitem{APS1} M. Asaeda, J. Przytycki, A. Sikora, {\em Categorification of the Kauffman bracket skein module of I-bundles over surfaces}. Algebr. Geom. Topol. 4 (2004), 1177–1210 (electronic).

\bibitem{APS} M. Asaeda, J. Przytycki, A. Sikora, {\em Categorification of the skein module of tangles}. Primes and knots, 1–8,
Contemp. Math., 416, Amer. Math. Soc., Providence, RI, 2006.

\bibitem{Bar1}D. Bar-Natan, {\em On Khovanov's categorification of the Jones polynomial}. {\em Alg. \& Geom. Top.} 2:337--370 (2002).

\bibitem{Bar2}D. Bar-Natan, {\em Khovanov's homology for tangles and cobordisms}. {\em Geom. Topol.} 9:1443-1499 (2005).

\bibitem{BKel}B. Keller, {\em Introduction to $A$-infinity algebras and modules}. {\em Homology Homotopy Appl.} 3(1):1--35 (2001).

\bibitem{Khov}M. Khovanov, {\em A categorification of the Jones polynomial}. {\em Duke Math. J.}  101(3):359--426 (2000). 

\bibitem{Khta}M. Khovanov, {\em A functor-valued invariant of tangles}. {\em Algebr. Geom. Topol.} 2:665-741 (2002).

\bibitem{PLau}A. D. Lauda \& H. Pfeiffer,  {\em Open-closed TQFTS extend Khovanov homology from links to tangles}. {\em J. Knot Theory Ramifications} 18(1)87–150 (2009)

\bibitem{Bor1} R. Lipshitz, P. S. Ozsvath, \& D. P. Thurston. {\em Bordered Heegaard Floer homology: Invariance and pairing.} arXiv:0810.0687

\bibitem{Lee}E. S. Lee, {\em An endomorphism of the Khovanov invariant}. Adv. Math. 197(2):554-–586 (2005).

\bibitem{typeD}L. P. Roberts, {\em A type D structure in Khovanov homology}, preprint.

\bibitem{Viro}O. Viro, {\em Khovanov homology, its definition and ramifications}. {\em Fund. Math.} 184:317--342 (2004).



\end{thebibliography}
\end{document}